\documentclass{article}

\usepackage{hyperref} 
\usepackage[caption=true]{subfig}

\usepackage{multirow}
\usepackage{amsmath}
\usepackage{amssymb}
\usepackage{amsfonts}
\usepackage{dsfont}
\usepackage{bm}
\usepackage{graphicx}
\usepackage{epstopdf}
\usepackage{nicefrac}
\usepackage{url}

\usepackage{geometry}
\usepackage{marginnote}

\usepackage{algpseudocode}
\usepackage{algorithm}

\usepackage{xcolor}
\usepackage{color}
\usepackage{mathabx}

\usepackage{lipsum}
\usepackage{amsopn}
\usepackage{epstopdf}
\ifpdf
  \DeclareGraphicsExtensions{.eps,.pdf,.png,.jpg}
\else
  \DeclareGraphicsExtensions{.eps}
\fi

\graphicspath{ {figures/} }

\newcommand{\RR}{\mathbb{R}}
\newcommand{\K}{\mathcal{K}}

\newcommand{\bigO}{\mathcal{O}}

\newcommand{\D}{\mathcal{D}}
\newcommand{\mcA}{\mathcal{A}}
\newcommand{\mcR}{\mathcal{R}}
\newcommand{\mcK}{\mathcal{K}}
\newcommand{\mcS}{\mathcal{S}}

\newcommand{\bfb}{\bm b}
\newcommand{\bfc}{\bm c}
\newcommand{\bfd}{\bm d}
\newcommand{\bex}{{\bm b}^{\textrm{\footnotesize{ex}}}}

\newcommand{\bfeta}{\bm \eta}

\newcommand{\bfx}{\bm x}
\newcommand{\xex}{{\bm x}^{\textrm{\footnotesize{ex}}}}
\newcommand{\bfv}{\bm v}
\newcommand{\bfu}{\bm u}
\newcommand{\bfr}{\bm r}
\newcommand{\bfy}{\bm y}
\newcommand{\bfe}{\bm e}
\newcommand{\bfz}{\bm z}
\newcommand{\bfw}{\bm w}

\newcommand{\bfo}{\bm 0}
\newcommand{\mxA}{\bm A}
\newcommand{\mxC}{\bm C}

\newcommand{\mxZ}{\bm Z}
\newcommand{\mxU}{\bm U}
\newcommand{\mxV}{\bm V}
\newcommand{\mxSig}{\bm \Sigma}
\newcommand{\mxH}{\bm H}
\newcommand{\mxX}{\bm X}
\newcommand{\mxW}{\bm W}
\newcommand{\mxS}{\bm S}
\newcommand{\mxM}{\bm M}
\newcommand{\mxT}{\bm T}
\newcommand{\mxB}{\bm B}
\newcommand{\mxP}{\bm P}

\newcommand{\mxY}{\bm Y}
\newcommand{\mxI}{\bm I}
\newcommand{\mxWpg}{\mxW_p^\gamma}

\newcommand{\eps}{\varepsilon}

\newcommand{\whx}{\widehat{\bfx}}
\newcommand{\wbB}{\widebar{\mxB}}

\newcommand{\hl}{\widehat{\lambda}}

\newcommand{\bfzero}{\bm 0}

\DeclareMathOperator*{\argmin}{arg\,min}

\DeclareMathOperator*{\vectorize}{vec}
\DeclareMathOperator*{\diag}{diag}
\DeclareMathOperator*{\Tr}{Tr}

\newcommand{\vect}{\mbox{vec}}

\title{
{Krylov Methods for Low-Rank Regularization}}

\author{Silvia Gazzola\thanks{Department of Mathematical Sciences, University of Bath, United Kingdom. }
\and Chang Meng\thanks{Department of Mathematics, 
Emory University, Atlanta, GA, USA. }
\and James Nagy\footnotemark[3]
}

\date{}

\begin{document}

\maketitle

\begin{abstract}
This paper introduces new solvers for the computation of low-rank approximate solutions to large-scale linear problems, with a particular focus on the regularization of linear inverse problems. Although Krylov methods incorporating explicit projections onto low-rank subspaces are already used for well-posed systems that arise from discretizing stochastic or time-dependent PDEs, we are mainly concerned with algorithms that solve the so-called nuclear norm regularized problem, where a suitable nuclear norm penalization on the solution is imposed alongside a fit-to-data term expressed in the 2-norm: this has the effect of implicitly enforcing low-rank solutions. By adopting an iteratively reweighted norm approach, the nuclear norm regularized problem is reformulated as a sequence of quadratic problems, which can then be efficiently solved using Krylov methods, giving rise to an inner-outer iteration scheme. Our approach differs from the other solvers available in the literature in that: (a) Kronecker product properties are exploited to define the reweighted 2-norm penalization terms; (b) efficient preconditioned Krylov methods replace gradient (projection) methods; (c) the regularization parameter can be efficiently and adaptively set along the iterations. Furthermore, we reformulate within the framework of flexible Krylov methods both the new inner-outer methods for nuclear norm regularization and some of the existing Krylov methods incorporating low-rank projections. This results in an even more computationally efficient (but heuristic) strategy, that does not rely on an inner-outer iteration scheme. Numerical experiments show that our new solvers are competitive with other state-of-the-art solvers for low-rank problems, and deliver reconstructions of increased quality with respect to other classical Krylov methods. 
\end{abstract}
\section{Introduction} \label{sec:intro}

Consider the following linear system 
\begin{equation}\label{eq:Axeqb}
\mxA\bfx = \bfb,\quad\mbox{where}\quad\mxA\in\RR^{M\times N},\;\bfx\in\RR^N,\;\bfb= {\bex}+\bfeta\in\RR^M\,.
\end{equation}
We are mainly interested in large-scale linear systems (\ref{eq:Axeqb}) arising from inverse problems, where $\mxA$ is a discretization of the linear forward operator, $\bfx$ is a quantity of interest, and $\bfb$ is the observed perturbed data ($\bex= \mxA{\xex}$ being the ideally exact data, and $\bfeta$ being unknown Gaussian white noise). Our focus is on two-dimensional imaging problems, where the unknown vector $\bfx\in\RR^N$ is obtained by stacking the columns of an unknown true image $\mxX$ of size $n\times n$, with $n=\sqrt{N}$ (this operation and its inverse are denoted by $\bfx=\vectorize(\mxX)$ and $\mxX=\vectorize^{-1}(\bfx)$, respectively). 

Discrete inverse problems are ill-posed in nature \cite{Hansen2010} and, because of the {presence} of noise in (\ref{eq:Axeqb}), regularization needs to be applied so that the solution of (\ref{eq:Axeqb}) is a meaningful approximation to $\xex$. 
One typically achieves regularization by replacing the original problem (\ref{eq:Axeqb}) with a closely related one that is less sensitive to perturbations: effective regularization methods do so by incorporating known or desired properties of $\bfx$ into the solution process. 

In this paper we consider regularization methods that compute a low-rank approximate solution 
$\mxX=\vect^{-1}(\bfx)$ of (\ref{eq:Axeqb}): this is generally 
meaningful when the unknown $\bfx$ encodes a high-dimensional quantity and, in particular, in the case of a two-dimensional image. Indeed, two-dimensional images are often assumed to have low-rank or to be well-approximated by low-rank two-dimensional arrays (see \cite{recht2010} and the references therein). 

Numerical linear algebra solvers 
for the estimation of low-rank solutions to linear systems have been developed in the literature, mainly targeting well-posed linear discrete problems, such as those arising when considering the numerical solution of stochastic PDEs (see \cite{elman2017} and the references therein). 
In particular, the authors of \cite{elman2017} devise a restarted GMRES-like method (RS-LR-GMRES) that involves low-rank projections of the basis vectors of the solution subspace, 
as well as a low-rank projection of the current solution at the end of each cycle. Since, in general, the basic operations involved in standard GMRES (such as matrix-vector products and vector sums) increase the ranks of the computed quantities, low-rank projections are needed to assure that the computed solution is low-rank. In the framework of compressive sensing, the authors of  \cite{blanchard2015} consider a modified version of the conjugate gradient method that incorporates appropriate rank-truncation operations. All the methods mentioned so far employ, often in a heuristic way, Krylov subspace methods together with rank-reduction operations (e.g., projections onto a chosen set of low-rank matrices). Since many Krylov subspace methods are iterative regularization methods for (\ref{eq:Axeqb}), this brings us to the question of how incorporating rank-reduction operations 
would affect the solution of the discrete inverse problem (\ref{eq:Axeqb}), with a particular focus on imaging applications. 
Low-rank matrix estimation can be naturally formulated as a nonconvex optimization problem having either: (a) a least-squares data fitting term as objective function 
and a rank constraint; 
(b) the rank of $\mxX=\vect^{-1}(\bfx)$ as objective function and a constraint on the least-squares data fitting term. The last instance is commonly referred to as 
{\it affine rank minimization problem},  
and both formulations are in general NP-hard \cite{recht2010}. 
%
In this paper we consider the unconstrained and convex optimization problem
\begin{equation}\label{eq:nnr}
\min_{\bfx} \|\mxA\bfx-\bfb\|_2^2 + \lambda\|\vect^{-1}(\bfx)\|_\ast\,,
\end{equation}
where $\lambda>0$ is a regularization parameter and $\|\cdot\|_\ast$ denotes the nuclear norm of $\vect^{-1}(\bfx)=\mxX$, 
defined as the sum of the singular values of $\mxX$. Indeed, if the singular value decomposition (SVD) of $\mxX$ is given by $\mxX=\mxU_{\mxX}\mxSig_{\mxX}\mxV_{\mxX}^T$, where $\mxU_{\mxX},\mxV_{\mxX}\in\RR^{n\times n}$ are orthogonal matrices, and $\mxSig_{\mxX}\in\RR^{n\times n}$ is the diagonal matrix whose diagonal entries are $\sigma_1(\mxX)\geq\dots\geq\sigma_n(\mxX)\geq 0$, then 
\[
\|\mxX\|_*=\sum_{i=1}^n\sigma_i(\mxX)\,.
\]
Problem (\ref{eq:nnr}) is refered to as a {\it nuclear norm regularized (NNR) problem}. In particular, the nuclear norm is a convex function that has been proven to be the best convex lower approximation of the rank function over the set of matrices $\mxX$ such that $\|\mxX\|_2\leq 1$ 
(see \cite{recht2010} and the references therein). 
The nuclear norm has been used in many applications, such as low-rank matrix completion and compressed sensing; see, e.g., \cite{cai2010,goldfarb2011,keshavan2010,ma2011,recht2010}, where the constrained formulation of problem (\ref{eq:nnr}) has also been considered (note that, for a proper choice of $\lambda>0$, constrained and unconstrained formulations are equivalent; see, e.g., \cite{rockafellar1970convex}). In the framework of compressive sensing, under the assumption that the matrix $\mxA$ satisfies a certain null-space property, recovery guarantees for the affine rank minimization problem are proven in \cite{fornasier2011,mohan2012}. We also consider the following formulation
\begin{equation}\label{eq:nnrp}
\min_{\bfx} \|\mxA\bfx-\bfb\|_2^2 + \lambda\|\vect^{-1}(\bfx)\|_{\ast,p}\,,\quad\mbox{where}\quad
\|\mxX\|_{*,p}=\sum_{i=1}^n(\sigma_i(\mxX))^p,\quad 0<p\leq1\,.
\end{equation}
Problem (\ref{eq:nnrp}) is refereed to as {\it NNRp problem}, and it generalizes problem (\ref{eq:nnr}) (which is obtained taking $p=1$ in (\ref{eq:nnrp})). The constrained version of (\ref{eq:nnrp}) is already considered in \cite{mohan2012}, where the authors empirically show an improved recovery performance of the constrained formulation of problem (\ref{eq:nnrp}) with $p<1$ with respect to $p=1$. Note, however, that the choice $p<1$ in (\ref{eq:nnrp}) results in a nonconvex minimization problem. 

Many different optimization methods, such as singular value thresholding (i.e., projected gradient descent) and  continuation methods \cite{goldfarb2011}, have been proposed for the solution of problem (\ref{eq:nnr}) or its constrained counterpart. In particular, the so-called IRLS(-$p$) (i.e., iteratively reweighted least squares) family of methods has recently attracted a lot of attention \cite{fornasier2011,mohan2012,recht2010}. IRLS(-$p$) solves the affine rank minimization problem by solving a sequence of problems whose objective function only involves an iteratively updated weighted 2-norm term. The authors of \cite{lu2015} apply the IRLS(-$p$) framework to the unconstrained problem (\ref{eq:nnrp}), requiring the solution of a sequence of sub-problems
\begin{equation}\label{eq:nnrprw}
\min_{\bfx} \|\mxA\bfx-\bfb\|_2^2 + \lambda\|\mxW_k\vect^{-1}(\bfx)\|_F^2\,,
\end{equation}
where $\mxW_k$ is an appropriate weight matrix to be employed to solve the $k$th sub-problem, 
and $\|\cdot\|_F$ denotes the Frobenius norm of a matrix. A gradient (projection) algorithm is typically used to solve each 
sub-problem (\ref{eq:nnrprw}). Since an IRLS(-$p$) approach is also commonly applied to 
objective functions involving a quadratic fit-to-data term and a general $p$-norm penalization on $\bfx$, 
and since efficient strategies based on Krylov methods have been devised to solve each quadratic sub-problem in the IRLS(-$p$) sequence \cite{renaut,rodriguez2008}, this brings us to the question of how Krylov methods can be best employed to solve each problem (\ref{eq:nnrprw})  
(recall that $\|\vect^{-1}(\bfx)\|_{\ast,p}$ can be regarded as a $p$-norm of the vector whose entries are the singular values of $\mxX=\vect^{-1}(\bfx)$).

The goal of this paper is to propose new efficient Krylov methods for the estimation of low-rank solutions to (\ref{eq:Axeqb}). We will mainly consider an IRLS(-$p$) approach to problem (\ref{eq:nnrp}) (rather than incorporating low-rank projections into a linear solver for (\ref{eq:Axeqb})), the upside being that low-rank is implicitly enforced into the solution by penalizing the $p$-norm of the singular values for a suitable choice of $\lambda$. Our main contributions are the new IRN-GMRES-NNR$p$ and IRN-LSQR-NNR$p$ methods for (\ref{eq:nnrp}), where automatic strategies for choosing a suitable $\lambda$ are naturally incorporated. Here and in the following, the IRN acronym indicates an iteratively reweighted norm (rather than an iteratively reweighted least squares problem, \cite{rodriguez2008}). One of the key points in deriving the new methods is expressing in matrix form the invertible linear operator mapping $\bfx$ to the reweighted 2-norm of the singular values of $\mxX=\vect^{-1}(\bfx)$: this can be achieved in a computationally affordable way by exploiting Kronecker product properties. Each iteratively reweighted quadratic sub-problem of the form (\ref{eq:nnrprw}) 
can then be expressed as a Tikhonov regularization problem in general form, which can be straightforwardly transformed into standard form. In this way, the inverse of the linear operator mapping $\bfx$ into the reweighted 2-norm of the singular values of $\mxX=\vect^{-1}(\bfx)$ formally acts as a preconditioner for $\mxA$, and the so-called hybrid methods \cite{OS} based on the preconditioned Arnoldi (if $\mxA$ is square) or Golub-Kahan bidiagonalization algorithms can be used to efficiently approximate the solution of each problem of the form (\ref{eq:nnrprw}). 
Once a hybrid method is adopted, many automatic, adaptive, and efficient parameter choice strategies can be employed to choose a suitable $\lambda$; see \cite{KO} for an overview. Therefore, contrarily to many existing methods for (\ref{eq:nnrp}), IRN-GMRES-NNR$p$ and IRN-LSQR-NNR$p$ have the advantage of not requiring 
a regularization parameter (either $\lambda$ or the desired rank of the solution) to be available in advance of the iterations, nor the repeated solution of (\ref{eq:nnrp}) for different regularization parameters. 

Although inherently efficient, both the IRN-GMRES-NNR$p$ and IRN-LSQR-NNR$p$ methods are inner-outer iteration schemes, where each outer iteration requires running a ``preconditioned'' Krylov subspace method until convergence (inner iteration) before updating the weights (and therefore the ``preconditioner'') in the next outer iteration. In order to avoid inner-outer iterations and with the aim of generating only one approximation subspace for the solution of (\ref{eq:nnrp}), where a new ``preconditioner'' is incorporated as soon as a new approximate solution becomes available (i.e., at each iteration), we propose to solve (\ref{eq:nnrp}) using flexible Krylov subspace methods, such as those based on the flexible Arnoldi \cite{saad1993} and 
Golub-Kahan \cite{chunggazzola2018} algorithms. 
The use of flexible Krylov methods for 
$p$-norm regularization of inverse problems was already proposed in \cite{chunggazzola2018,gazzolanagy2014}; however, differently from the available solvers, our new approach involves iteratively defining both weights and transform matrices (i.e., the linear operator mapping $\vect^{-1}(\bfx)$ into its singular values). Switching from IRN-GMRES-NNR$p$ and IRN-LSQR-NNR$p$ to their flexible counterparts (dubbed FGMRES-NNR$p$ and FLSQR-NNR$p$, respectively) allows for savings in computations and, although FGMRES-NNR$p$ and FLSQR-NNR$p$ 
are purely heuristic, it leads to approximate solutions whose accuracy on many test problems is comparable to the ones of 
other well established solvers for (\ref{eq:nnr}). Motivated by the same idea of avoiding inner-outer iteration cycles while adaptively incorporating (low-rank) information into the approximation subspace for the solution, we also propose a flexible version of the projected and restarted Krylov subspace methods (such as RS-LR-GMRES, \cite{elman2017}) that were originally devised for square, 
considering also extensions to rectangular matrices $\mxA$.

This paper is organized as follows. In Section \ref{sec:bg} we review the available low-rank Krylov methods for square linear systems and, after surveying the available flexible Krylov solvers, we formulate new low-rank flexible Krylov solvers for both square and rectangular problems, where the basis vectors for the approximation subspace are truncated to low-rank. 
In Section \ref{sec:method} we derive the new iteratively reweighted methods for (\ref{eq:nnrp}) as fixed-point methods, and we describe how to efficiently solve each reweighted problem of the form (\ref{eq:nnrprw}) 
using preconditioned Krylov methods: this leads to the IRN-GMRES-NNR$p$ and IRN-LSQR-NNR$p$ methods; their flexible counterparts (FGMRES-NNR$p$ and FLSQR-NNR$p$, respectively) are also derived. Some implementation details, such as stopping criteria and regularization parameter choice strategies for the new methods, are unfolded in Section \ref{sec:implementation}. 
Numerical results are presented in Section \ref{sec:experiments}, including comparisons between the proposed methods, low-rank projection methods, projected gradient methods, and standard Krylov subspace methods. Conclusions are drawn in Section \ref{sec:conclusions}.
\medskip

\emph{Definitions and notations}. Matching lower and upper case letters are used to denote the ``vectorized'' and ``matricized'' versions of a given quantity, respectively; e.g., $\bfc=\vect(\mxC)$ and $\mxC=\vect^{-1}(\bfc)$. We denote the $i$th entry of a vector $\bfc$ by $[\bfc]_i$, and the $(i,j)$th entry of a matrix $\mxC$ by $[\mxC]_{ij}$ or, using MATLAB-like notations, $[\bfc]_i=\bfc(i)$, $[\mxC]_{ij}=\mxC(i,j)$. Using again MATLAB-like notations, $\bfd=\diag(\mxC)$ defines a vector $\bfd$ whose entries are the diagonal elements of the matrix $\mxC$. $\Tr(\mxC)$ denotes the trace of a matrix $\mxC$. $\mcR(\mxC)$ denotes the range (or column space) of the matrix $\mxC$, and 
$\mcK_m(\mxA,\bfb)$ denotes  the $m$-dimensional Krylov subspace defined by $\mxA$ and $\bfb$, i.e.,
$\mcK_m(\mxA,\bfb) = \mbox{span}\left\{\bfb,\, \mxA\bfb,\, \mxA^2\bfb, \, \ldots, \, \mxA^{m-1}\bfb\right\}$.
We denote by $\mxI\in\RR^{d\times d}$ the identity matrix of order $d$, and by $\bfe_i$ the $i$th canonical basis vector of $\RR^d$, where $d$ should be clear from the context. Note that, in the following, we will quite often interchange $\bfx$ and $\mxX$ and, with a slight abuse of notations, we will denote the action of a linear operator on $\bfx$ or $\mxX$ by $\mcA(\mxX)=\mcA\mxX=A\bfx$, and the action of the adjoint operator by $\mcA^\ast(\mxY)=\mcA^\ast\mxY=A^T\vect(\mxY)$. 

\section{Low-rank projection methods: classical and new approaches}\label{sec:bg}


As recalled in Section \ref{sec:intro}, when solving square well-posed linear systems coming from the discretization of some instances of stochastic or time-dependent PDEs, a suitable rearrangement of the solution is expected to be low-rank: for this reason, schemes that incorporate low-rank projections within the basis vectors and the approximate solution obtained by a Krylov method have been proposed in the literature. In the following we summarize the working ideas underlying the so-called restarted low-rank-projected GMRES (RS-LR-GMRES) method proposed in \cite{elman2017}. 

The starting points for the derivation of RS-LR-GMRES are the basic properties and relations underlying GMRES. Indeed, one can define GMRES for the solution of (\ref{eq:Axeqb}) with a square $\mxA\in\RR^{N\times N}$ and initial guess $\bfx_0=\bfzero$ by generating a matrix $\mxV_m=[\bfv_1,\dots,\bfv_m]\in\RR^{N\times m}$ with orthonormal columns, such that $\mcR(\mxV_m)=\mcK_m(\mxA,\bfb)$, and imposing that the residual $\bfr_m=\bfb-\mxA\bfx_m$ is orthogonal to $\mxU_m=\mxA\mxV_m$. In practice, at the $k$th iteration of GMRES, one computes
\begin{equation}\label{GMRESbasis}
\bfu_k=\mxA\bfv_{k-1}\quad\mbox{and}\quad 
\|\bfv_{k}\|_2\bfv_{k}=(\mxI - \mxV_{k-1}(\underbrace{\mxV_{k-1}^T\mxV_{k-1}}_{=\mxI})^{-1}\mxV_{k-1}^T)\bfu_{k}\,,
\end{equation}
and the approximate solution is computed as
\begin{equation}\label{GMRESvar}
\bfx_k=\mxV_k\bfy_k\,,\quad\mbox{where}\quad (\mxU_k^T\mxA\mxV_k)\bfy_k=\mxU_k^T\bfb\,.
\end{equation}
This procedure is mathematically equivalent to the somewhat more standard procedure that, at the $k$th iteration of GMRES, updates the partial Arnoldi factorization and computes the approximate solution as follows:
\begin{equation}\label{GMRESstd}
\mxA\mxV_k=\mxV_{k+1}\mxH_k\,,\quad\bfx_k=\mxV_k\bfy_k\,,\;\mbox{where}\; \bfy_k=\arg\min_{\bfy\in\RR^k}\|\mxH_k\bfy-\|\bfb\|_2\bfe_1\|_2\,.
\end{equation}
Note that, in particular, the matrix $\mxV_k$ appearing in (\ref{GMRESvar}) coincides with the matrix $\mxV_k$ appearing in (\ref{GMRESstd}). However, since matrix-vector products and vector sums of low-rank vectorized matrices increase the rank of the latter, relations (\ref{GMRESbasis}) and (\ref{GMRESvar}) obviously do not guarantee that the new basis vectors $\bfv_k$ for the solution nor the new solution $\bfx_k$ are low-rank. To force the basis vector for the solution and the approximate solution to be low-rank, a truncation operator should be incorporated into the GMRES algorithm. Given a vectorized matrix $\bfc=\vect(\mxC)$, and given a desired low-rank $\kappa$ for $\mxC$, one can define a truncation operator $\tau_{\kappa}(\bfc)$ by the following standard operations: 
%
 \begin{equation}\label{truncop}
 \left[\begin{array}{cl}
1. & \mbox{Take $\mxC=\vect^{-1}(\bfc)$;}\\
2. & \mbox{Compute the SVD of $\mxC$, $\mxC = \mxU_{\mxC}\mxSig_{\mxC}\mxV_{\mxC}^T$;}\\
3. & \mbox{Compute $\mxC_{\kappa} = \mxU_{\mxC}(:,1:\kappa)\mxSig_{\mxC}(1:\kappa, 1:\kappa)\mxV_{\mxC}(:, 1:\kappa)^T$;}\\
4. & \mbox{Take $\tau_{\kappa}(\bfc)=\vect(\mxC_{\kappa})$.}
 \end{array}
 \right.
 \end{equation}
RS-LR-GMRES is a restarted version of the standard GMRES method where the basis vectors for the solution are truncated at each inner iteration, and the solution itself is truncated at the beginning of each outer iteration. More precisely, at the $\ell$th outer iteration of RS-LR-GMRES, one takes $\bfv_1=\bfr_{\ell-1}/\|\bfr_{\ell-1}\|_2$, where $\bfr_{\ell-1}=\bfb-\mxA\bfx_{\ell-1}$, and, at the $k$th inner iteration, one computes
\begin{equation}\label{LRGMRESbasis}
\bfu_k=\mxA\bfv_{k-1}\quad\mbox{and}\quad 
\|\bfv_{k}\|_2\bfv_{k}=\tau_{\kappa}\left((\mxI - \mxV_{k-1}(\mxV_{k-1}^T\mxV_{k-1})^{-1}\mxV_{k-1}^T)\bfu_{k}\right)\,.
\end{equation}
Once $m$ inner iteartions are performed, the approximate solution at the $\ell$th outer iteration is computed as
\begin{equation}\label{LRGMRESvar}
\bfx_{\ell}=\tau_{\kappa}\left(\bfx_{\ell-1}+\mxV_m\bfy_m\right)\,,\quad\mbox{where}\quad (\mxU_m^T\mxA\mxV_m)\bfy_m=\mxU_m^T\bfr_{\ell-1}\,.
\end{equation}
%
%
%
%
%
%
The operations in (\ref{LRGMRESbasis}) and (\ref{LRGMRESvar}) heavily depend on the value $\kappa$ of the truncated rank, which eventually coincides with the rank of the approximate solution. In the framework of stochastic PDEs, a suitable estimate for $\kappa$ can be obtained by first performing coarse-grid computations (see \cite{elman2017} for details, and \cite{KreTob11,StoBre15} for similar approaches). 
Comparing (\ref{LRGMRESbasis}) and (\ref{GMRESbasis}) one can see that, as in standard GMRES, RS-LR-GMRES computes a new basis vector for the solution by applying the linear operator $\mxA$ to the previous basis vector $\bfv_{k-1}$ and orthogonalizing it against the previous basis vectors $\bfv_i$, $i=1,\dots,k-1$. However, since the basis vectors are truncated to low rank, the matrix $\mxV_k$ does not have orthonormal columns anymore, and $\mcR(\mxV_m)$ is not a Krylov subspace anymore. This remark leads us to the derivation of alternative low-rank projection solvers, which can be (re)casted into the framework of flexible Krylov methods and can work with both square and rectangular systems (\ref{eq:Axeqb}).

\medskip
{\it Low-rank flexible GMRES (LR-FGMRES) and low-rank flexible LSQR (LR-FLSQR)}. 
Flexible Krylov methods are a class of linear solvers that can handle iteration-dependent preconditioners: they were originally introduced in \cite{saad1993} for FGMRES, where a preconditioner for GMRES was allowed to change from one iteration to the next (either because at each iteration the preconditioner is implicitly defined by applying an iterative linear solver, or because the preconditioner can be updated with newly-computed information; see \cite{simoncini07} for an overview). 
%
In the framework of regularizing linear solvers, flexible Krylov methods were proposed in \cite{chunggazzola2018,gazzolanagy2014,gazzolasabate2019}, where the iteration-dependent ``preconditioner'' was associated to an iteratively reweighted norm approach to Tikhonov-like regularized problems involving penalization terms expressed in some $p$-norm, $0<p\leq 1$ (and, indeed, these ``preconditioners'' have the effect of enforcing specific regularity into the approximation subspace for the solution, rather than accelerating the convergence of the iterative solvers). Leveraging flexible Krylov subspaces in this setting comes with the upside of avoiding restarts of the iterative solver, which is the approach commonly used when adopting an iteratively reweighted norm method. When considering low-rank projections of the basis vectors within RS-LR-GMRES, we enforce the basis vectors to have low-rank, so to better reproduce available information about the solution of (\ref{eq:Axeqb}) (i.e., the solution should be low-rank). It is therefore natural to consider flexible Krylov methods that involve truncation of the basis vectors at each iteration, as a computationally cheaper alternative to RS-LR-GMRES that does not involve restarts. 

Considering first the case of a square $\mxA\in\RR^{N\times N}$, we can use the flexible Arnoldi algorithm \cite{saad1993} to naturally incorporate low-rank basis vectors for the solution of (\ref{eq:Axeqb}). 
In general, starting with $\bfx_0=\bfo$, at the $k$th iteration, FGMRES updates a partial flexible Arnoldi factorization and computes the $k$th approximate solution as follows:
\begin{equation}\label{eq:arnoldifac}
\mxA\mxZ_k = \mxV_{k+1}\mxH_k\,,\quad\bfx_k=\mxZ_k\bfy_k\,,\;\mbox{where}\; \bfy_k=\arg\min_{\bfy\in\RR^k}\|\mxH_k\bfy-\|\bfb\|_2\bfe_1\|_2\,.
\end{equation}
where $\mxV_{k+1}= \left[\bfv_1,\ \dots\ ,\bfv_{k+1} \right]\in\RR^{N\times (k+1)}$ has orthonormal columns, $\mxH_k\in\RR^{(k+1)\times k}$ is upper Hessenberg, and $\mxZ_k=[\mxP_1\bfv_1,\dots, \mxP_k\bfv_k] \in\RR^{N\times k}$ has columns that span the approximation subspace for the solution ($\mxP_i$ is an iteration-dependent preconditioner that is applied to $\bfv_i$ and, in the particular case of low-rank truncation, $\mxP_i\bfv_i=\tau_{\kappa_B}(\bfv_i)$, is the truncation operator defined in (\ref{truncop}), so that $\mbox{rank}(\vect^{-1}(\mxZ\bfe_i)) = \kappa_B,\ i=1,\dots,k$). 
%
The resulting algorithm is dubbed ``LR-FGMRES'', and it is summarized in Algorithm \ref{alg:lr-fgmres}. Note that the approximate solution computed as in (\ref{eq:arnoldifac}) is also truncated to guarantee rank $\kappa$ (in general, we assume $\kappa_B\neq\kappa$). Also LR-FGMRES is started with $\bfx_0=\bfzero$, to guarantee that the basis vectors for the solution (rather than a correction thereof) are low-rank.
\begin{algorithm}[H]
\caption{LR-FGMRES}\label{alg:lr-fgmres}
\begin{algorithmic}[1]
\State{Inputs: $\mxA$, $\bfb$, $\tau_{\kappa_B}$, $\tau_\kappa$}
\State{Take $\bfv_1 = \bfb/\|\bfb\|_2$}
\For{$i=1,2,\dots$ until a stopping criterion is satisfied}
  \State{Compute $\bfz_i = \tau_{\kappa_B}(\bfv_i)$ and $\bfw = A\bfz_i$}
  \State{Compute $h_{ji} = \bfw^T\bfv_j$ for $j = 1,\dots,i$ and set $\bfw = \bfw - \sum_{j=1}^i h_{ji}\bfv_j$}
  \State{Compute $h_{i+1,i} = \|\bfw\|_2$, and if $h_{j+1,j}\neq 0$, take $\bfv_{i+1} = \bfw/h_{i+1,i}$}
\EndFor
\State{Compute $\bfy_k = \argmin_{\bfy} \| H_k \bfy - \|\bfb\|_2\bfe_1\|_2^2$ and take $\bfx_k = \tau_\kappa(Z_k\bfy_k)$}
 \end{algorithmic}
\end{algorithm}
A few remarks are in order. Differently from the $k$th iteration in the inner cycle of the RS-LR-GMRES method (\ref{LRGMRESbasis}), the $k$th iteration of LR-FGMRES expands the approximation subspace by modifying (i.e., truncating) the previous orthonormal basis vector for the space $\mcR([\bfb, \mxA\mxZ_k])$. Analogously to RS-LR-GMRES, the basis vectors for the approximate LR-FGMRES solution are all of rank $\kappa$, are not orthogonal, and do not span a Krylov subspace. Differently from RS-LR-GMRES, the basis vector for the space $\mcR([\bfb, \mxA\mxZ_k])$ are orthogonal. Also, the $k$th LR-FGMRES approximate solution is obtained by solving an order-$k$ projected least squares problem that is formally analogous to the GMRES one (see (\ref{GMRESstd}) and (\ref{eq:arnoldifac})).  


With LR-FGMRES in place, the extension to more general matrices $\mxA\in\RR^{M\times N}$, with $M$ not necessarily equal to $N$, can be naturally devised considering the flexible Golub-Kahan (FGK) process  \cite{chunggazzola2018}.  Taking $\bfx_0=\bfo$ as initial guess, the $k$th iteration, FGK updates partial factorizations of the form
\begin{equation}\label{eq:gkfac}
\mxA\mxZ_k = \mxU_{k+1}\mxM_k\ \ \mathrm{and}\ \ \mxA^T\mxU_{k+1} = \mxV_{k+1}\mxT_{k+1},
\end{equation}
where the columns of $\mxU_{k+1}\in\RR^{M\times(k+1)}$, $\mxV_{k+1}\in\RR^{N\times(k+1)}$ are orthonormal, $\mxM_k\in\RR^{(k+1)\times k}$ is upper Hessenberg, $\mxT_{k+1}\in\RR^{(k+1)\times (k+1)}$ is upper triangular, and $\mxZ_k=[\mxP_1\bfv_1,\dots, \mxP_k\bfv_k] \in\RR^{N\times k}$ has columns that span the approximation subspace for the solution ($\mxP_i$ is an iteration-dependent preconditioner that is applied to $\bfv_i$ and, in the particular case of low-rank truncation, $\mxP_i\bfv_i=\tau_{\kappa_B}(\bfv_i)$, as defined in (\ref{truncop}), so that $\mbox{rank}(\vect^{-1}(\mxZ\bfe_i)) = \kappa_B,\ i=1,\dots,k$). 
The flexible LSQR method (FLSQR) uses the FGK process (\ref{eq:gkfac}) to generate iterates of the form $\bfx_k = Z_k\bfy_k$, where the vector $\bfy_k$ is computed as $\bfy_k = \argmin_{\bfy} \Big\Vert \mxM_k\bfy - \|\bfb\|_2\bfe_1\Big\Vert_2^2$. When rank-truncation of the basis vectors takes place at each iteration, and the final approximate solution is rank-truncated as well, the resulting algorithm is dubbed ``LR-FLSQR'', and it is summarized in Algorithm \ref{alg:lr-flsqr}. 
%
\begin{algorithm}
\caption{LR-FLSQR}\label{alg:lr-flsqr}
\begin{algorithmic}[1]
\State{Inputs: $\mxA$, $\bfb$, $\tau_{\kappa_B}$, $\tau_\kappa$}
\State{Take $\bfu_1 = \bfb/\|\bfb\|_2$}
\For{$i=1,2,\dots,$ until a stopping criterion is satisfied}
  \State{Compute $\bfw = \mxA^T\bfu_i$, $t_{ji} = \bfw^T\bfv_j$ for $j = 1,\dots,i-1$}
  \State{Set $\bfw = \bfw - \sum_{j=1}^{i-1} \ t_{ji}\bfv_j$, compute $t_{ii} = \|\bfw\|$ and take $\bfv_i = \bfw/t_{ii}$}
  \State{Compute $\bfz_i = \tau_{\kappa_B}(\bfv_i)$ and $\bfw = \mxA\bfz_i$}
  \State{Compute $m_{ji} = \bfw^T\bfu_j$ for $j = 1,\dots,i$ and set $\bfw = \bfw - \sum_{j=1}^{i} \ m_{ji}\bfu_j$}
  \State{Compute $m_{i+1,i} = \|\bfw\|$ and take $\bfu_{i+1} = \bfw/m_{i+1,i}$}
\EndFor
\State{Compute $\bfy_k = \argmin_{\bfy} \|M_k\bfy - \|\bfb\|_2\bfe_1\|_2^2$ and take $\bfx_k = \tau_\kappa(Z_k\bfy_k)$}
\end{algorithmic}
\end{algorithm}
Note that, similarly to RS-LR-GMRES, 
both LR-FGMRES and LR-FLSQR are quite heuristic. Although the low-rank projection idea can be formulated in the flexible framework, we lack a formal formulation of the problem that is being solved, and also a justification of why they work. Strategies of selecting $\kappa_B$ and $\kappa$ are not so clear either. To stabilize the behavior of LR-FGMRES as the iterations proceed, one may consider imposing additional Tikhonov regularization on the projected least-squares problem in (\ref{eq:arnoldifac}), in a hybrid fashion; the same holds for LR-FLSQR (see Sections \ref{sec:flexinnr} and \ref{sec:experiments} for more details). 



%

\section{Proposed Method} \label{sec:method}
In this section, we first derive the IRN method for the solution of the NNR$p$ problem (\ref{eq:nnrp}). The starting point for our derivations is the approximation of the nondifferentiable nuclear norm regularizer by a smooth Schatten function (similarly to what is proposed in \cite{mohan2012} for the affine rank minimization problem). The optimality conditions associated to the smoothed problem give rise to a nonlinear system of equations in $\mxX$, which is handled by a fixed-point iteration scheme. We show that each iteration amounts to the solution of a Tikhonov-regularized problem involving an iteratively reweighted 2-norm regularization term, which can be efficiently solved employing ``preconditioned'' Krylov methods. Flexible Krylov methods are introduced to approximate the solution of the IRN problem within only one adaptively defined approximation subspace for the solution, bypassing the inner-outer iteration scheme required by standard Krylov methods.  
%

\subsection{Derivation}\label{ssec:deriv}
Define the smooth Schatten-$p$ function 
%
%
as
\[
\mcS_p^{\gamma}(\mxX)=\Tr((\mxX^T\mxX + \gamma \mxI)^{p/2})\,,\quad\mbox{with $\gamma>0$}\,.
\]
Note that $\mcS_p^{\gamma}(\mxX)$ is differentiable for $p>0$ and convex for $p\geq 1$. In particular, for $p=1$ and $\gamma = 0$ (i.e., no smoothing),
\[
\mcS_1^{0}(\mxX)=\Tr((\mxX^T\mxX)^{1/2})=\|\mxX\|_\ast\,.
\]
We start by considering the following smooth approximation to (\ref{eq:nnrp}):
\begin{equation}\label{eq:nnr_smooth}
\min_{\mxX\in\RR^{n\times n}}\|\mcA(\mxX)-\mxB\|_F^2+\lambda \mcS_p^{\gamma}(\mxX)\,.
\end{equation}
The following derivations are valid for $p>0$ (and we keep them generic, being aware that $p=1$ approximates (\ref{eq:nnr})). 
The optimality conditions associated to (\ref{eq:nnr_smooth}) read
\begin{eqnarray}
0 &=& \nabla_{\mxX} \left( \|\mcA(\mxX) - \mxB\|_F^2 + \lambda \mcS_p^{\gamma}(\mxX)\right)\nonumber\\
   &=& 2\mcA^\ast(\mcA(\mxX)-\mxB)+\lambda\,p(\mxX\mxX^T+\gamma\mxI)^{p/2-1} \mxX\,,
\label{eq:nlnnr}
\end{eqnarray}
where we have used that 
\[
\nabla_{\mxX}\Tr((\mxX^T\mxX + \gamma \mxI)^{p/2})=p\mxX(\mxX^T\mxX+\gamma\mxI)^{p/2-1}=p(\mxX\mxX^T+\gamma\mxI)^{p/2-1} \mxX\,.
\]
Equivalently, the nonlinear system of equations (\ref{eq:nlnnr}) with respect to $\mxX$ can be expressed as
\begin{eqnarray*}
\mxX&=& \left(\mcA^\ast\mcA+\widehat{\lambda}(\mxX\mxX^T+\gamma\mxI)^{p/2-1}\right)^{-1}\mcA^\ast\mxB\\
&=&\left(\mcA^\ast\mcA+\widehat{\lambda}((\mxX\mxX^T+\gamma\mxI)^{p/4-1/2})^T(\mxX\mxX^T+\gamma\mxI)^{p/4-1/2}\right)^{-1}\mcA^\ast\mxB\,,\quad\mbox{with $\widehat{\lambda}=\lambda\,p/2$}\,,
\end{eqnarray*}
%
%
which is naturally associated to the following fixed-point iteration scheme
\begin{equation}\label{eq:fpnnr}
\mxX_{k+1}= \left(\mcA^\ast\mcA+\widehat{\lambda}((\mxX_k\mxX_k^T+\gamma\mxI)^{p/4-1/2})^T(\mxX_k\mxX_k^T+\gamma\mxI)^{p/4-1/2}\right)^{-1}\mcA^\ast\mxB\,,
\end{equation}
which leads to the solution of (\ref{eq:nnr_smooth}). Equivalently,
\[
\mxX_{k+1}=\arg\min_{\mxX}\left\|\left[\begin{array}{c}
\mcA\\
\sqrt{\hl}(\mxX_k\mxX_k^T+\gamma\mxI)^{p/4-1/2}
\end{array}
\right]\mxX -
\left[\begin{array}{c}
\mxB\\
\bfzero
\end{array}
\right]
\right\|^2_F\,,
\]
i.e., (\ref{eq:fpnnr}) are the normal equations associated to the penalized least squares problem written above or, equivalently,
\begin{equation}\label{eq:nnr_irn_matr}
\mxX_{k+1}=\arg\min_{\mxX}\left\|
\mcA\mxX-\mxB\right\|_F^2 + \hl\left\|(\mxX_k\mxX_k^T+\gamma\mxI)^{p/4-1/2}\mxX\right\|_F^2\,.
\end{equation}
We now reformulate problem (\ref{eq:nnr_irn_matr}) in vectorial form.

Let $\mxU_{\mxX_k}\mxSig_{\mxX_k}\mxV_{\mxX_k}^T=\mxX_{k}$ be the SVD of $\mxX_k$; thanks to the invariance of the Frobenius norm under orthogonal transformations, the regularization term in the above problem can be rewritten as 
\begin{eqnarray*}
\left\|(\mxX_{k}\mxX_{k}^T+\gamma\mxI)^{p/4-1/2}\mxX\right\|_F^2 
\!\! = \! \left\| \mxU_{\mxX_k}(\mxSig_{\mxX_k}^2+\gamma\mxI)^{p/4-1/2}\mxU_{\mxX_k}^T\mxX\right\|_F^2\!\!
= \!   \left\| (\mxSig_{\mxX_k}^2+\gamma\mxI)^{p/4-1/2}\mxU_{\mxX_k}^T\mxX\mxV_{\mxX_k}\right\|_F^2.
\end{eqnarray*}
Using well-known Kronecker product properties
\begin{eqnarray*}
\left\| (\mxSig_{\mxX_k}^2+\gamma\mxI)^{p/4-1/2}\mxU_{\mxX_k}^T\mxX\mxV_{\mxX_k}\right\|_F^2=
\left\| \vect\left((\mxSig_{\mxX_k}^2+\gamma\mxI)^{p/4-1/2}\mxU_{\mxX_k}^T\mxX\mxV_{\mxX_k}\right)\right\|_2^2\\
=\left\| \left(\mxV_{\mxX_k}^T\otimes\left((\mxSig_{\mxX_k}^2+\gamma\mxI)^{p/4-1/2}\mxU_{\mxX_k}^T\right)\right)\bfx\right\|_2^2
=\left\| \left(\mxI\otimes(\mxSig_{\mxX_k}^2+\gamma\mxI)^{p/4-1/2}\right)\left(\mxV_{\mxX_k}^T\otimes\mxU_{\mxX_k}^T\right)\bfx\right\|_2^2\,.
\end{eqnarray*}
Problem (\ref{eq:nnr_irn_matr}) is therefore equivalent to
\begin{equation}\label{eq:nnr_irn}
\bfx_{k+1}=\arg\min_{\bfx}\left\|
\mxA\bfx-\bfb\right\|_2^2 + \hl\| \underbrace{\left(\mxI\otimes(\mxSig_{\mxX_k}^2+\gamma\mxI)^{p/4-1/2}\right)}_{=:(\mxW_p^\gamma)_k}\overbrace{\left(\mxV_{\mxX_k}^T\otimes\mxU_{\mxX_k}^T\right)}^{=:\mxS_k}\bfx\|_2^2\,.
\end{equation}
In the above formulation, $(\mxW_p^\gamma)_k$ is a diagonal weighting matrix 
and $\mxS_k$ is an orthogonal matrix; both $(\mxW_p^\gamma)_k$ and $\mxS_k$ depend on the current approximation $\bfx_k$ of the solution $\bfx$. Intuitively, the matrix $\mxS_k$ maps $\bfx$ into the ``singular value domain'' of $\mxX_k$ (and acts as an iteration-dependent sparsity transform), and the matrix $(\mxW_p^\gamma)_k$ assigns suitable weights that allow to approximate a $p$-norm of the singular values. Therefore, the penalization term in (\ref{eq:nnr_irn}) can be interpreted as a reweighted vectorial 2-norm, with respect to a transformation of the solution $\bfx$. 
%
%
For this reason, the proposed approach is dubbed ``IRN-NNR$p$'' and is summarized in Algorithm \ref{alg:irn_nnr}.
\begin{algorithm}[H]
\caption{IRN-NNR$p$}\label{alg:irn_nnr}
\begin{algorithmic}[1]
\State{Inputs: $\mxA$, $\bfb$, $(\mxW_p^\gamma)_0=\mxI$, $\mxS_0=\mxI$}
\For{$k=0,1,\dots$ until a stopping criterion is satisfied}
\State{Solve problem (\ref{eq:nnr_irn})}
\State{``Decrease'' $\gamma$}
\State{Update $(\mxW_p^\gamma)_{k+1}$ and $\mxS_{k+1}$}
\EndFor
 \end{algorithmic}
\end{algorithm}
The next subsection derives new strategies for the efficient solution of the sequence of sub-problems (\ref{eq:nnr_irn}) appearing in Algorithm \ref{alg:irn_nnr}.

\subsection{Solution of problem (\ref{eq:nnr_irn}) via Krylov methods} \label{sec:irnmethods} 

First rewrite problem (\ref{eq:nnrp}) using an appropriate change of variable as
\begin{equation}\label{eq:nnr4}
\widehat{\bfx}_{k+1}=\arg\min_{\widehat{\bfx}} \|\mxA\mxS_k^T(\mxW_p^\gamma)_k^{-1}\widehat{\bfx} - \bfb\|_2^2 + \hl\| \widehat{\bfx}\|_2^2,\ \mbox{with}\ \widehat{\bfx}=(\mxWpg)_k\mxS_k\bfx\,.
\end{equation}
Note that 
\begin{equation}\label{SWmatrix}
\mxS_k^T=\mxS_k^{-1}=\mxV_{\mxX_k}\otimes \mxU_{\mxX_k}\quad\mbox{and}\quad (\mxW_p^\gamma)_k^{-1}=\mxI\otimes(\mxSig_{\mxX_k}^2+\gamma\mxI)^{1/2-p/4}\,,
\end{equation}
so that the above transformations (inversion of an orthogonal and a diagonal matrix) are numerically affordable by exploiting properties of Kronecker products. The Tikhonov-regularized problem (\ref{eq:nnr4}) in standard form is equivalent to the Tikhonov-regularized problem (\ref{eq:nnr_irn}) in general form. Many Krylov subspace methods based on the Golub-Kahan Bidiagonalization (GKB) or Arnoldi algorithms 
can be employed to approximate the solution of (\ref{eq:nnr4}). Moreover, if the regularization parameter $\hl$ is not known a priori, many efficient strategies to set its value adaptively within the sequence of projected problems can be used (i.e., in the framework of hybrid methods; see \cite{KO,gazzola2014}). 
The matrices $\mxS_k$ and $(\mxWpg)_k^{-1}$ can be formally thought of as preconditioners for the original problem (\ref{eq:Axeqb}), whose purpose is to enforce additional regularization into the solution subspace, 
rather than speeding-up the convergence of linear solvers applied to (\ref{eq:Axeqb}). 

\paragraph{Methods based on the GKB algorithm} The $m$th step of the GKB algorithm applied to the matrix $\mxA\mxS_k^T(\mxW_p^\gamma)_k^{-1}$ with starting vector $\bfb$ (i.e., taking $\bfx_0=\bfo$) can be expressed by the following partial matrix factorizations
\begin{equation}\label{GKB}
(\mxA\mxS_k^T(\mxW_p^\gamma)_k^{-1})\mxV_m=\mxU_{m+1}\wbB_m\quad \mbox{and}\quad
((\mxW_p^\gamma)_k^{-1}\mxS_k\mxA^T)\mxU_{m+1}=\mxV_{m+1}\mxB_{m+1}^T,
\end{equation}
where $\mxU_{j}\in\RR^{M\times j}$ and $\mxV_{j}\in\RR^{N\times j}$ (with $j=m,m+1$ and $\mxU_j\bfe_1=\bfb/\|\bfb\|_2$) have orthonormal columns, and $\mxB_{m+1}\in\RR^{(m+1)\times (m+1)}$ is lower bidiagonal (with $\wbB_m$ obtained by removing the last column of $\mxB_{m+1}$). The orthonormal columns of $\mxV_m$ are such that 
\[
\mcR(\mxV_{m}) = \K_m\left(
((\mxW_p^\gamma)_k^{-1}\mxS_k\mxA^T)(\mxA\mxS_k^T(\mxW_p^\gamma)_k^{-1}),
((\mxW_p^\gamma)_k^{-1}\mxS_k\mxA^T)\bfb
\right)\,.
\]
We find an approximate solution of (\ref{eq:nnr4}) by imposing ${\whx\in\mcR(\mxV_{m})}$, i.e., $\whx_m=\mxV_{m}\bfy_m$, where, by exploiting the first decomposition in (\ref{GKB}) and the properties of the matrices appearing therein, $\bfy_m\in\RR^m$ is such that
\begin{equation}\label{GKBhybrid}
\bfy_m=\arg\min_{\bfy\in\RR^m} \|\wbB_m{\bfy} - \|\bfb\|_2\bfe_1\|_2^2 + \hl_m\| \bfy\|_2^2\,.
\end{equation}
We used the notation $\hl_m$ for the regularization parameter to highlight that its value can be adaptively set within the iterations. 
The approximate solution to problem (\ref{eq:nnr_irn}) is such that
\begin{equation}
\bfx=\mxS_k^T(\mxW_p^\gamma)_k^{-1}\widehat{\bfx} \in \K_m\left(
(\mxS_k^T(\mxW_p^\gamma)_k^{-2}\mxS_k)\mxA^T\mxA,
(\mxS_k^T(\mxW_p^\gamma)_k^{-2}\mxS_k)\mxA^T\bfb\right)\,.
\label{GKBspace}
\end{equation}
%
%
%
Looking at the above approximation subspace for the solution $\bfx$, it is evident that the ``preconditioner'' acts by first mapping into the ``singular value domain'' (by applying $\mxS_k$), enforcing sparsity in the singular values (by reweighting with $(\mxWpg)_k^{-2}$), and eventually transforming back into the ``solution domain'' (by applying $\mxS_k^T$).

\paragraph{Methods based on the Arnoldi algorithm} If $\mxA$ is square, the $m$th step of the Arnoldi algorithm applied to the matrix $\mxA\mxS_k^T(\mxW_p^\gamma)_k^{-1}$ with starting vector $\bfb$ (i.e., taking $\bfx_0=\bfo$) can be expressed by the following partial matrix factorization
\begin{equation}\label{Arnoldi}
(\mxA\mxS_k^T(\mxW_p^\gamma)_k^{-1})\mxV_m=\mxV_{m+1}\mxH_m
\end{equation}
where $\mxV_{j}\in\RR^{N\times j}$ (with $j=m,m+1$ and $\mxV_j\bfe_1=\bfb/\|\bfb\|_2$) 
have orthonormal columns such that
\[
\mcR(\mxV_{m}) = \K_m\left(
\mxA\mxS_k^T(\mxW_p^\gamma)_k^{-1},
\bfb\right)\,,
\]
and $\mxH_{m}\in\RR^{(m+1)\times m}$ is upper Hessenberg. Similarly to the GKB case, we find an approximate solution of (\ref{eq:nnr4}) by imposing ${\whx\in\mcR(\mxV_{m})}$ and by solving a projected Tikhonov problem of order $m$. 
The approximate solution to problem (\ref{eq:nnr_irn}) is such that
$$
\bfx=\mxS_k^T(\mxW_p^\gamma)_k^{-1}\widehat{\bfx} \in \mxS_k^T(\mxWpg)_k^{-1} \K_m\left(
\mxA\mxS_k^T(\mxW_p^\gamma)_k^{-1},
\bfb\right)\,,
$$
where  
$$
\mxS_k^T(\mxWpg)_k^{-1} \K_m\left( \mxA\mxS_k^T(\mxW_p^\gamma)_k^{-1} \right)=\mbox{span}\{\mxS_k^T(\mxW_p^\gamma)_k^{-1}\bfb,\dots,
\left(\mxS_k^T(\mxW_p^\gamma)_k^{-1}\mxA\right)^{m-1}\mxS_k^T(\mxW_p^\gamma)_k^{-1}\bfb\}.
$$
Contrarily to the GKB case, 
we immediately notice that $\bfx$ does not belong to a meaningful approximation subspace. Indeed, just by looking at the first vector: $\bfb$ is in the image space and $(\mxWpg)_k^{-1}$ is supposed to act on the singular value space of $\mxX_k$, so $(\mxWpg)_k^{-1}\bfb$ does not make any sense; furthermore, $\mxS_k^T$ is supposed to link the singular value space of $\mxX_k$ to the image space, so $\mxS_k^T(\mxWpg)_k^{-1}\bfb$ does not make sense either. 
Similarly to what is proposed in \cite{belge2000wavelet,chunggazzola2018}, where the Arnoldi algorithm is applied to a regularized problem that enforces sparsity in the wavelet domain, we propose to fix this issue by incorporating $\mxS_k$ also as an orthogonal left ``preconditioner'' for the original system (\ref{eq:Axeqb}) so that, by exploiting the invariance of the vectorial 2-norm under orthogonal transformations, problem (\ref{eq:nnr4}) can be equivalently reformulated as
\begin{equation}\label{eq:nnr5}
\widehat{\bfx}_{k+1}=\arg\min_{\widehat{\bfx}} \|\mxS_k(\mxA\mxS_k^T(\mxW_p^\gamma)_k^{-1}\widehat{\bfx} - \bfb)\|_2^2 + \hl\| \widehat{\bfx}\|_2^2,\ \mbox{with}\ \widehat{\bfx}=(\mxWpg)_k\mxS_k\bfx\,.
\end{equation}
The (right and left) preconditioned Arnoldi algorithm applied to problem (\ref{eq:nnr5}) can now be expressed by the following partial matrix factorization
\begin{equation}\label{eq:Arnoldi-new}
(\mxS_k\mxA\mxS_k^T(\mxW_p^\gamma)_k^{-1})\mxV_m=\mxV_{m+1}\mxH_m\,.
\end{equation}
%
%
%
We find an approximate solution of (\ref{eq:nnr5}) by imposing ${\whx\in\mcR(\mxV_{m})}=\K_m(\mxS_k\mxA\mxS_k^T\mxW^{-1},\mxS_k\bfb)$, i.e., $\whx_m=\mxV_{m}\bfy_m$, where, by exploiting 
(\ref{eq:Arnoldi-new}) and the properties of the matrices appearing therein, $\bfy_m\in\RR^m$ is such that
\begin{equation}\label{Arnoldihybrid}
\bfy_m=\arg\min_{\bfy\in\RR^m} \|\mxH_m{\bfy} - \|\bfb\|_2\bfe_1\|_2^2 + \hl_m\| \bfy\|_2^2\,.
\end{equation}
Hence
\begin{equation}
\bfx 
\in \nonumber
\mxS_k^T(\mxWpg)_k^{-1}\K_m(\mxS_k\mxA\mxS_k^T(\mxWpg)_k^{-1},\mxS_k\bfb)\,.
\label{Arnoldispace}
\end{equation}
which is suitable for approximating the solution. 
%
%
%
The new methods based on the GKB algorithm (for generic matrices) and Arnoldi algorithm (only if $\mxA\in\RR^{N\times N}$) are dubbed ``IRN-LSQR-NNR$p$'' and ``IRN-GMRES-NNR$p$'', respectively, and are summarized in Algorithm \ref{alg:irn_krylov}. 
\begin{algorithm}[H]
\caption{IRN-LSQR-NNR$p$ and IRN-GMRES-NNR$p$}\label{alg:irn_krylov}
\begin{algorithmic}[1]
\State{Inputs: $\mxA$, $\bfb$, $(\mxW_p^\gamma)_0=\mxI$, $\mxS_0=\mxI$}
\For{$k=0,1,\dots$ until a stopping criterion is satisfied}
\For{$m=1,2,\dots$ until a stopping criterion is satisfied}
\State{Update the factorizations (\ref{GKB}) and (\ref{eq:Arnoldi-new}), respectively}
\State{Solve the projected problem (\ref{GKBhybrid}) and (\ref{Arnoldihybrid}), respectively, tuning $\hl_m$ if necessary}
\EndFor
\State{``Decrease'' $\gamma$}
\State{Update the new $(\mxW_p^\gamma)_{k+1}$ and $\mxS_{k+1}$}
\EndFor
 \end{algorithmic}
\end{algorithm}

\subsection{Solution through flexible Krylov subspaces}\label{sec:flexinnr}


Problem (\ref{eq:nnrp}) reformulated as (\ref{eq:nnr4}) allows us to naturally apply the flexible Golub-Kahan (FGK) and flexible Arnoldi algorithms. Indeed, instead of updating the ``preconditioners'' $\mxS_k$ and $(\mxWpg)_k$ at the $k$th outer iteration of the nested iteration schemes of Algorithm \ref{alg:irn_krylov}, we propose to consider new ``preconditioners'' as soon as a new approximation of the solution is available, i.e., at each iteration of a Krylov subspace solver. Therefore, at the $(i+1)$th iteration of the new solvers, the ``preconditioners'' $(\mxWpg)_i$ and $\mxS_i$ are computed as in (\ref{SWmatrix}), but using the SVD of the $i$th approximate solution
\[
\mxX_{i}=\vect^{-1}(\bfx_{i})=\mxU_{\mxX_i}\mxSig_{\mxX_i}\mxV_{\mxX_i}^T,\quad\mbox{for $i=1,\dots,k-1$}\,,
\]
with $(\mxWpg)_0=\mxI$ and $\mxS_0=\mxI$. In order to incorporate iteration-dependent preconditioning, the flexible versions of the Golub-Kahan and Arnoldi factorizations have to be used. 

Namely, at the $i$th iteration, 
the new instance of the FGK algorithm updates partial factorizations of the form (\ref{eq:gkfac}), i.e., $\mxA\mxZ_i = \mxU_{i+1}\mxM_i$ and $\mxA^T\mxU_{i+1} = \mxV_{i+1}\mxT_{i+1}$, where
\begin{equation}\nonumber
\mxZ_i=[\mxS_0^T(\mxWpg)_0^{-2}\mxS_0\bfv_1,\dots,\mxS_{i-1}^T(\mxWpg)_{i-1}^{-2}\mxS_{i-1}\bfv_i]\,,\quad \bfv_1=\mxA^T\bfb/\|\mxA^T\bfb\|_2\,.
\end{equation}
Taking $\bfx_0=\bfo$, the $i$th approximate solution is such that $\bfx_i=\mxZ_i\bfy_i$, where 
\begin{equation}\label{FGKhybrid}
\bfy_i=\arg\min_{\bfy\in\RR^i} \|\mxM_i{\bfy} - \|\bfb\|_2\bfe_1\|_2^2 + \hl_i\| \bfy\|_2^2\,.
\end{equation}
Note that the subspace for the solution $\mcR(\mxZ_i)$ can be regarded as a generalization of the subspace (\ref{GKBspace}) computed when considering preconditioned GKB within the IRN-LSQR-NNR$p$ method. The new method is dubbed ``FLSQR-NNR$p$'', and is summarized in Algorithm \ref{alg:fnnrp}. 

For $\mxA\in\RR^{N\times N}$ and $\bfx_0=\bfo$, at the $i$th iteration, the new instance of the flexible Arnoldi algorithm updates a partial factorization of the form (\ref{eq:arnoldifac}), with $k=i$, and generates
\[
\mxZ_i=[\mxS_0^T(\mxWpg)_0^{-1}\mxS_0\bfv_1,
\dots,\mxS_{i-1}^T(\mxWpg)_{i-1}^{-1}\mxS_{i-1}\bfv_i]\,,\quad \bfv_1=\bfb/\|\bfb\|_2\,,
\]
where both right and left preconditioners are used analogously to IRN-GMRES-NNR$p$. The $i$th approximate solution is such that $\bfx_i=\mxZ_i\bfy_i$, where 
\begin{equation}\label{FArnoldihybrid}
\bfy_i=\arg\min_{\bfy\in\RR^i} \|\mxH_i{\bfy} - \|\bfb\|_2\bfe_1\|_2^2 + \hl_i\| \bfy\|_2^2\,.
\end{equation}
Note that the subspace for the solution $\mcR(\mxZ_i)$ can be regarded as a generalization of the subspace (\ref{Arnoldispace}) computed when considering the preconditioned Arnoldi algorithm within the IRN-GMRES-NNR$p$ method. The new method is dubbed ``FGMRES-NNR$p$'', and is summarized in Algorithm \ref{alg:fnnrp}.
%
\begin{algorithm}[H]
\caption{FLSQR-NNR$p$ and FGMRES-NNR$p$}\label{alg:fnnrp}
\begin{algorithmic}[1]
\State{Inputs: $\mxA$, $\bfb$, $(\mxW_p^\gamma)_0=\mxI$, $\mxS_0=\mxI$}
\For{$i=1,2,\dots$ until a stopping criterion is satisfied}
\State{Update a factorization of the form (\ref{eq:gkfac}) and (\ref{eq:arnoldifac}), respectively, to expand the space $\mcR(\mxZ_i)$}
\State{Solve the projected problem (\ref{FGKhybrid}) and (\ref{FArnoldihybrid}), respectively, tuning $\hl_i$ if necessary}
\State{``Decrease'' $\gamma$}
\State{Update the new $(\mxW_p^\gamma)_{i}$ and $\mxS_{i}$, using the SVD $\mxX_{i}=\vect^{-1}(\bfx_i)=\mxU_{\mxX_i}\mxSig_{\mxX_i}\mxV_{\mxX_i}^T$.}
\EndFor
 \end{algorithmic}
 \end{algorithm}
Note that, although the approach of Algorithm \ref{alg:fnnrp} is quite heuristic, it avoids nested iteration cycles and computes only one approximation subspace for the solution of (\ref{eq:nnrp}), where low-rank penalization is adaptively incorporated. Because of this, in many situations, Algorithm \ref{alg:fnnrp} computes solutions of quality comparable to the ones computed by Algorithm \ref{alg:irn_krylov}, with a significant reduction in the number of iterations. We should also mention that, in the framework of affine rank minimization problems, \cite{mohan2012} outlines an algorithm that avoids inner projected gradient iterations for the solution of each quadratic subproblem in the sequence generated within the IRN strategy. 

Finally, we underline that, within the framework of flexible Krylov subspaces, the approximation subspaces $\mcR(\mxZ_i)$ for the $i$th approximate solution can be further modified, with some insight into the desired properties of the solution. Indeed, since the $i$th basis vector for the solution is of the form \linebreak[4]$\bfz_i=\mxS_{i-1}^T(\mxWpg)_{i-1}^{-2}\mxS_{i-1}\bfv_i$ for FLSQR-NNR$p$, and $\bfz_i=\mxS_{i-1}^T(\mxWpg)_{i-1}^{-1}\mxS_{i-1}\bfv_i$ for FGMRES-NNR$p$, one can consider alternative ``preconditioners'' $\mxS_{i-1}$ and $(\mxWpg)_{i-1}$ that are still effective in delivering low-rank solutions. For instance, focusing on FGMRES, and given $\bfv_i=\mxV_i\bfe_i$, where $\mxV_i$ is the matrix appearing on the right-hand side of the factorization (\ref{eq:gkfac}), 
and given the SVD of $\vect^{-1}(\bfv_i)=\mxU_{\mxV_i}\mxSig_{\mxV_i}\mxV_{\mxV_i}^T$, one can take
\begin{equation}\label{SWmatrix-var}
\mxS_{i-1}=\mxV_{\mxV_i}^T\otimes \mxU_{\mxV_i}^T\quad
\mbox{and}
\quad(\mxW_p^\gamma)_{i-1}^{-1}=\mxI\otimes(\mxSig_{\mxV_i})^{1/2-p/4}\,,
\end{equation}
and as a result,
\begin{eqnarray*}
\mxS_{i-1}\bfv_i &=& \vectorize(\mxU_{\mxV_i}^T\vect^{-1}(\bfv_i)\mxV_{\mxV_i}) = \vectorize(\mxSig_{V_i}),\\
(\mxWpg)_{i-1}^{-1}\mxS_{i-1}\bfv_i&=&\vectorize((\mxSig_{\mxV_{i}})^{1/2-p/4}\mxSig_{\mxV_i}) = \vectorize((\mxSig_{\mxV_{i}})^{3/2-p/4}),\\
\mxS_{i-1}^T(\mxWpg)_{i-1}^{-1}\mxS_{i-1}\bfv_i&=&
\vectorize(\mxU_{\mxV_i}((\mxSig_{\mxV_{i}})^{3/2-p/4})\mxV_{\mxV_i}^T)=\bfz_i,
\end{eqnarray*}
so that the singular values of $\vect^{-1}(\bfv_i)$ are rescaled: taking $0<p\leq 1$, the power of $\mxSig_{\mxV_{i}}$, $3/2-p/4$, is always larger than 1, which means that large singular values get magnified and small singular values become even smaller. In this way, the gaps between singular values are emphasized and to some extent contribute to the low rank properties of the basis vectors. Similar derivations hold for FLSQR. Hence, methods analogous to LR-FLSQR and LR-FGMRES are obtained, and are dubbed FGMRES-NNR$p$(v) and FLSQR-NNR$p$(v), respectively. 

\section{Implementation details} \label{sec:implementation}

All the methods considered in this paper are iterative, and therefore at least one suitable stopping criterion should be set for the iterations. When considering hybrid formulations (like the ones in Algorithms \ref{alg:irn_krylov} and \ref{alg:fnnrp}), one could simultaneously set a good value for the regularization parameter $\hl_j$ at the $j$th iteration, as well as properly stop the iterations. Strategies for achieving this are already available in the literature (see \cite{gazzola2017irtools,gazzola2014}). 

Assuming that a good estimate for the norm of the noise $\bfeta$ affecting the right-hand-side of (\ref{eq:Axeqb}) is available, i.e., $\eps\simeq\|\bfeta\|_2$, one can consider the discrepancy principle and stop the iterative scheme at the first iteration $j$ such that
\begin{equation}\label{discrP}
\|\bfb-\mxA\bfx_j\|_2\leq\theta\eps\,,\quad\mbox{where $\theta>1$, $\theta\simeq1$ is a safety threshold.}
\end{equation}
Applying the discrepancy principle to LR-FGMRES (Algorithm \ref{alg:lr-fgmres}) and LR-FLSQR (Algorithm \ref{alg:lr-flsqr}) is particularly convenient, as the norm of the residual on the left-hand side of (\ref{discrP}) can be monitored using projected quantities, i.e., 
\[
\left\|\|\bfb\|_2\bfe_1-\mxH_j\bfy_j\right\|_2\quad\mbox{for LR-FGMRES}\quad\mbox{and}\quad
\left\|\|\bfb\|_2\bfe_1-\mxM_j\bfy_j\right\|_2\quad\mbox{for LR-FLSQR},
\]
where decompositions (\ref{eq:arnoldifac}) and (\ref{eq:gkfac}), respectively, and the properties of the matrices appearing therein, have been exploited. When running hybrid methods (see Algorithms \ref{alg:irn_krylov} and \ref{alg:fnnrp}), we employ the so-called ``secant method'',  which updates the regularization parameter for the projected problem in such a way that stopping by the discrepancy principle is ensured. We highlight again that the quantities needed to implement the ``secant method'' (namely, the norm of the residual and the discrepancy associated to (\ref{eq:nnr4}) at each iteration) can be conveniently monitored using projected quantities: this is obvious for IRN-LSQR-NNR$p$ and FLSQR-NNR$p$, as only right-``preconditioning'' is employed; it is less obvious for IRN-GMRES-NNR$p$ and FGMRES-NNR$p$, but since the left-``preconditioner'' is orthogonal, one can still write
\[
\|\bfb - \mxA\bfx_j\|_2=\|\mxS_k\bfb - \mxS_k\mxA\mxS_k^T(\mxWpg)_k^{-1}\widehat{\bfx}_j\|_2=\| \|\bfb\|_2\bfe_1 - \mxH_j\bfy_j \|_2\,.
\]
Note that all the methods in Algorithm \ref{alg:irn_krylov} and \ref{alg:fnnrp} can also run with $\hl=0$, and still achieve low-rank approximate solutions: this is because the approximation subspace for the solution incorporates regularizing ``preconditioning'' (see \cite{HH,HanJen06} for details on this approach in the case on smoothing ``preconditioning'' with finite-difference approximations of derivatives operators). Finally, when dealing with the inner-outer iteration scheme of Algorithm \ref{alg:irn_krylov}, in addition to a parameter choice strategy and stopping criterion for the hybrid projected problems (\ref{GKBhybrid}) and (\ref{Arnoldihybrid}), one should also consider a
stopping criterion for the outer iterations. 
We propose to do this by monitoring the norm of the difference of the singular values (normalized by the largest singular value so that $\sigma_1(\mxSig_{\mxX_{k+1}}) = \sigma_1(\mxSig_{\mxX_k})=1$) of two approximations of the solution of (\ref{eq:nnrp}) obtained at two consecutive outer iterations of Algorithm \ref{alg:irn_krylov}, i.e., we stop as soon as
\begin{equation}\label{stopOUT}
\|\mbox{diag}(\mxSig_{\mxX_{k+1}})-\mbox{diag}(\mxSig_{\mxX_{k}})\|_2<\tau_{\sigma},
\quad k=1,2,\dots,\quad
\end{equation}
where 
$\vect^{-1}(\bfx_i)=\mxX_i=\mxU_{\mxX_i}\mxSig_{\mxX_i}\mxV_{\mxX_i}^T$ ($i=k,k+1$), and $\tau_\sigma>0$ is a user-specified threshold. 
If no significant changes happen in the rank and singular values of two consecutive approximations of the solution, then (\ref{stopOUT}) is satisfied. 


We conclude this section with a few remarks about the computational cost of the proposed methods. Note that, if $\mxA\in\RR^{N\times N}$, IRN-GMRES-NNR$p$ is intrinsically cheaper than IRN-LSQR-NNR$p$ (since, at each iteration, the former requires only one matrix-vector product with $\mxA$, while the latter requires one matrix-vector product with $\mxA$ and one with $\mxA^T$). However, methods based on the Arnoldi algorithm are typically less successful than methods based on the GKB algorithm for regularization; see \cite{JenHan07}.  
Other key operations for implementing our proposed methods are the computation of the SVDs of relevant quantities, and/or the application of the ``preconditioners'' in (\ref{SWmatrix-var}). Namely, 
each iteration of LR-FGMRES, LR-FLSQR, FLSQR-NNR$p$, and FGMRES-NNR$p$ requires the computation of the SVD of an $n\times n$ matrix, which amounts to $O(n^3)$ floating point operations. When considering IRN-LSQR-NNR$p$ and IRN-GMRES-NNR$p$, only the SVD of the approximate solution should be computed once at each outer iteration. However, each inner iteration of IRN-LSQR-NNR$p$ and IRN-GMRES-NNR$p$, as well as each iteration of FLSQR-NNR$p$ and FGMRES-NNR$p$, requires the computation of matrix-vector products of the form $\mxS_k^T(\mxWpg)_k^{-1}\bfv_i$: 
this can be achieved within a two-step process, where first the rescaling $\widetilde{\bfv}_i=(\mxWpg)_k^{-1}\bfv_i$ is applied with $O(N)=O(n^2)$ floating-point operations, and then $\mxS_k^T\widetilde{\bfv}_i=(\mxV_{\mxX_{k}}\otimes\mxU_{\mxX_{k}})\widetilde{\bfv}_i$ is computed. While a straightforward implementation of the latter would require $\bigO(N^2)=\bigO(n^4)$ floating-point operations, exploiting Kronecker product properties can bring down the cost of this operation to $\bigO(n^3) = \bigO(N^{3/2})$, by computing $\mxS_k^T\widetilde{\bfv}_i=\vectorize(\mxU_{\mxX_{k}}^T\vect^{-1}(\bfv_i)\mxV_{\mxX_{k}})$. 
%



\section{Experimental Results}\label{sec:experiments}
In this section, we present results of numerical experiments on several image processing problems to demonstrate the performance of the new  IRN-GMRES-NNR$p$, IRN-LSQR-NNR$p$, FGMRES-NNR$p$, and FLSQR-NNR$p$ methods. Variants of FGMRES-NNR$p$ and FLSQR-NNR$p$ (marked with ``(v)") are also tested. To shorten the acronyms, we omit $p$ when $p=1$, which means IRN-GMRES-NNR denotes IRN-GMRES-NNR$p$ when $p=1$, etc. Examples are generated using \emph{IR Tools} \cite{gazzola2017irtools}. 

In general, we compare the performances of the proposed methods to standard Krylov subspace methods GMRES and LSQR, also used in a hybrid fashion. We also test against the low-rank projection methods described in Section \ref{sec:bg} and the singular value thresholding (SVT) algorithm  \cite{cai2010}, which was originally proposed for low-rank matrix completion problems, and can be extended to problems with linear constraints of the form
%
%
%
\begin{equation}\label{eq:SVT}
\min_{\bfx}   \tau\|\vect^{-1}(\bfx)\|_*+\frac{1}{2}\|\vect^{-1}(\bfx)\|_F^2\quad
 \mbox{subject\ to}\quad
 \mxA\bfx=\bfb,\mbox{ where $\vect^{-1}(\bfx)=\mxX$}.
\end{equation}
The $k$th iteration of the SVT algorithm for (\ref{eq:SVT}) reads
 \begin{equation}\label{eq:SVTalg}
\begin{cases} \mxX_k = \D_\tau(\mxA^T\bfy_{k-1})\\ \bfy_k = \bfy_{k-1}+\delta_k(\bfb-A\bfx_k) \end{cases}, 
\end{equation}
where $\delta_k$ is a step size and $\D_\tau$ is the {\it singular value shrinkage operator}, defined as 
\begin{equation}\nonumber
\D_\tau(\mxX) = \mxU_{\mxX}\D_\tau(\mxSig_{\mxX}) \mxV_{\mxX}^T,\ \ \D_\tau(\mxSig_{\mxX}) = 
\max\{\mxSig_{\mxX}-\tau\mxI,\bfo\},
\end{equation}
where $\mxX=\mxU_{\mxX}\mxSig_{\mxX}\mxV_{\mxX}^T$ is the SVD of $\mxX$, $\bfo$ is a matrix of zeros, and the maximum is taken component-wise. Although (\ref{eq:SVTalg}) is not the same problem as (\ref{eq:nnr}), they are similar in that both penalize the nuclear norm of $\vect^{-1}(\bfx)$ and they respect the constraint $\mxA\bfx=\bfb$.

Regarding the comparisons with the low-rank projection methods presented in Section \ref{sec:bg}, 
there are no established and theoretically informed ways of choosing the truncation ranks for the solutions and for the basis vectors of the solution subspace. Hence, for all test problems, we experiment on a reasonable number of trials, each with different truncation rank choices, and select the best performing rank out of all ranks tested. For simplicity, we consider the same truncation rank for basis vectors and solutions ($\tau_{\kappa_B} = \tau_\kappa$). We follow the same process to choose the number of restarts and the number of iterations for each restart for RS-LR-GMRES, as well as the shrinkage threshold $\tau$ in SVT; strategies to select the step size for SVT are described in \cite{cai2010}. 

\paragraph{Example 1: Binary Star}

\begin{figure}[htbp]
\begin{center}
\begin{tabular}{cc}
\hspace{-0.8cm}\includegraphics[height=1.3in]{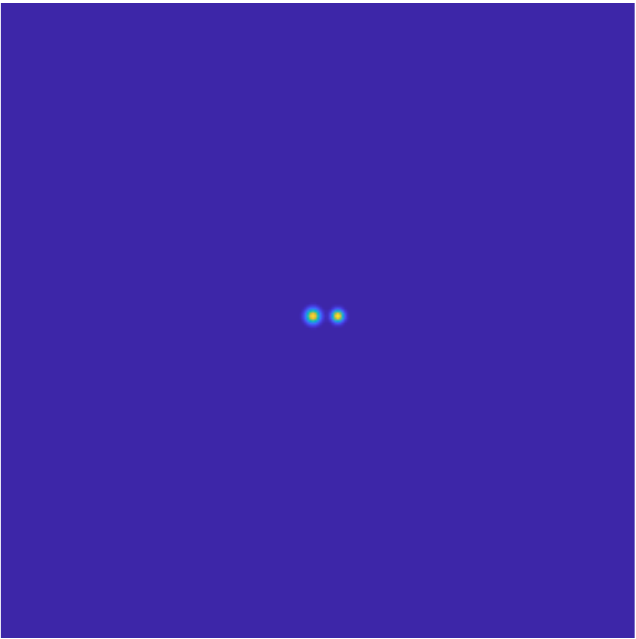}  &
\hspace{0.8cm}\includegraphics[height=1.3in]{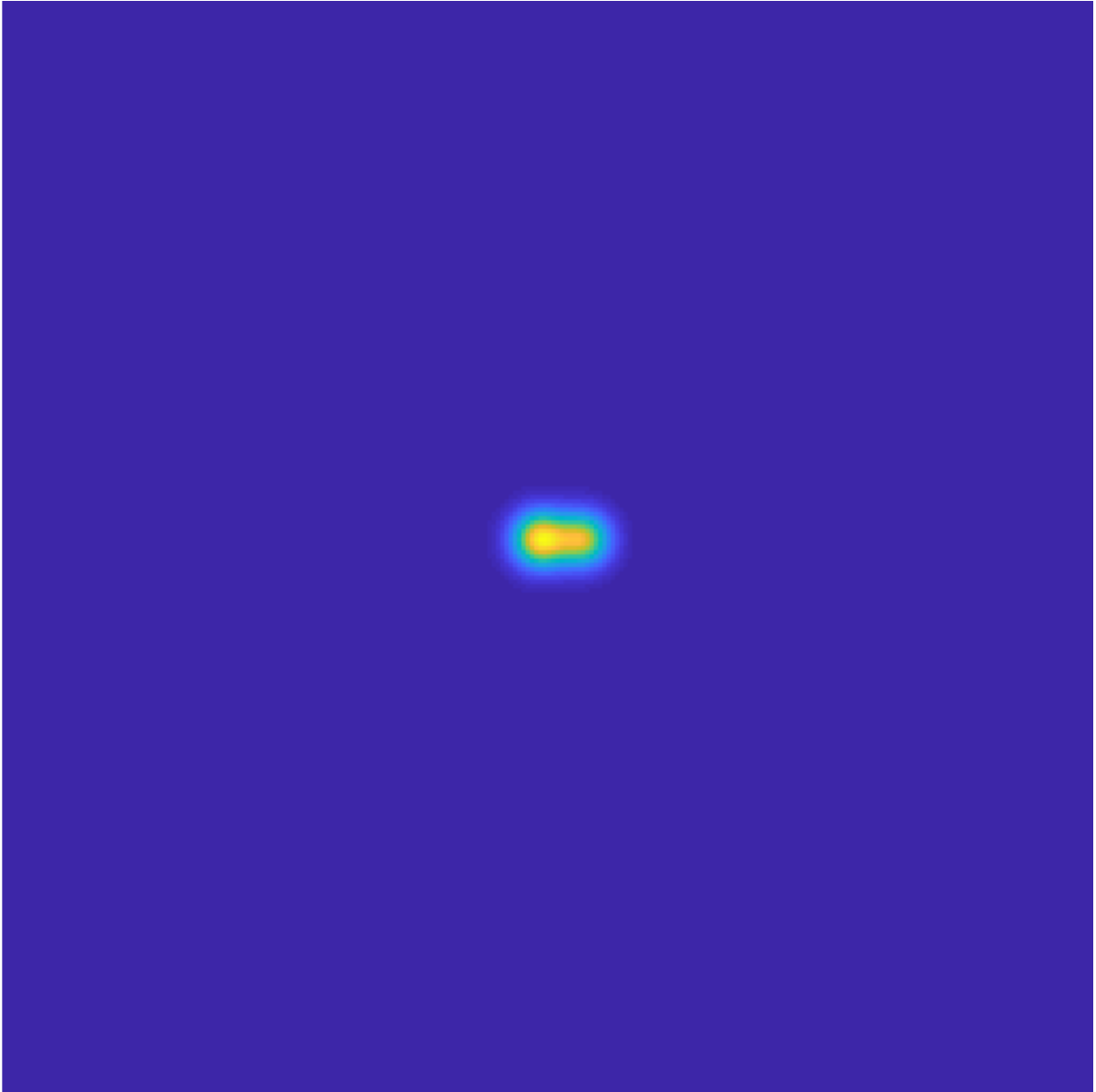}  \\
\hspace{-0.8cm}{\small exact} & 
\hspace{0.8cm}{\small blurred \& noise}\vspace{0.3cm}\\
\hspace{-0.8cm}\includegraphics[height=1.3in]{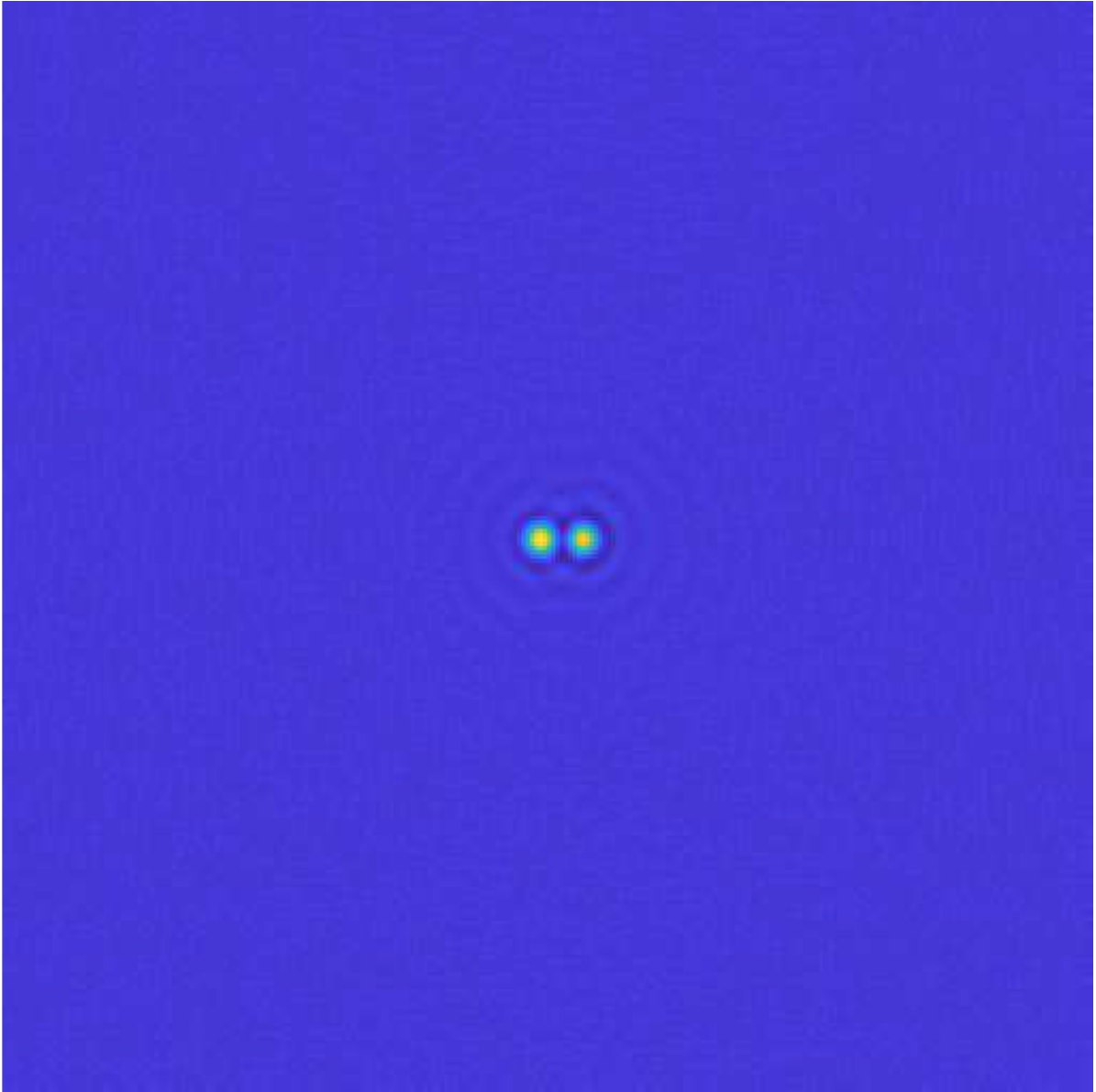} &
\hspace{0.8cm}\includegraphics[height=1.3in]{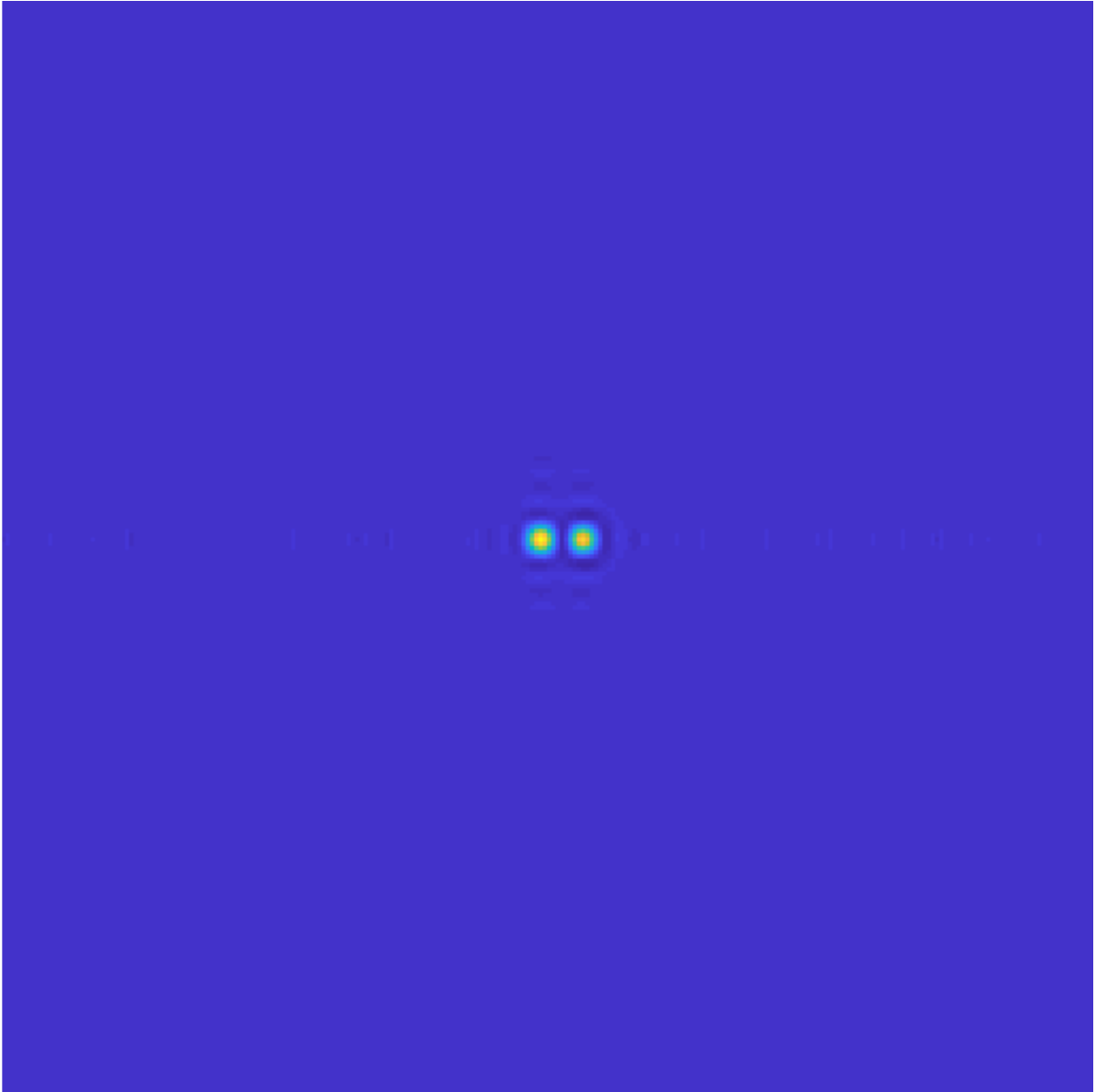} \\
\hspace{-0.8cm}{\small LR-FGMRES} & 
\hspace{0.8cm}{\small IRN-GMRES-NNR}\\
\end{tabular}
\end{center}
\caption{\emph{Example 1}. Exact and corrupted test images, together with the best reconstructions obtained by the LR-FGMRES and the IRN-GMRES-NNR methods.}
\label{fig:star_images}
\end{figure}

We consider an image deblurring problem involving a binary star test image of size $256\times 256$: this test image has rank 2. The true image is displayed in the leftmost frame of Figure \ref{fig:star_images}.  A standard Gaussian blur is applied to the test image, and Gaussian white noise of level $\|\bfeta\|_2/\|\bex\|_2=10^{-3}$ is added. The blurred and noisy image is shown in Figure \ref{fig:star_images}, second frame from the left. Due to the existence of noise, the blurred image has full rank. 
%
For this example, the blurring operator $\mxA$ is square of size $65536\times 65536$, hence GMRES-related methods are used for comparison, namely: GMRES, IRN-GMRES-NNR, FGMRES-NNR, LR-FGMRES and RS-LR-GMRES (i.e., we only consider the case $p=1$ here). SVT is also taken into consideration. The truncation rank for LR-FGMRES and RS-LR-GMRES is set to 30 for both basis vectors and approximate solutions (i.e., $\tau_{\kappa_B}=\tau_\kappa=30$). RS-LR-GMRES is restarted every 40 iterations. The step size for SVT is set to be $\delta_k=\delta = 2$ and the singular value shrinkage threshold $\tau$ is 1. 
Note that, although the true solution has only rank 2, setting truncation rank to 2 for low rank methods produces solutions of worse quality (compared to setting the rank to 30). This might be because of the inherent ill-posedness of the problem, which makes it harder to obtain solutions with desired properties (e.g., with rank 2): indeed, if we do truncate to rank 2, a lot of information about the solution might be lost.  

\begin{figure}[H]
\begin{center}
\begin{tabular}{c}
\includegraphics[width=12cm]{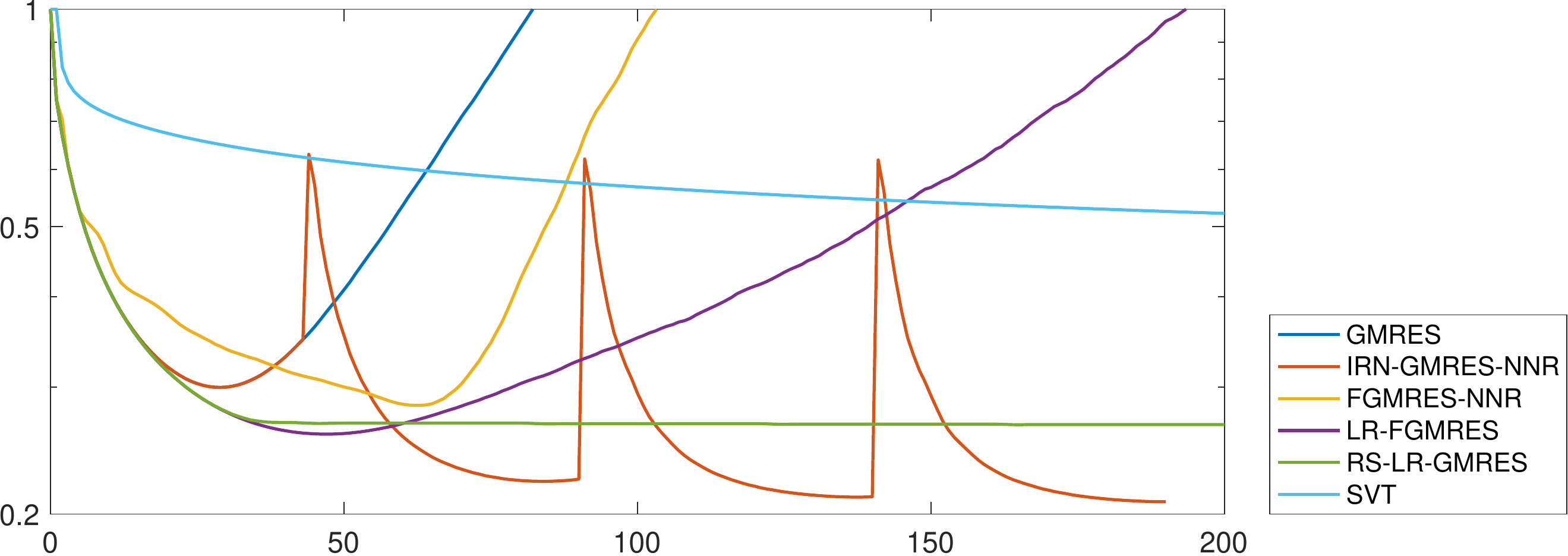}
\end{tabular}
\end{center}
\caption{{\em Example 1}. Relative errors vs. number of iterations for GMRES-based methods and SVT.}
\label{fig:star_error}
\end{figure}

Figure \ref{fig:star_error} displays the histories of relative errors $\|\xex-\bfx_m\|_2/\|\xex\|_2$ for the first 200 iterations (i.e., $m=1,\dots,200$) of these methods. For IRN-GMRES-NNR, 4 outer cycles were run, each with a maximum of 50 iterations: a new outer cycle is initiated as soon as the discrepancy principle is satisfied in the inner cycle. No additional regularization is used (i.e., $\widehat{\lambda} = 0$ for all methods). 

We can observe from Figure \ref{fig:star_error} that when the truncation ranks are chosen reasonably, LR-FGMRES and RS-LR-GMRES both produce a less pronounced semi-convergence behavior than GMRES, with LR-FGMRES attaining a smaller relative error than RS-LR-GMRES.  FGMRES-NNR, on the other hand, shows slower semi-convergence than GMRES, but it also converges to a slightly better relative error. IRN-GMRES-NNR behaves especially well in this case, with significantly reduced relative errors even at the end of the second outer cycle. 
The ``jumps" at the beginning of each outer IRN-GMRES-NNR iteration are due to the strategy used for restarts (the older basis vectors are cleared at each restart).

Figure \ref{fig:star_images} displays the exact and the corrupted images, as well as the best reconstructions computed by LR-FGMRES and IRN-GMRES-NNR: these are obtained at the 47th and the 189th (total) iteration of LR-FGMRES and IRN-GMRES-NNR, respectively. By looking at relative errors in Figure \ref{fig:star_error}, we see that LR-FGMRES is the second best out of all methods, and yet the quality of the solution is inferior compared to IRN-GMRES-NNR. 
Compared to the LR-FGMRES solution, the IRN-GMRES-NNR one is a more truthful reconstruction of the exact image: it not only has less artifacts immediately around the stars, but also has less background noise, in the sense that the pixel intensities in the background are closer to the true ones (as it can be seen by looking at the background color). 

\begin{figure}[htbp]
\begin{center}
\begin{tabular}{cc}
\hspace{-0.8cm}\includegraphics[height=1.3in]{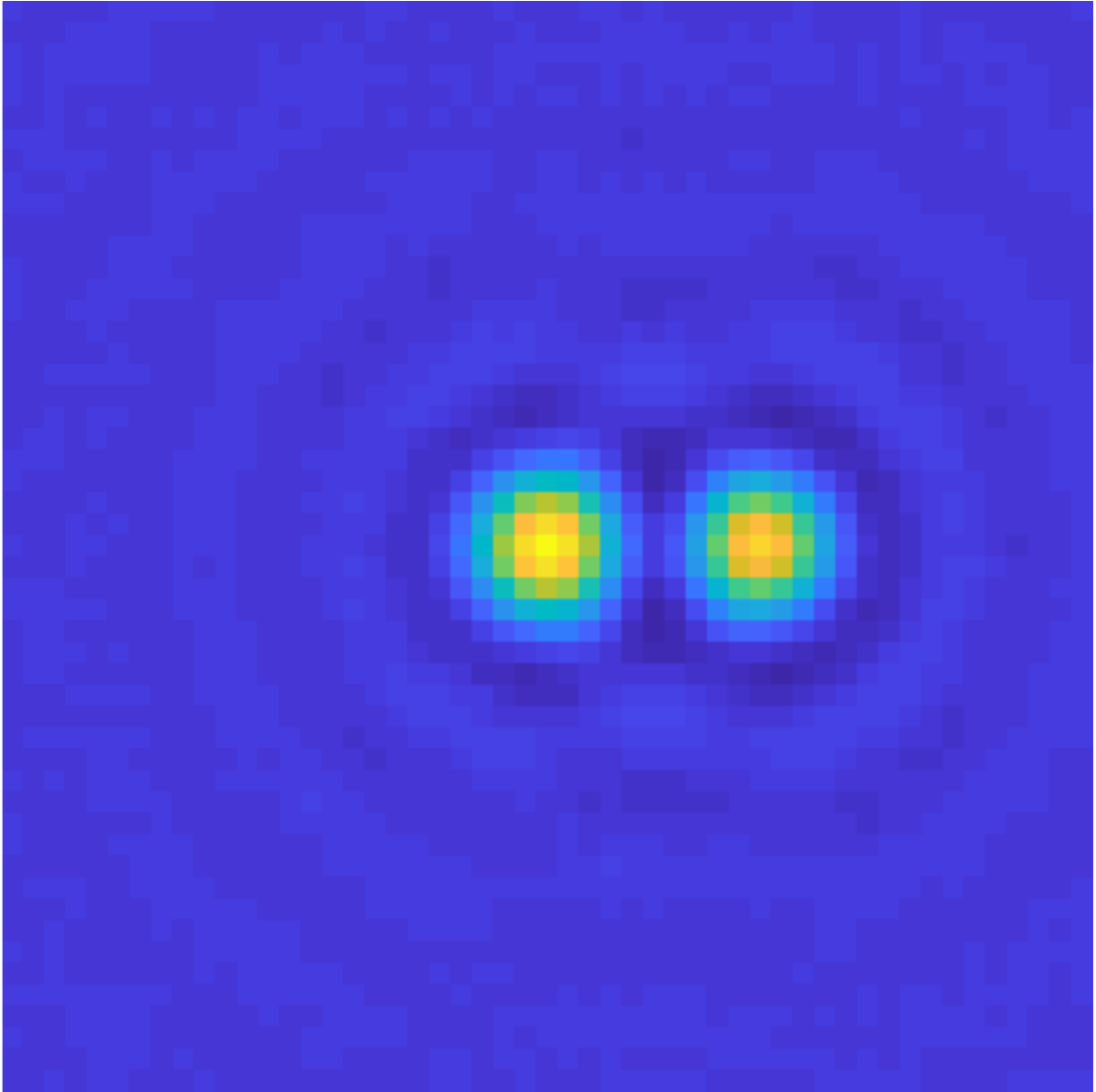}  &
\hspace{0.8cm}\includegraphics[height=1.3in]{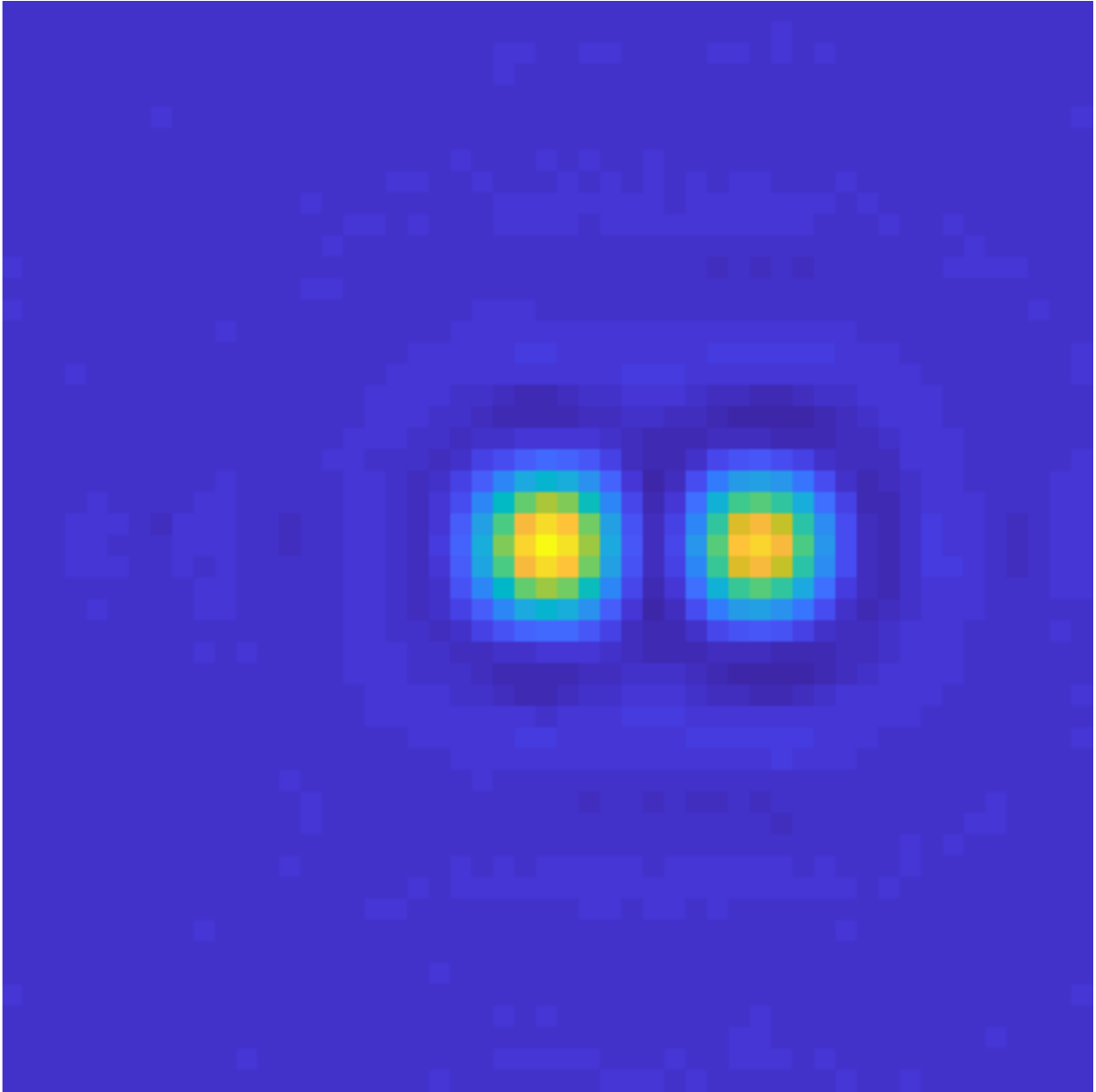} \\
\hspace{-0.8cm}{\small LR-FGMRES} &
\hspace{0.8cm}{\small 2nd cycle}\vspace{0.3cm}\\
\hspace{-0.8cm}\includegraphics[height=1.3in]{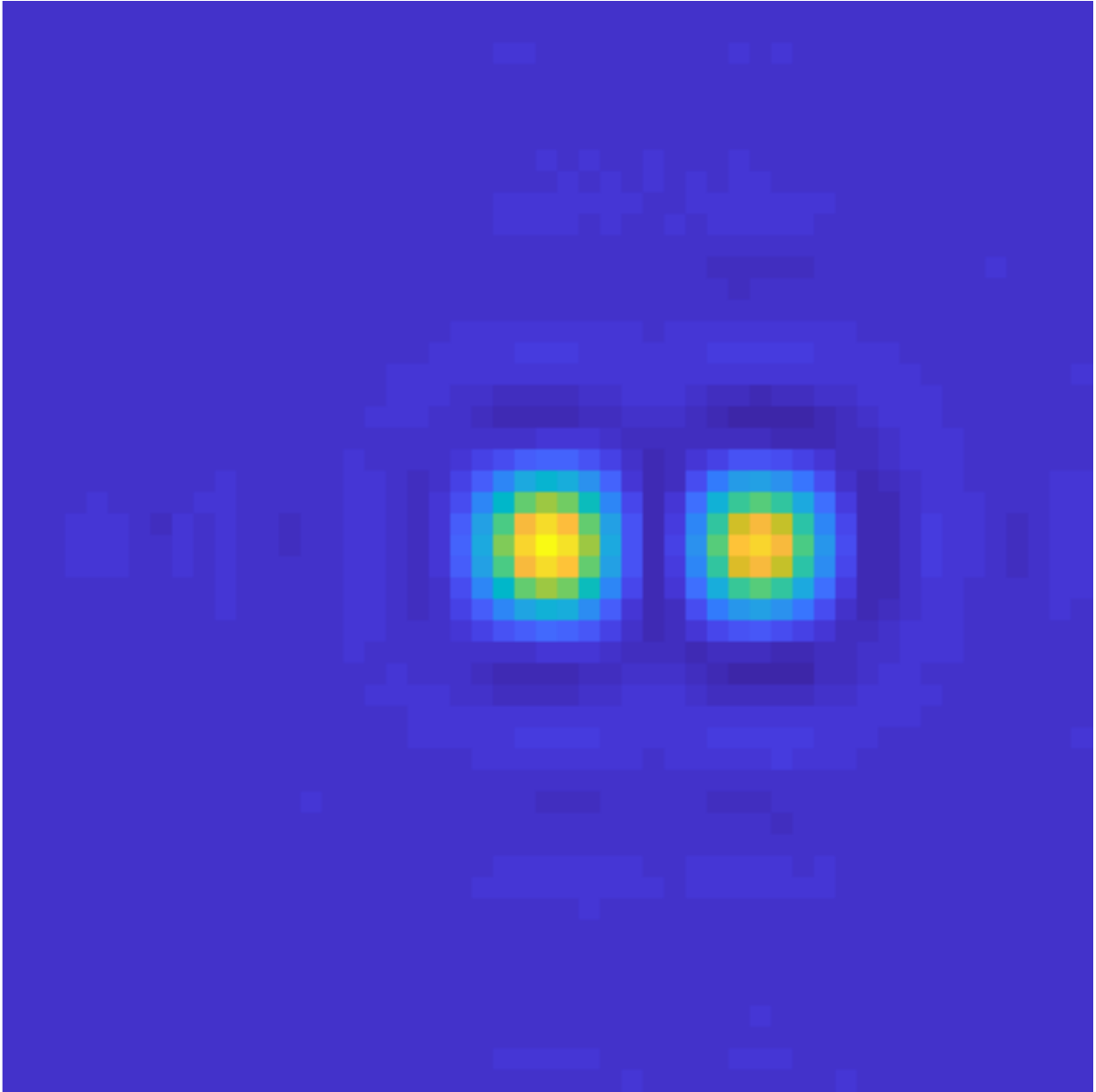}  &
\hspace{0.8cm}\includegraphics[height=1.3in]{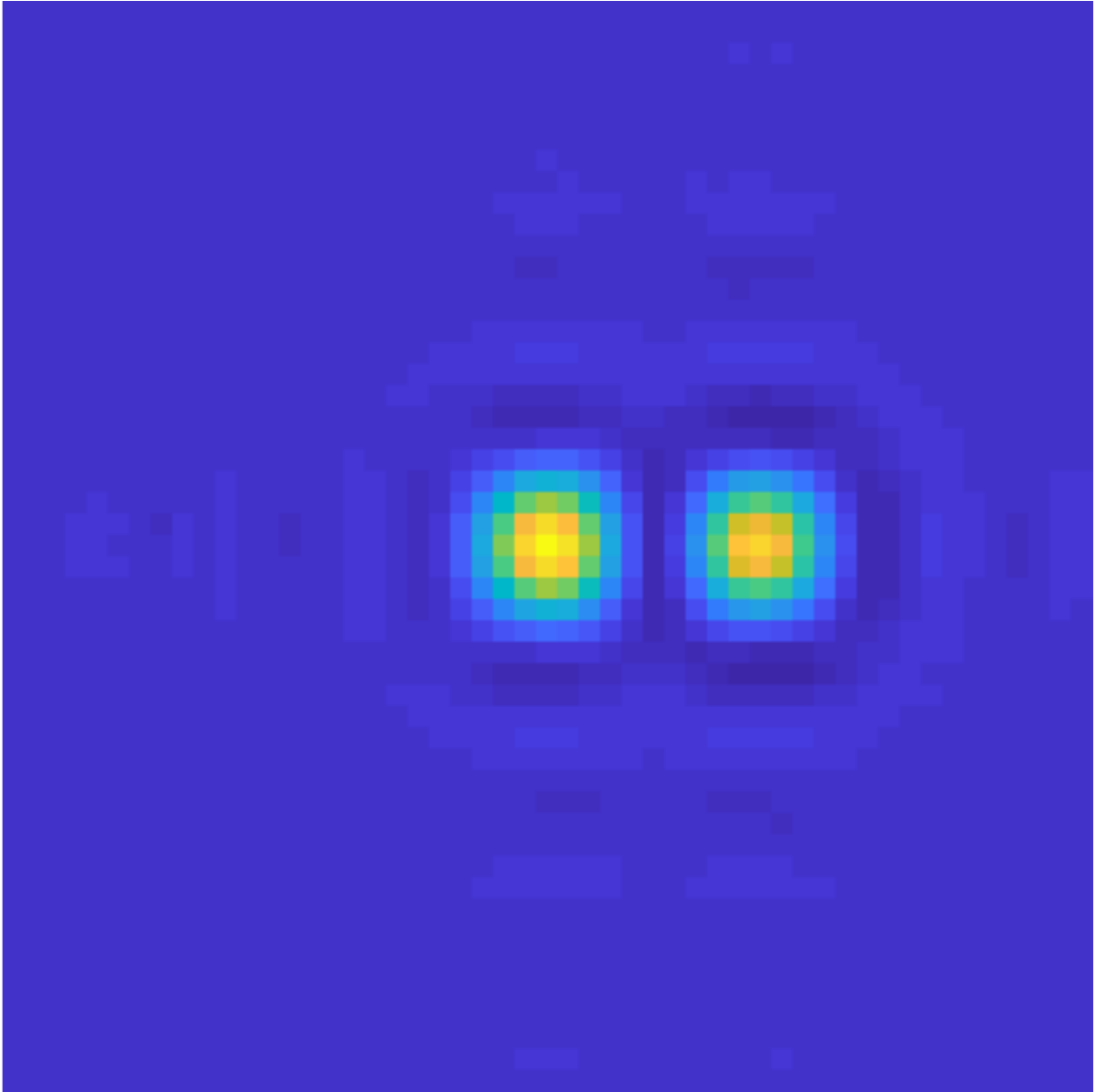} \\
\hspace{-0.8cm}{\small 3rd cycle} & 
\hspace{0.8cm}{\small 4th cycle}
\end{tabular}
\end{center}
\caption{\emph{Example 1}. Zoom-ins of the LR-FGMRES best solution, and the IRN-GMRES-NNR solutions at the end of each inner cycle.}
\label{fig:star_images_irncomp}
\end{figure}

More details can be spotted if we zoom into the central part ($51\times 51$ pixels) of the computed images, as shown in Figure \ref{fig:star_images_irncomp}: here the best LR-FGMRES reconstruction, as well as the IRN-GMRES-NNR reconstructions at the end of the 2nd, 3rd and 4th inner cycles are displayed. It is clear that the IRN reconstructions are improving over each outer cycle, 
and that even the solution at the end of the 2nd cycle is significantly better than the LR-FGMRES solution, which means that not all four outer iterations need to be run to achieve solutions of superior qualities (even if more outer iterations allow further improvement in the solution).
%
%

\begin{figure}[htbp]
\begin{center}
\begin{tabular}{ccc}
\includegraphics[height=1.0in]{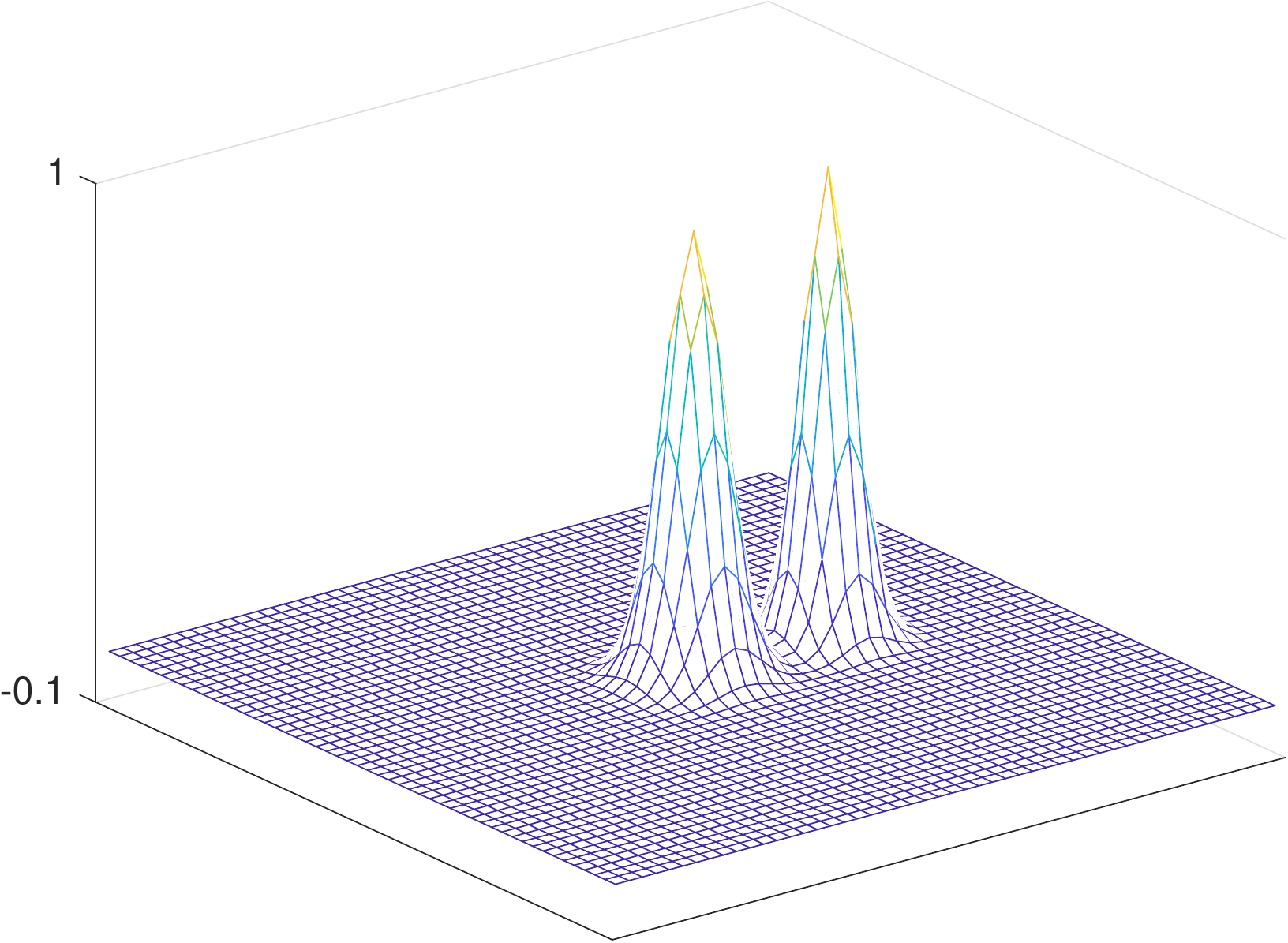}  &
\includegraphics[height=1.0in]{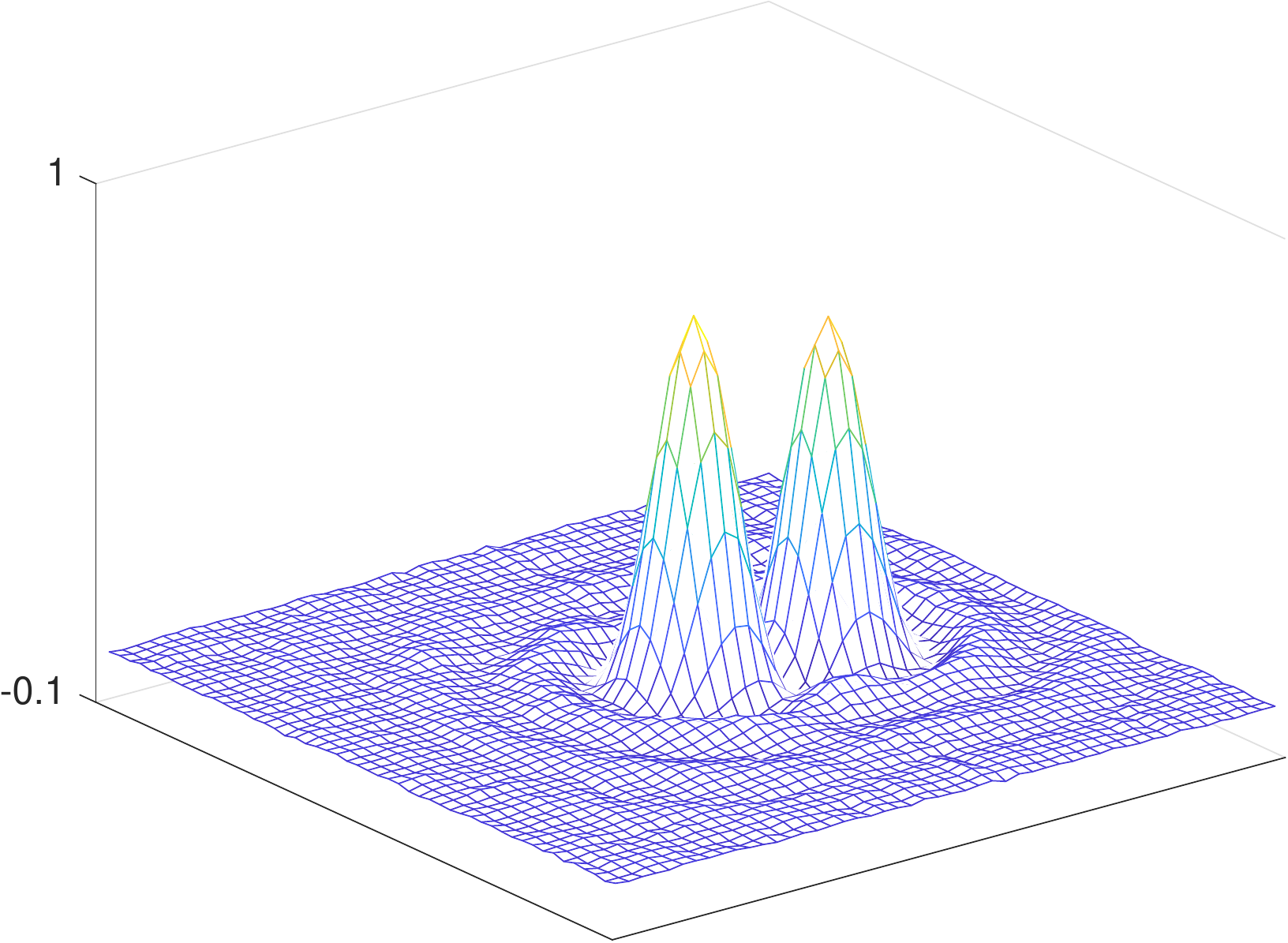} &
\includegraphics[height=1.0in]{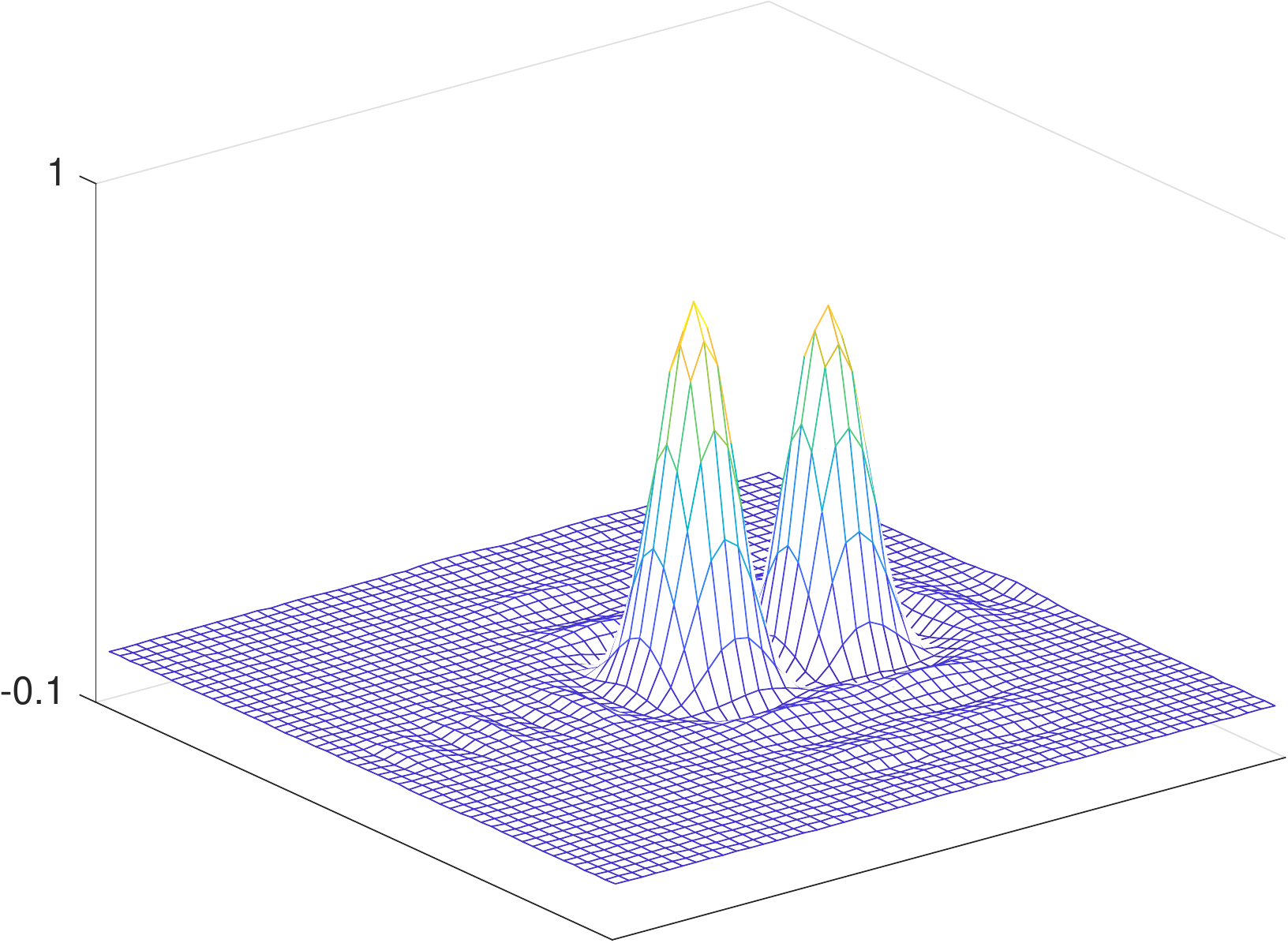} \\
{\small exact} & {\small LR-FGMRES} & {\small IRN-GMRES-NNR}\vspace{0.3cm}\\
\includegraphics[height=1.0in]{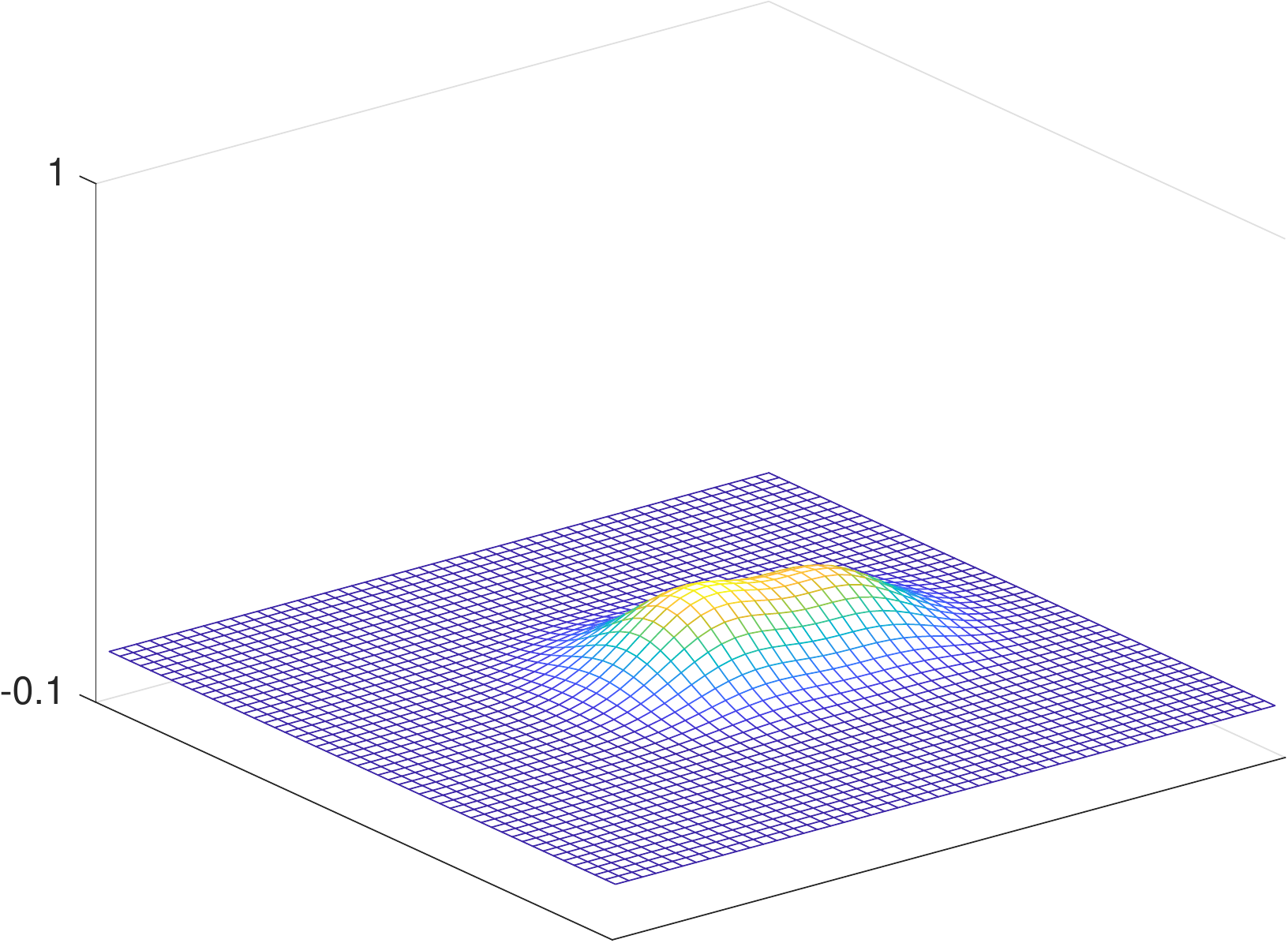}  &
\includegraphics[height=1.0in]{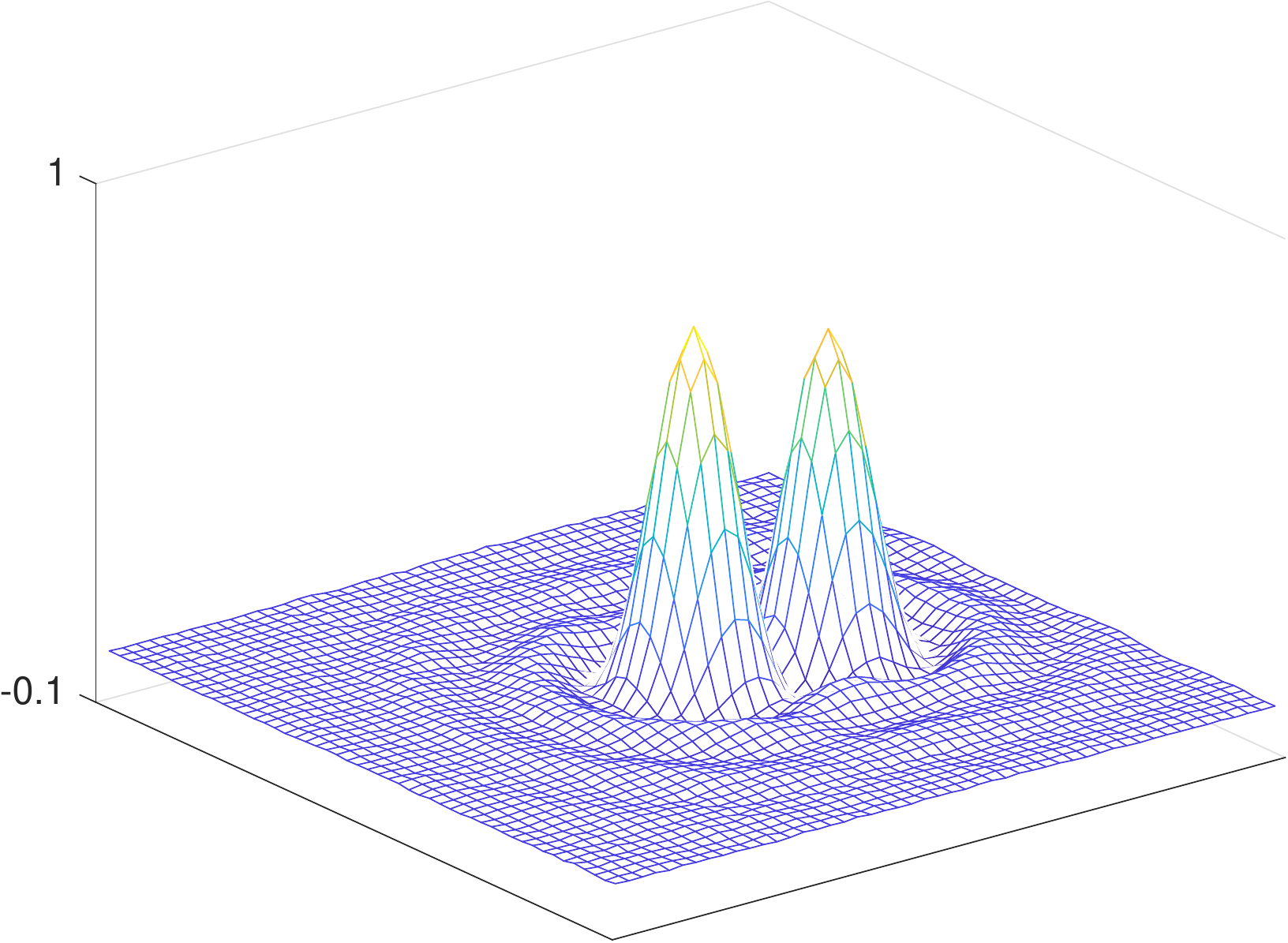} &
\includegraphics[height=1.0in]{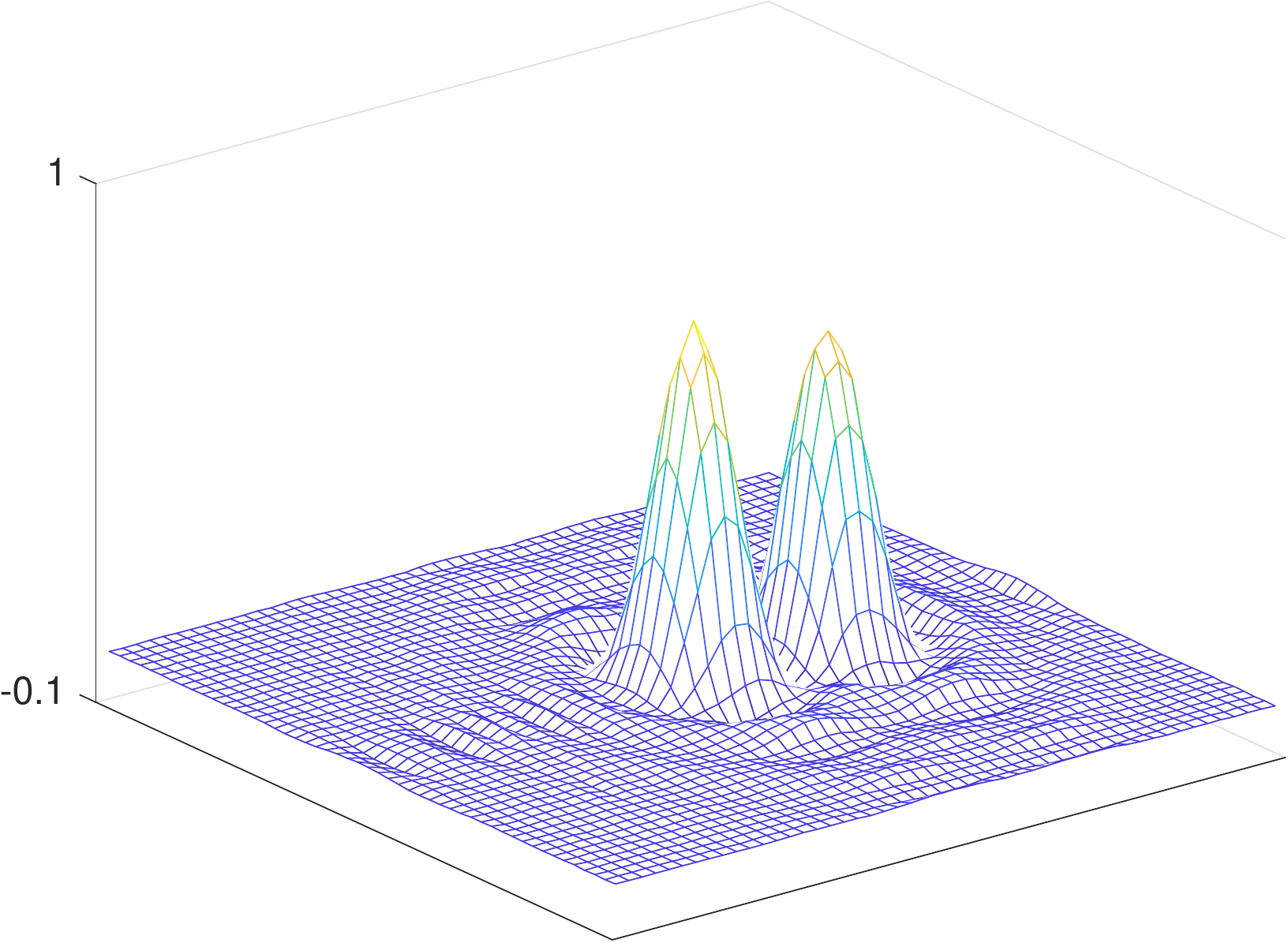} \\
{\small blurred \& noise} & {\small RS-LR-GMRES} & {\small FGMRES-NNR}
\end{tabular}
\end{center}
\caption{\emph{Example 1}. Zoomed-in surfaces of the exact solution and the available data, as well as the best reconstructions obtained by the new GMRES-based methods.}
\label{fig:star_surf}
\end{figure}

Figure \ref{fig:star_surf} displays surfaces plots of the central part ($51\times 51$ pixels) of the test problem data, as well as the best reconstructed images (for RS-LR-GMRES and FGMRES-NNR these are obtained at the 165th and the 63th (total) iterations, respectively). It can be seen that for all the solutions shown here, the reconstructed central two stars approximately have the same intensity, although they are somewhat less intense than in the exact image. These surface plots also confirm our earlier observation that IRN-GMRES-NNR does an exceptional job removing background noise. In addition, FGMRES-NNR also gives a good background reconstruction.

\begin{figure}[htbp]
\begin{center}
\begin{tabular}{cc}
\includegraphics[height=1.5in]{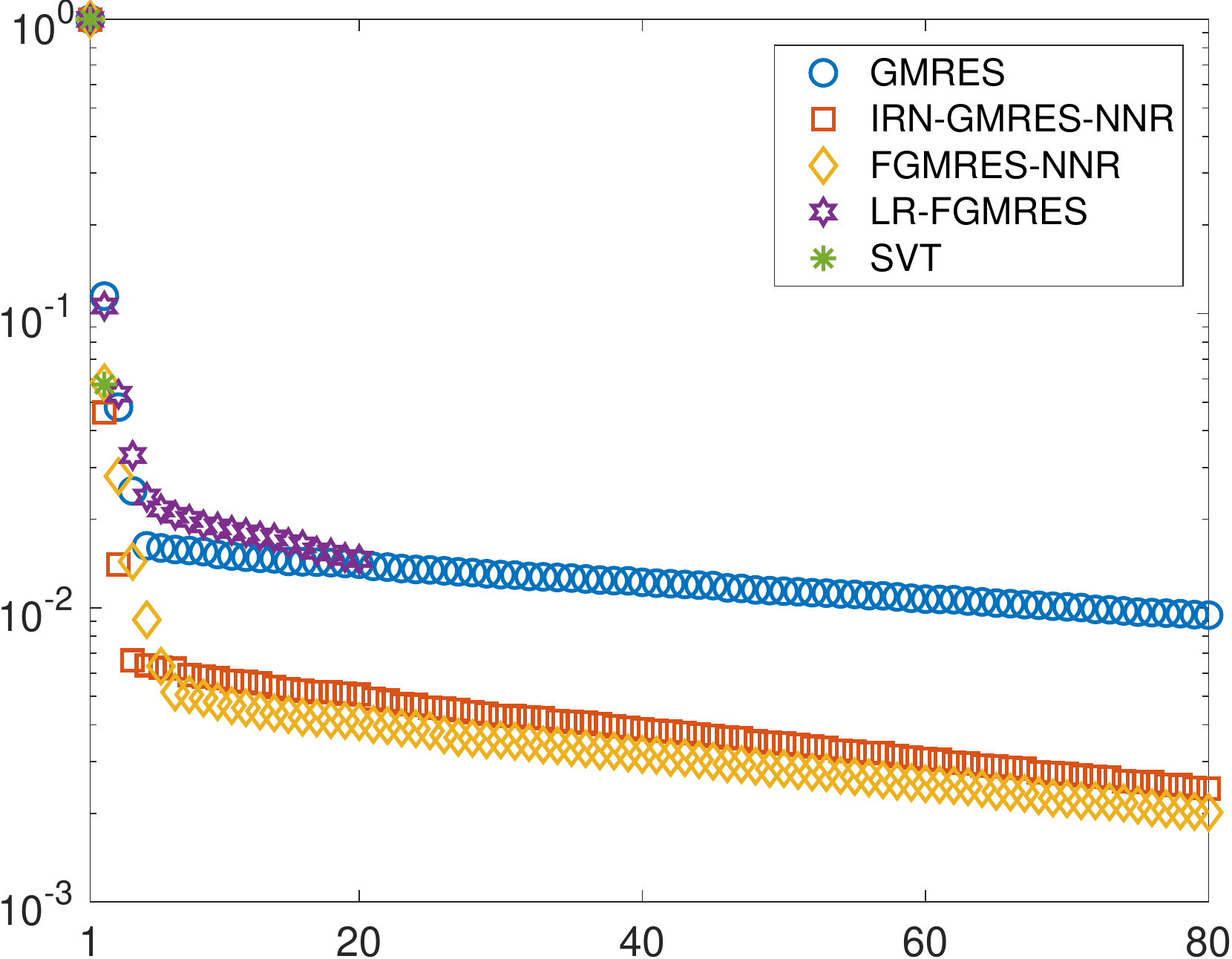}  &
\includegraphics[height=1.5in]{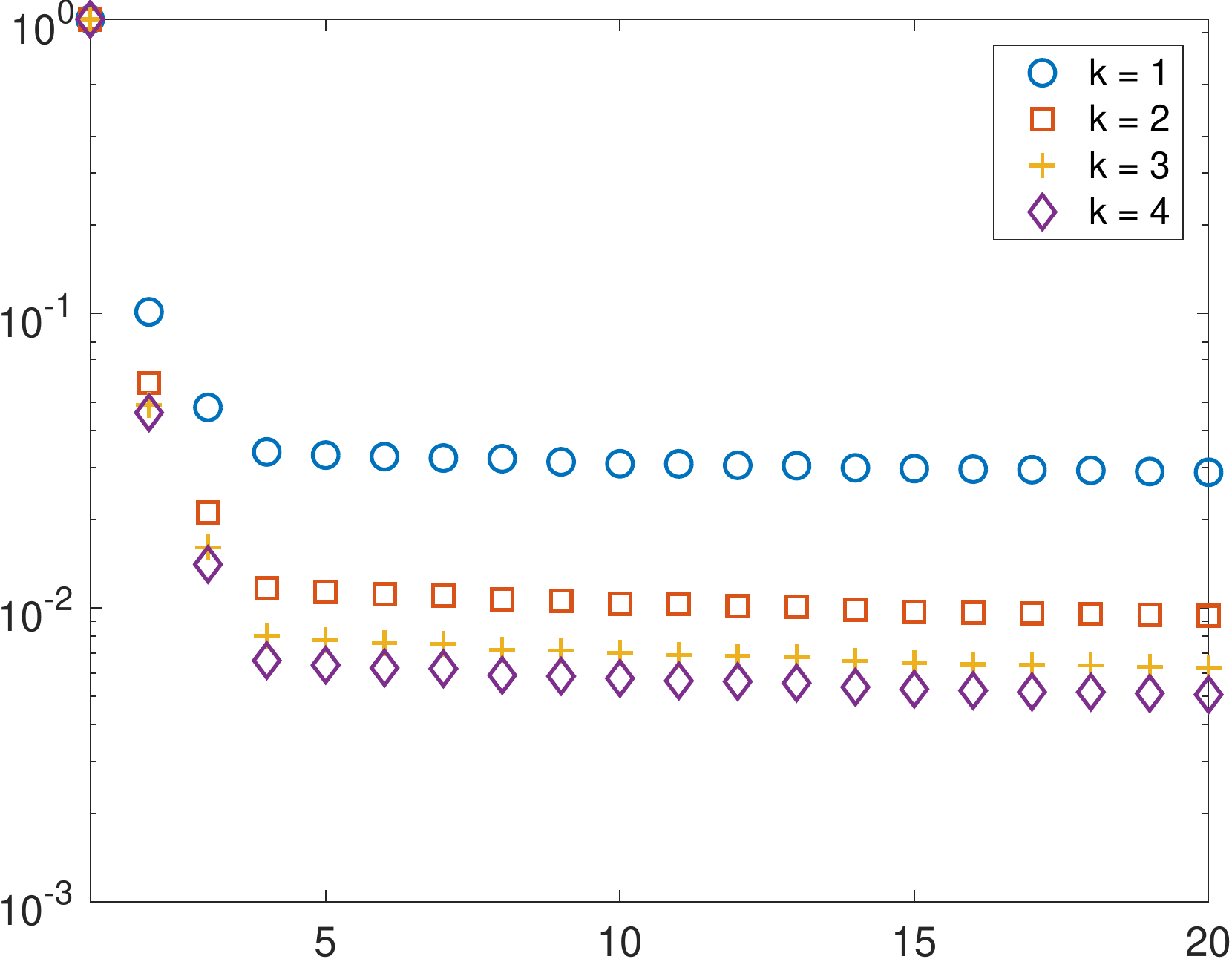} \\
{\small comparisons} & {\small IRN-GMRES-NNR} \\
\end{tabular}
\end{center}
\caption{\emph{Example 1}. Left frame: normalized singular values of the best solutions computed by each GMRES-based method. 
Right frame: evolution of the singular values of the solutions computed by IRN-GMRES-NNR at each outer iteration. Singular values less than $10^{-3}$ are omitted.}
\label{fig:star_sv}
\end{figure}
Finally, Figure \ref{fig:star_sv} displays the singular values of the best solutions obtained adopting different GMRES-based solvers, as well as the evolution of the singular values of the solution at the end of each inner IRN-GMRES-NNR cycle (matching the reconstructions displayed in Figure \ref{fig:star_images_irncomp}). The singular values are ``normalized" (i.e., divided by the largest one), and the graphs are cropped to focus on the relevant values. Looking at the displayed values, we can conclude that the solutions computed by all the low-rank solvers have indeed some low-rank properties, with very quickly-decaying large singular values followed by slowly-decaying smaller singular values. Compared to GMRES, the new FGMRES-NNR and IRN-GMRES-NNR methods give solutions that have a more pronounced low rank, as shown by the large gaps between the smaller singular values of the solutions computed by these methods. Regarding IRN-GMRES-NNR, the evolution of the singular values stabilizes as we move toward later outer iterations, which validates the stopping criterion 
proposed in Section \ref{sec:implementation}.

\begin{figure}[H]
\begin{center}
\begin{tabular}{cccc}
\hspace{-0.5cm}\includegraphics[height=1.3in]{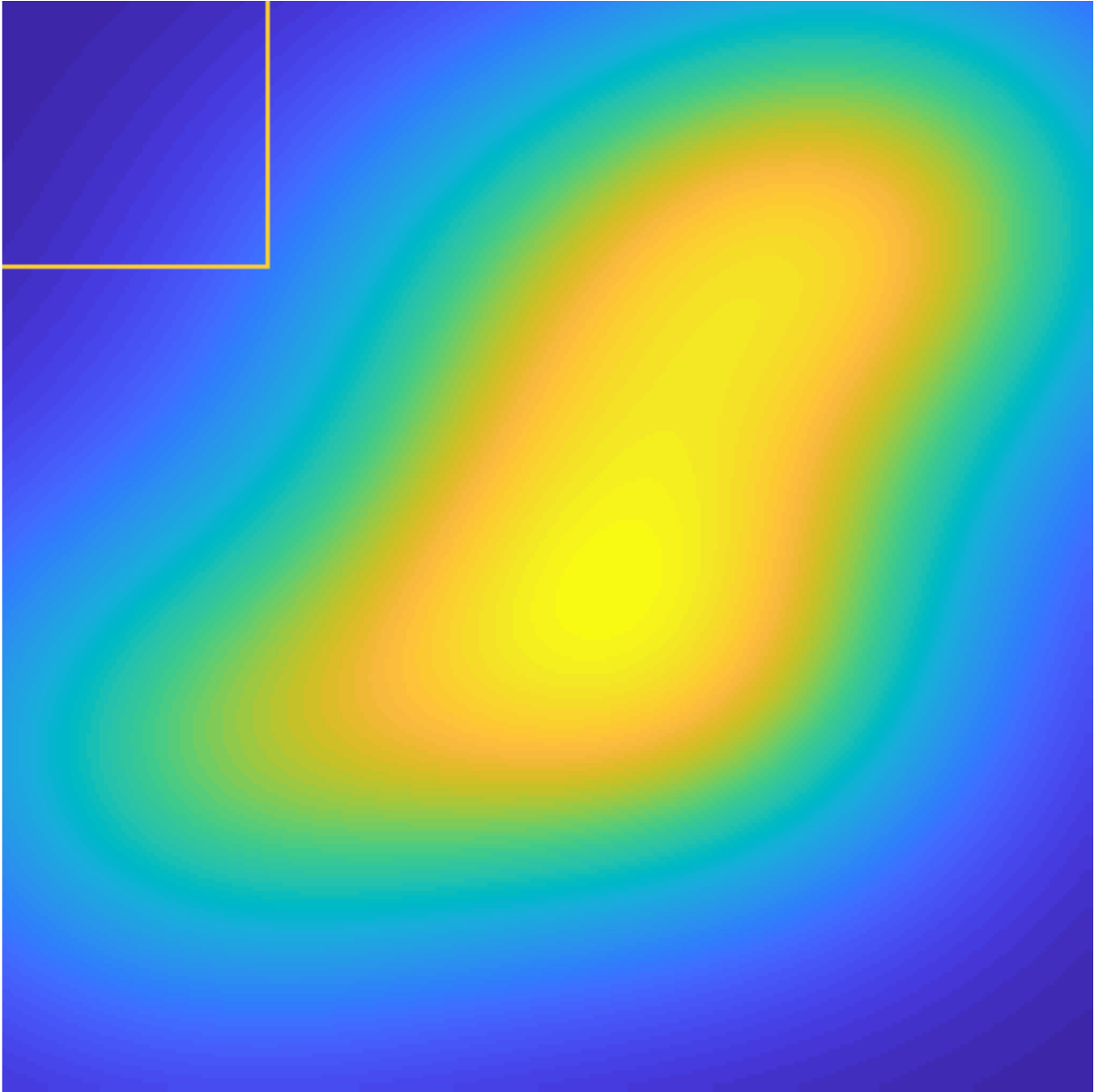}  &
\hspace{0.5cm}\includegraphics[height=1.3in]{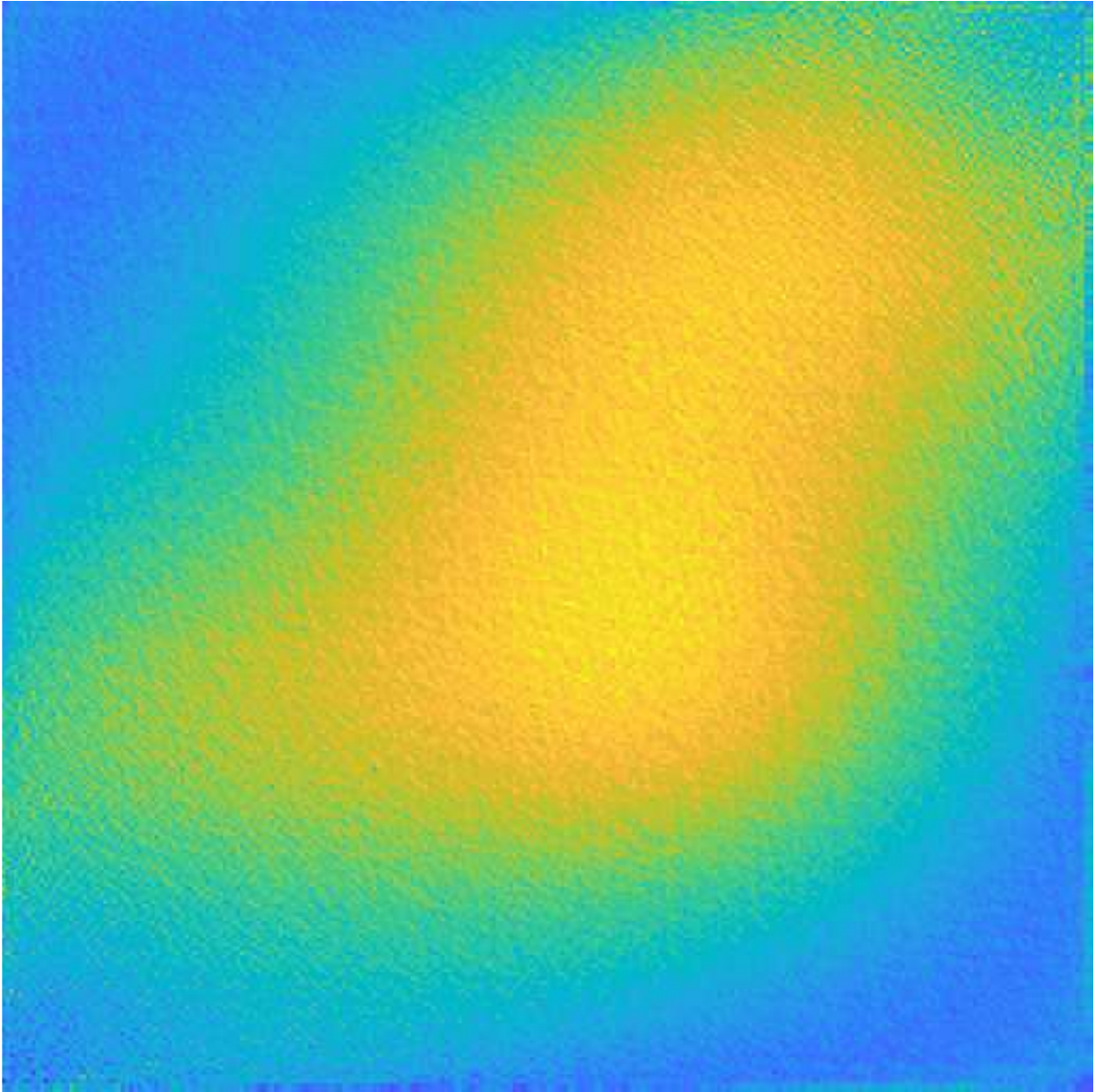}\\
\hspace{-0.5cm}{\small exact} & 
\hspace{0.5cm}{\small LSQR} \vspace{0.3cm}\\
\hspace{-0.5cm}\includegraphics[height=1.3in]{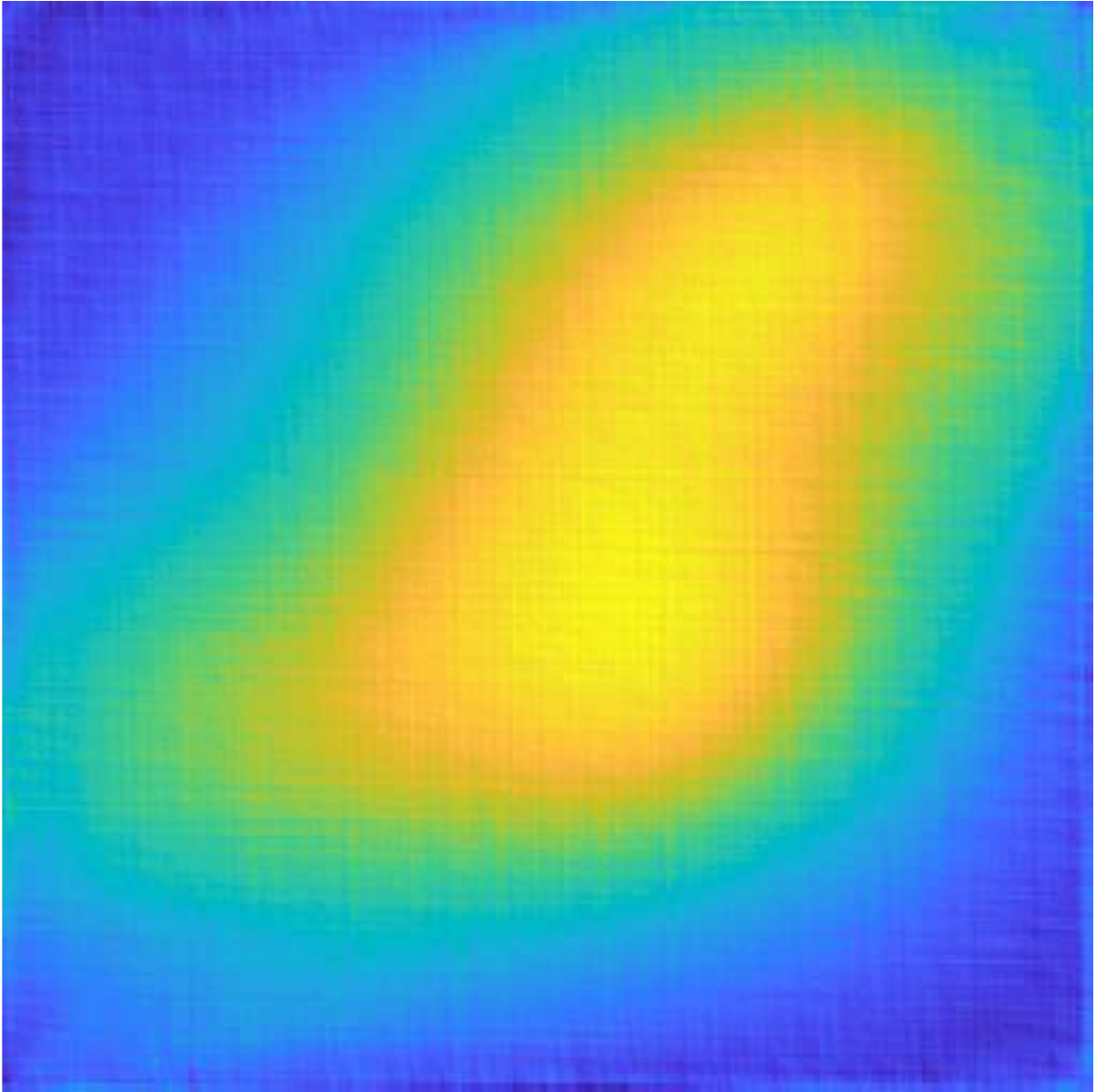}&
\hspace{0.5cm}\includegraphics[height=1.3in]{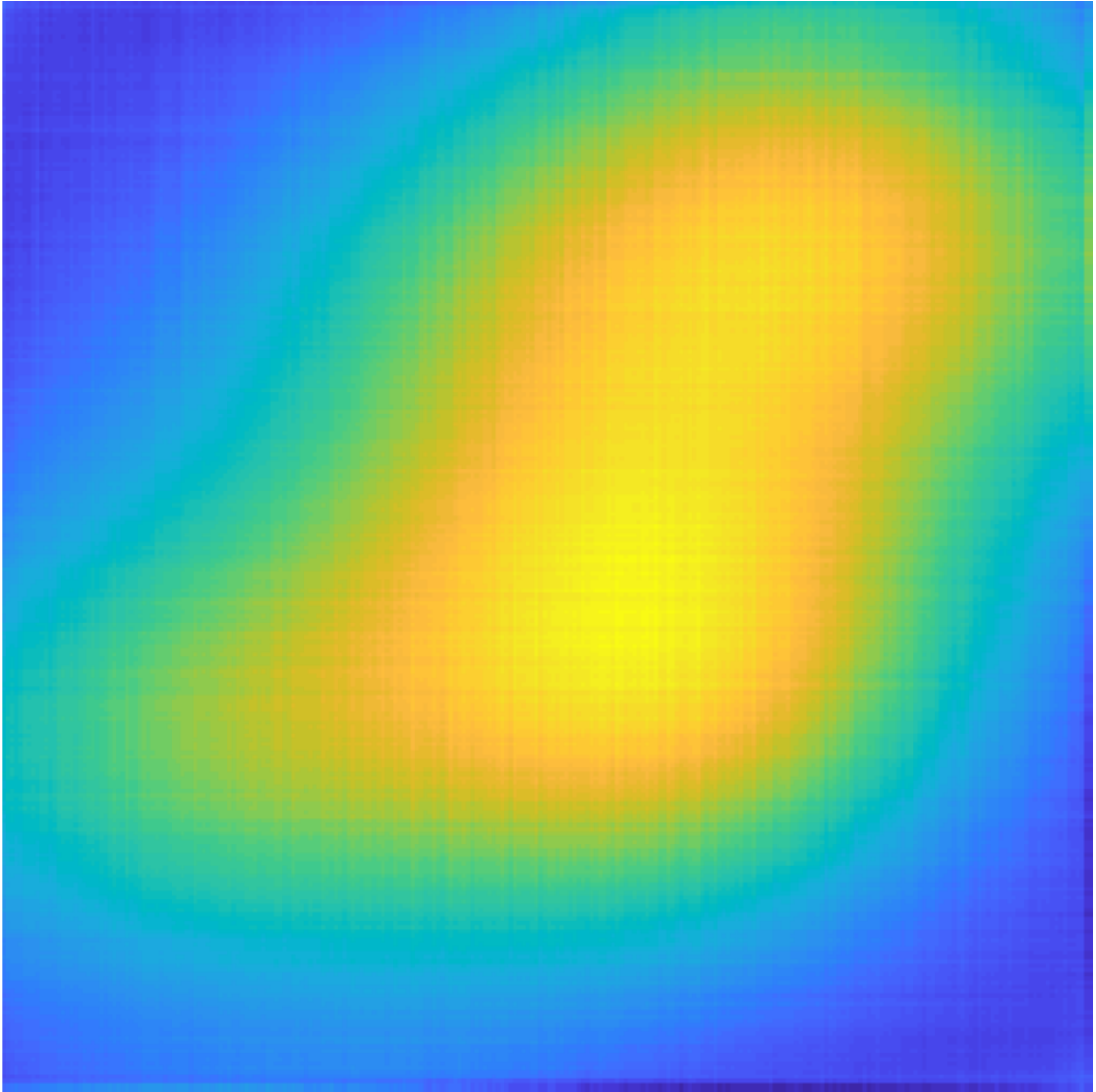} \\
\hspace{-0.5cm}{\small FLSQR-NNR$p$(v)} & 
\hspace{0.5cm}{\small IRN-LSQR-NNR}
\end{tabular}
\end{center}
\caption{{\em Example 2}. Exact phantom and best reconstructions obtained by different solvers.}
\label{fig:tomo_sol}
\end{figure}

\paragraph{Example 2:  Limited angle parallel-ray tomography}
We consider a computed tomography (CT) test problem, modeling an undersampled X-ray scan with parallel beam geometry. This is a so called ``limited angle" CT reconstruction problem, where the viewing angles for the object span less than 180 degrees. A smooth and rank-4 phantom is considered, as shown in the leftmost frame of Figure \ref{fig:tomo_sol} (note that the yellow straight lines in the northwestern corner do not belong to the phantom; they are shown for later purposes). Gaussian white noise of level $10^{-2}$ is added to the data. 
%
The coefficient matrix $\mxA$ has size $32942 \times 65536$. Because of this, among the new solvers, only LR-FLSQR, FLSQR-NNR$p$, FLSQR-NNR$p$(v), and IRN-LSQR-NNR$p$ will be tested, against their standard counterpart LSQR. Recall that FLSQR-NNR$p$(v) is the FLSQR-NNR$p$ variant that defines the preconditioners using the basis vectors of the solution subspace. The hybrid strategy is not used here, meaning that we set $\widehat{\lambda}=0$ for all methods. For this test problem, we consider both the values $p=1$ and $p=0.75$ (recall that, when $p=1$, we omit $p$ from the notation). 
%
The results obtained running the available low-rank solvers SVT and RS-LR-GMRES are shown, too. Note that RS-LR-GMRES only works for square matrices $\mxA$, hence this solver is tested on the normal equations $\mxA^T\mxA\bfx = \mxA^T\bfb$, which is not the problem solved by the other methods (therefore this comparison may not be completely fair). Parameters for SVT are chosen to be: step size $\delta_k = \delta = 8\times 10^{-5}$ and threshold $\tau =100$. RS-LR-GMRES is set to restart every 20 iterations. The truncation rank is 10 for both basis vectors and solutions, and for both the LR-FLSQR and the RS-LR-GMRES methods. The maximum number of iterations is 100 for all methods.

\begin{figure}[H]
\begin{center}
\begin{tabular}{c}
\includegraphics[height=1.8in]{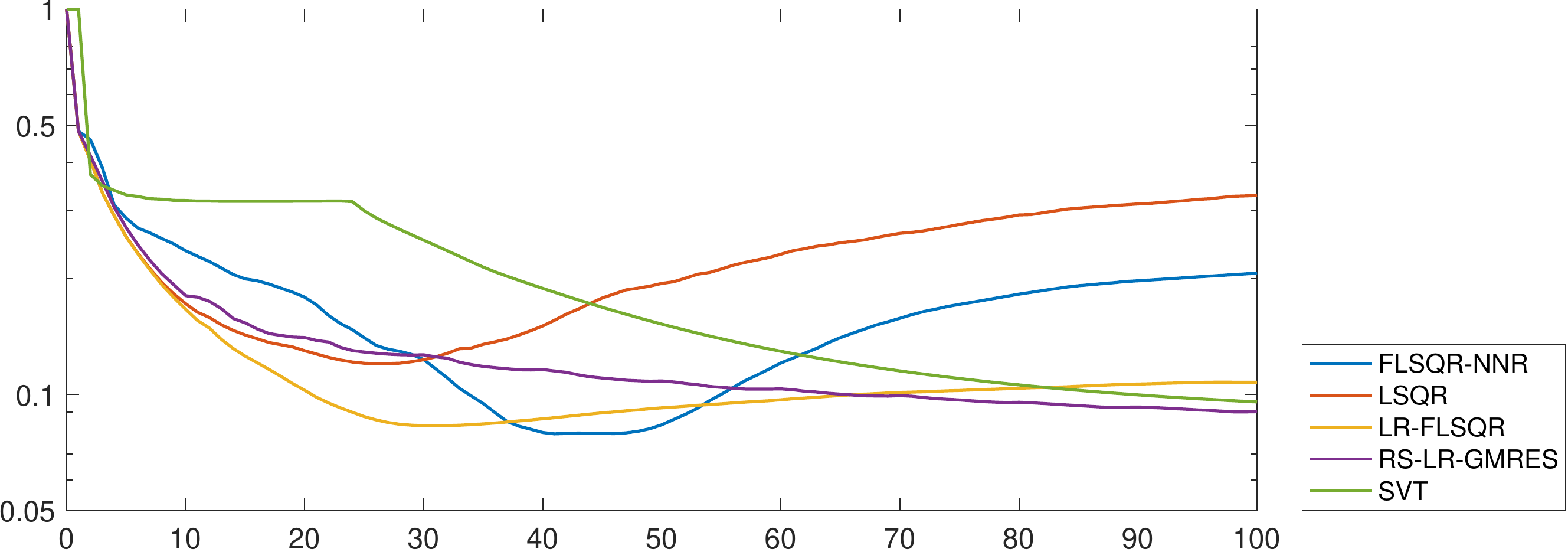} \\
\includegraphics[height=1.8in]{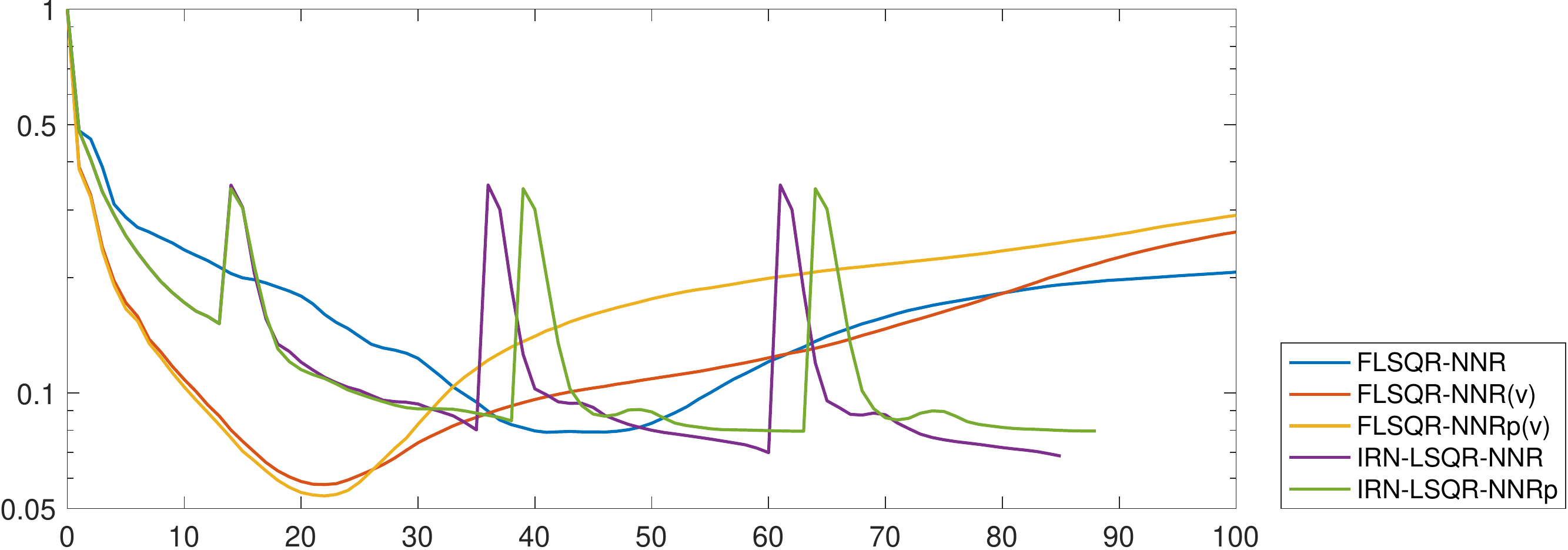}
\end{tabular}
\end{center}
\caption{{\em Example 2}. Relative errors vs. number of iterations for different solvers. Upper frame: some of the new solvers are compared to the already available solvers. Lower frame: comparisons of different instances of the new solvers (here $p=0.75$).}
\label{fig:tomo_error}
\end{figure}
\begin{figure}[H]
\begin{center}
\begin{tabular}{ccc}
\includegraphics[height=1.0in]{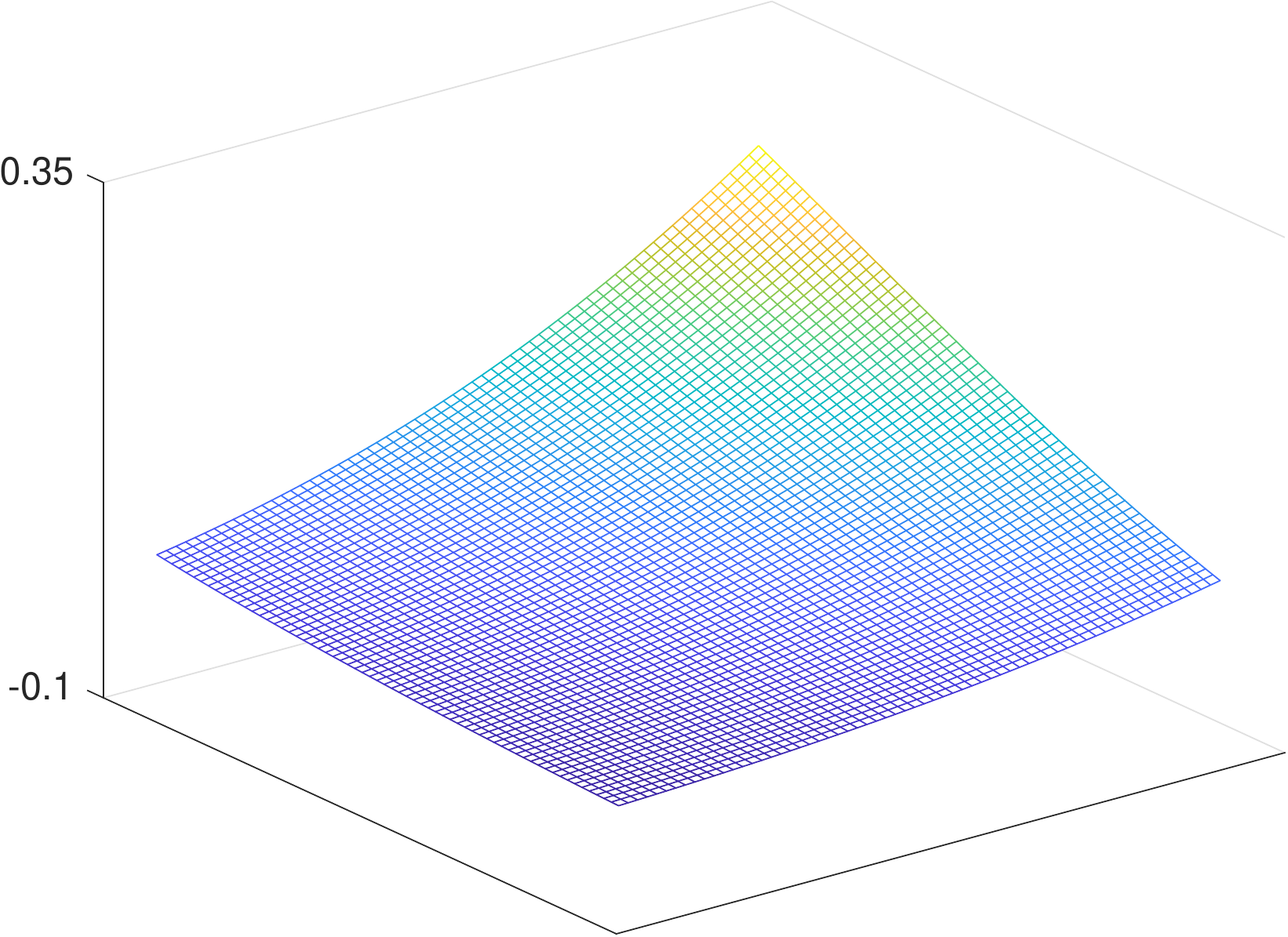}  &
\includegraphics[height=1.0in]{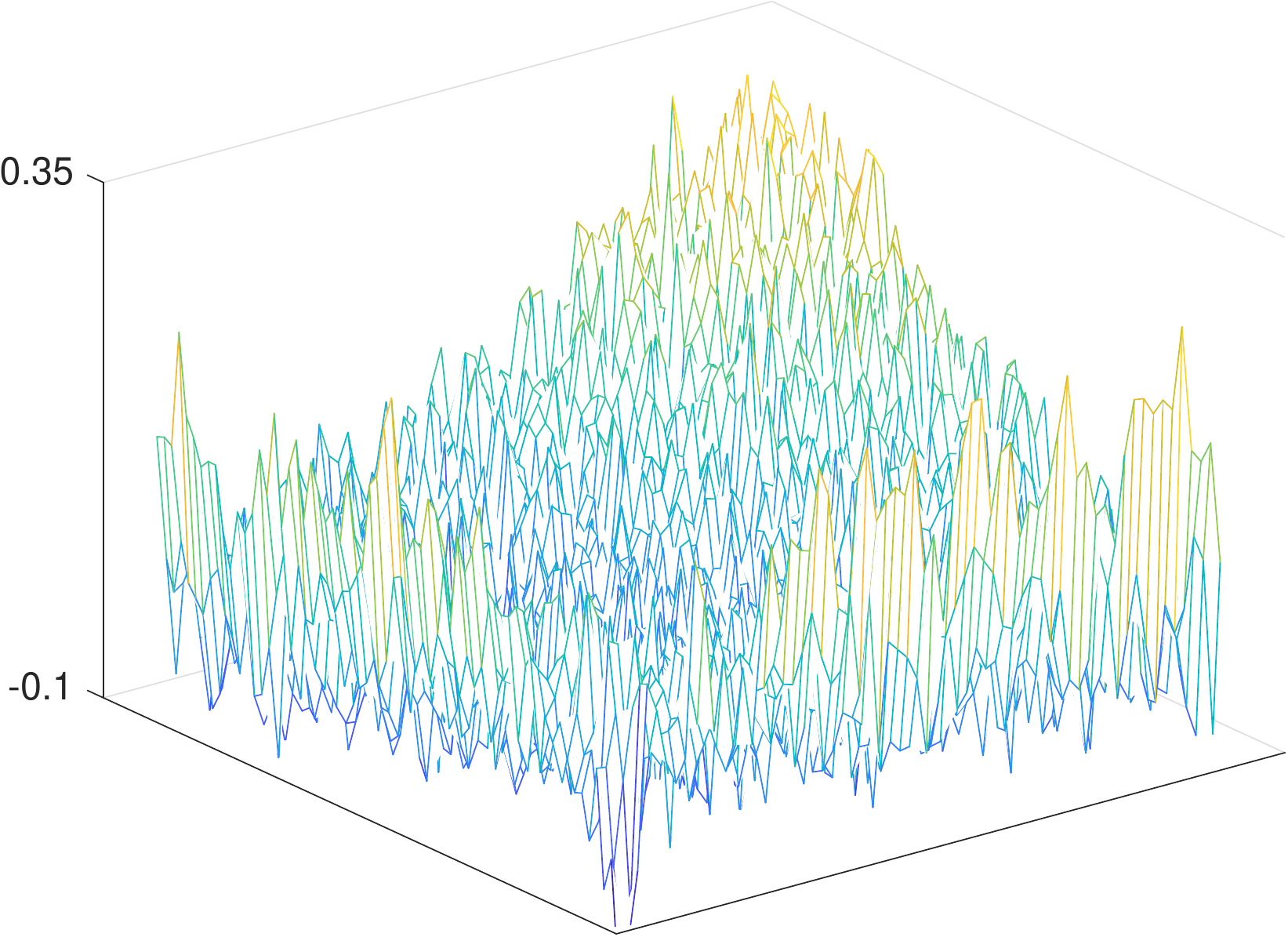}&
\includegraphics[height=1.0in]{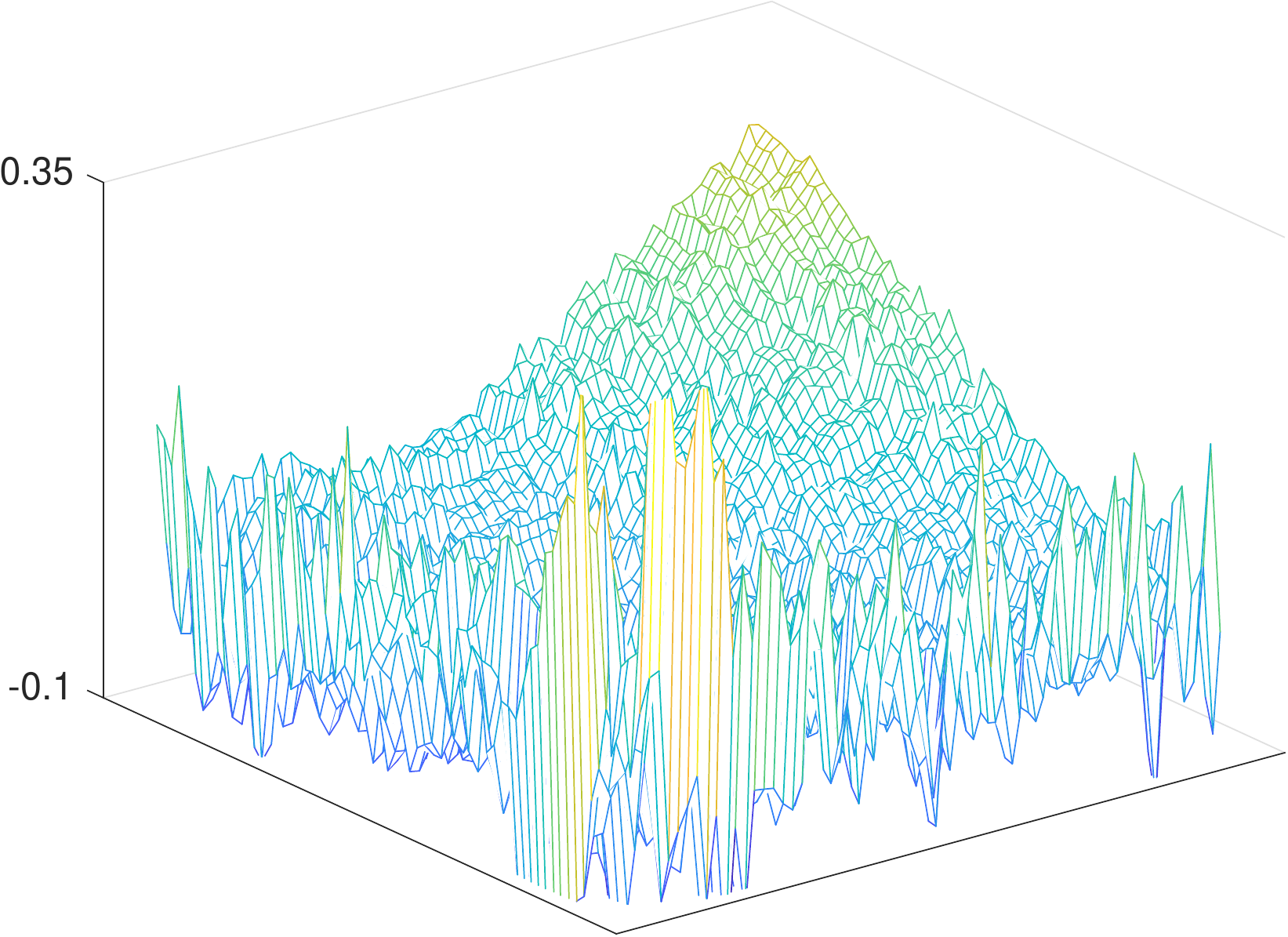}\\
{\small exact} & {\small LSQR} & {\small LR-FLSQR}\vspace{0.3cm}\\
\end{tabular}
\begin{tabular}{ccc}
\includegraphics[height=1.0in]{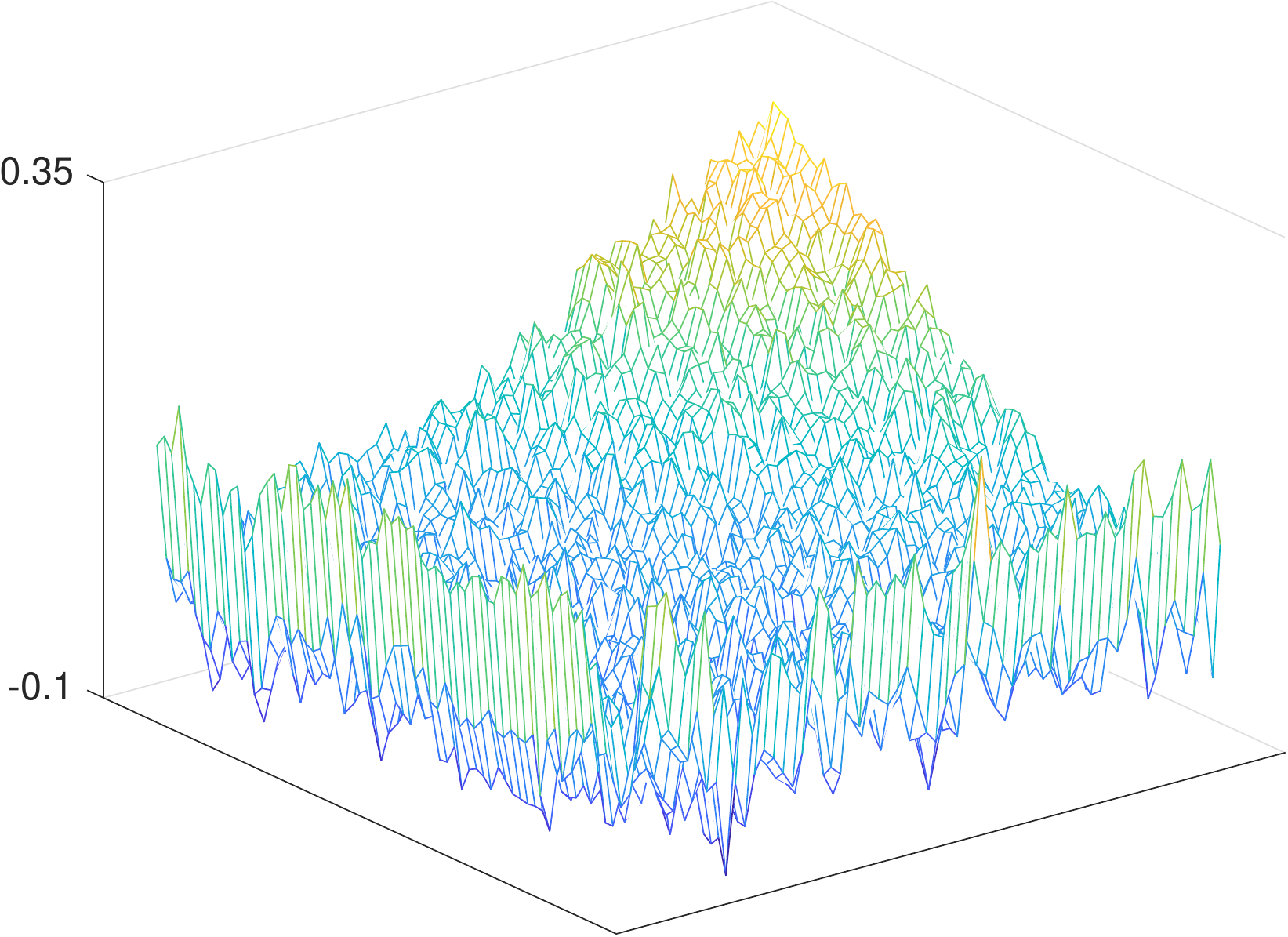} &
\includegraphics[height=1.0in]{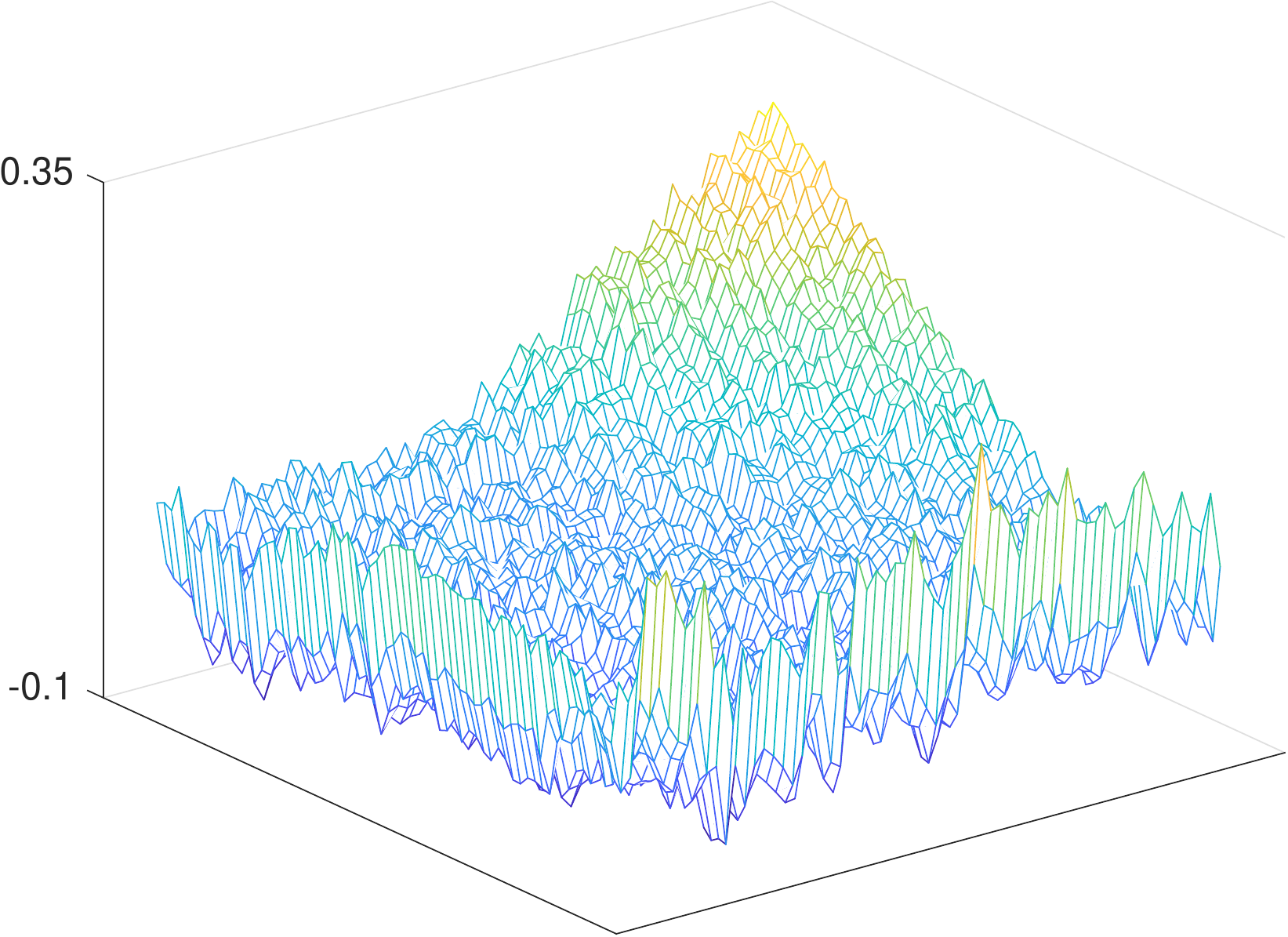} &
\includegraphics[height=1.0in]{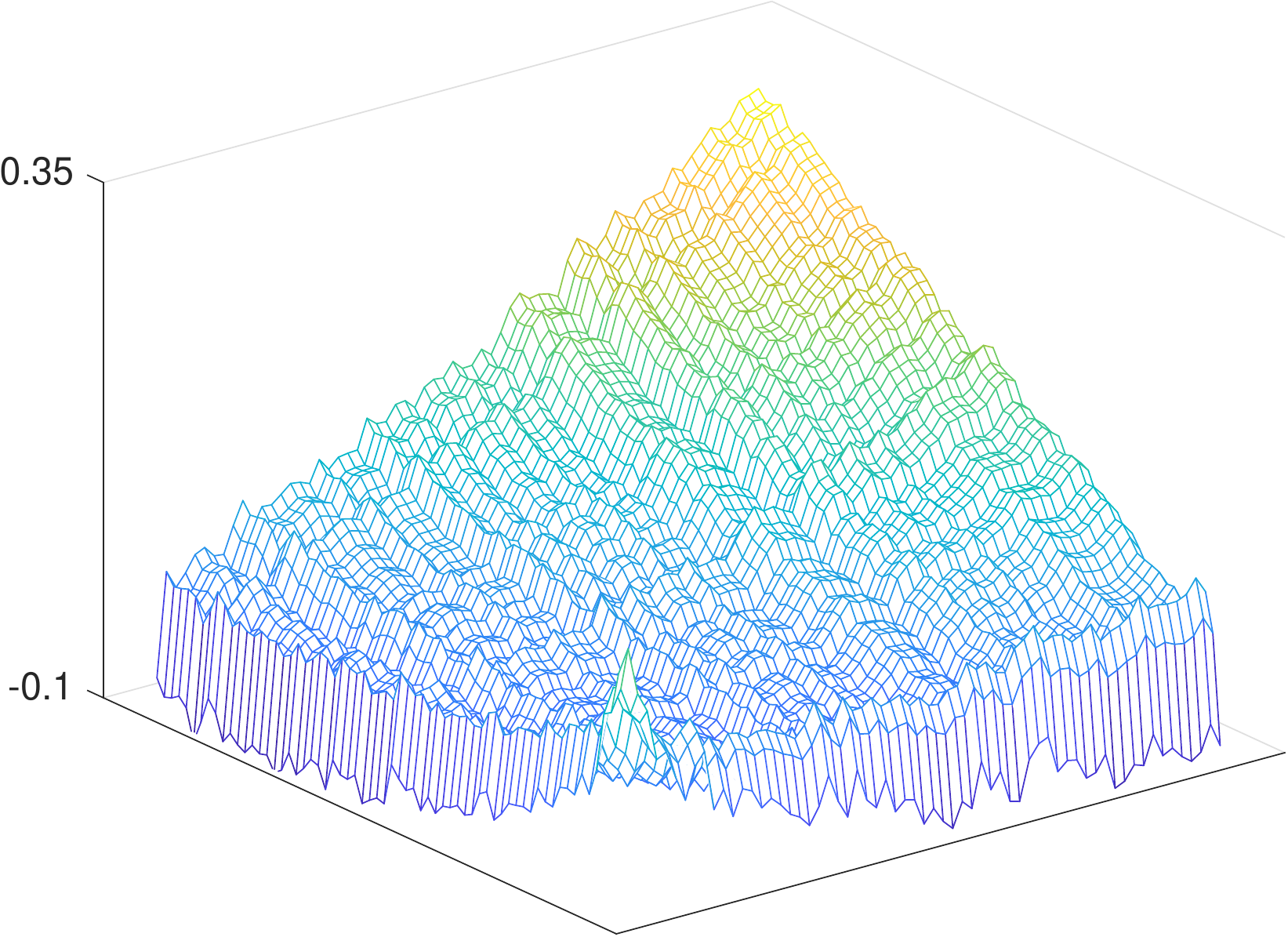}\\
 {\small FLSQR-NNR(v)}&{\small FLSQR-NNRp(v)}&{\small IRN-LSQR-NNR}\\
\end{tabular}
\end{center}
\caption{{\em Example 2}. Surface plots of the northwestern corner of the exact phantom (highlighted in Figure \ref{fig:tomo_sol}) and the best reconstructed phantoms computed by different solvers.}
\label{fig:tomo_surf}
\end{figure}

Figure \ref{fig:tomo_error} displays the history of the relative errors for LSQR, LR-LSQR, FLSQR-NNR$p$, FLSQR-NNR$p$(v), and IRN-LSQR-NNR$p$, for $p=1$ and $p=0.75$. 
Figure \ref{fig:tomo_sol} displays the exact phantom together with the best reconstructions obtained by LSQR, FLSQR-NNR$p$(v), and IRN-LSQR-NNR. Figure \ref{fig:tomo_surf} displays surface plots of the northwestern corner of the exact and reconstructed phantoms ($64\times 64$ pixels, as highlighted in the leftmost frame of Figure \ref{fig:tomo_sol}).

Looking at relative errors in Figure \ref{fig:tomo_error}, it is obvious that the winners are the FLSQR-NNR$p$(v) methods, with both $p=1$ and $p = 0.75$: they give the lowest relative errors, and the fastest semi-convergences. For this test problem, using a value of $p<1$ lowers the relative error of FLSQR-NNR$p$(v); however, the same does not hold for IRN-LSQR-NNR$p$. Therefore we can conclude that the the choice of $p$ is problem and solver dependent, and using $p<1$ does not necessarily improve the quality of the solution. We regard $p=1$ as a safe choice for this parameter. Although both the FLSQR-NNR$p$(v) methods with $p=1$ and $p=0.75$ perform well, the latter is able to further reduce the noise in the reconstructed solution, especially on the boundary.

Looking at all the displayed results, the advantages of our new FLSQR-NNR$p$(v) and IRN-LSQR-NNR methods are evident. Namely, they produce smooth solutions that preserve the original concave shape of the exact phantom, and they retain similar intensities of pixels at the same locations of the exact phantom (although the LR-FLSQR solution is smooth within the boundary, it fails to reconstruct intensity at the high point). Differences between FLSQR-NNR$p$(v) and IRN-LSQR-NNR reconstructions are clear, too: while both are smooth, the IRN-LSQR-NNR reconstruction has a less concave shape compared to that of FLSQR-NNR$p$(v), but a smoother boundary.

\paragraph{Example 3: Inpainting} 
We consider two different inpainting test problems. Inpainting is the process of restoring images that have missing or deteriorated parts. These images are likely to have quite a few lost pixels, either in the form of salt and pepper noise, or missing patches with regular or irregular shapes. The two examples considered here are of different nature: the first one has less structured and more randomly distributed missing patches, while the second one has more structured and regularly shaped missing parts. 
The corrupted images (shown in top-middle frames of Figures \ref{fig:recin} and \ref{fig:recinp}) are constructed by first applying a blur operator, and then superimposing the undersampling pattern. We follow this particular order of first blurring and then taking out pixels to simulate the real process of photo-taking. For both these test problems, white noise of level $10^{-2}$ is added to the data, and we consider purely iterative methods (i.e., $\hl=0$). We always take $p=1$, and we run 100 iterations of all the methods.

\begin{figure}[H]
\begin{center}
\includegraphics[height=1.8in]{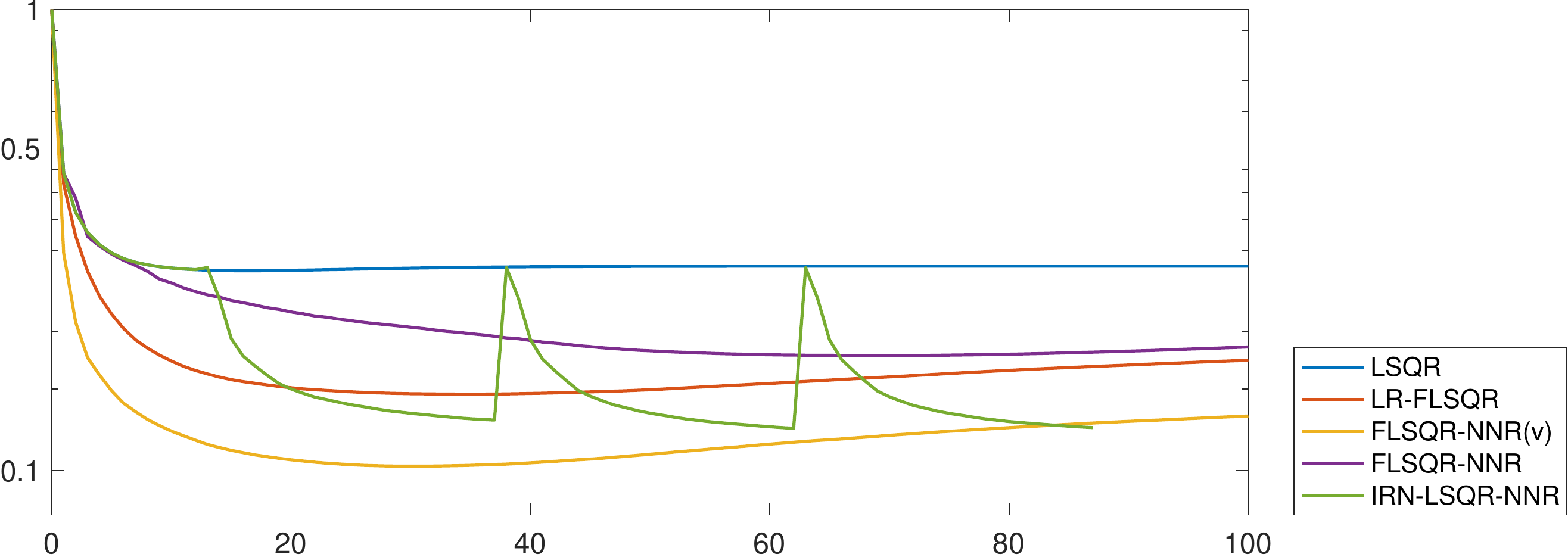}
\end{center}
\caption{{\em Example 3} (\texttt{house}). Relative errors vs. number of iterations for different solvers.}
\label{fig:relerrin}
\end{figure}

Firstly, we consider a test problem where 41.8\% of the pixels are missing (following some random and not very regular patterns). The exact image is commonly known as \texttt{house} test image: the version considered here has been truncated to rank 50, and has a total number of 65536 ($256\times 256$) pixels, out of which 27395 are non-zero. 
Correspondingly, the forward operator $\mxA$ is of size $27395\times 65536$, so we have an underdetermined linear system: $\mxA$ is obtained by first applying a shaking blur, and by then undersampling the blurred image. This can be easily coded within the \emph{IR Tools} framework.

\begin{figure}[H]
\begin{center}
\begin{tabular}{ccc}
\hspace{-1cm}\includegraphics[height=1.0in]{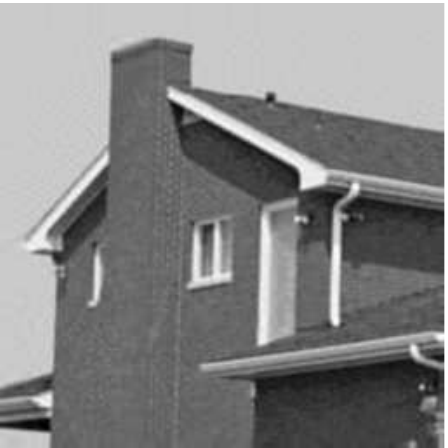}  &
\hspace{0.5cm}\includegraphics[height=1.0in]{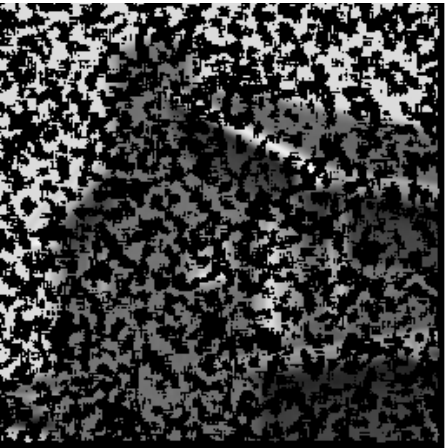} &
\hspace{0.5cm}\includegraphics[height=1.0in]{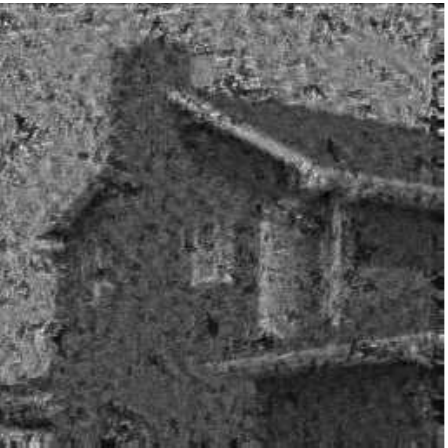} \\
\hspace{-1cm}{\small exact} & 
\hspace{0.5cm}{\small corrupted}&
\hspace{0.5cm}{\small LSQR}\vspace{0.3cm}
\end{tabular}
\begin{tabular}{cccc}
\hspace{-2.0cm}\includegraphics[height=1.0in]{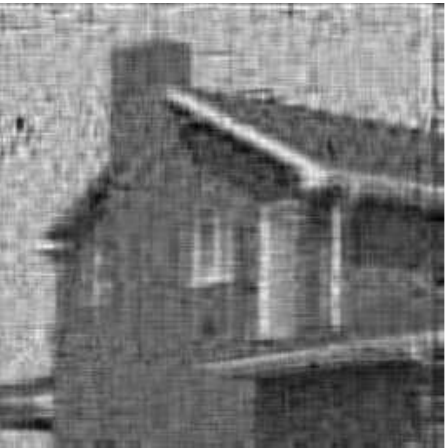} &
\hspace{0.8cm}\includegraphics[height=1.0in]{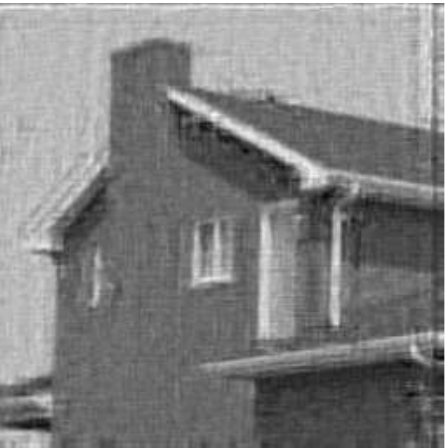}&
\hspace{0.8cm}\includegraphics[height=1.0in]{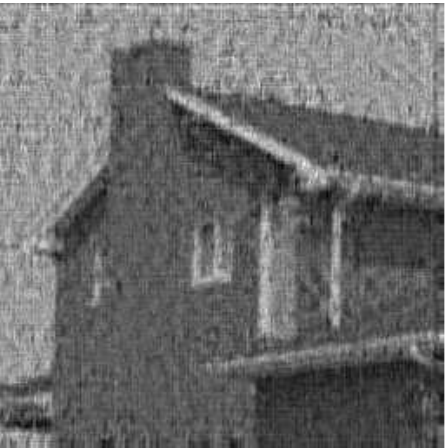}&
\hspace{0.8cm}\includegraphics[height=1.0in]{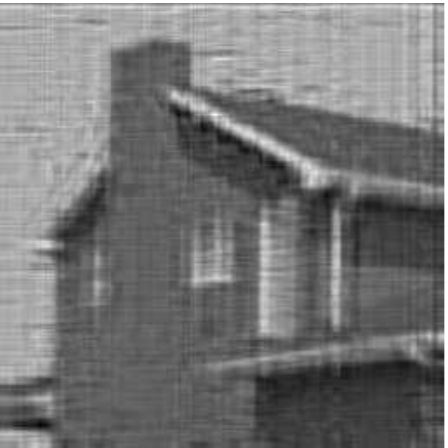}\\
\hspace{-2.0cm}{\small LR-FLSQR} & 
\hspace{0.8cm}{\small FLSQR-NNR(v)} &
\hspace{0.8cm}{\small FLSQR-NNR} & 
\hspace{0.8cm}{\small IRN-LSQR-NNR}
\end{tabular}
\end{center}
\caption{{\em Example 3} (\texttt{house}). Exact and corrupted images; best reconstructions obtained by standard and new solvers.}
\label{fig:recin}
\end{figure}

Figure \ref{fig:relerrin} displays the history of the relative errors for LSQR, LR-FLSQR (with truncation of the basis vectors for the solution, as well as the solution, to rank 20), FLSQR-NNR, FLSQR-NNR(v) and IRN-LSQR-NNR. Figure \ref{fig:recin} displays the exact and corrupted images, together with the best reconstructions obtained by the methods listed above: these correspond to the 16th, 35th, 70th, 30th and 62nd iterations of LSQR, LR-FLSQR, FLSQR-NNR, FLSQR-NNR(v) and IRN-LSQR-NNR, respectively (i.e., these are the iterations where the minimum relative error is attained over the total 100 iterations).

\begin{figure}[H]
\begin{center}
\includegraphics[height=1.8in]{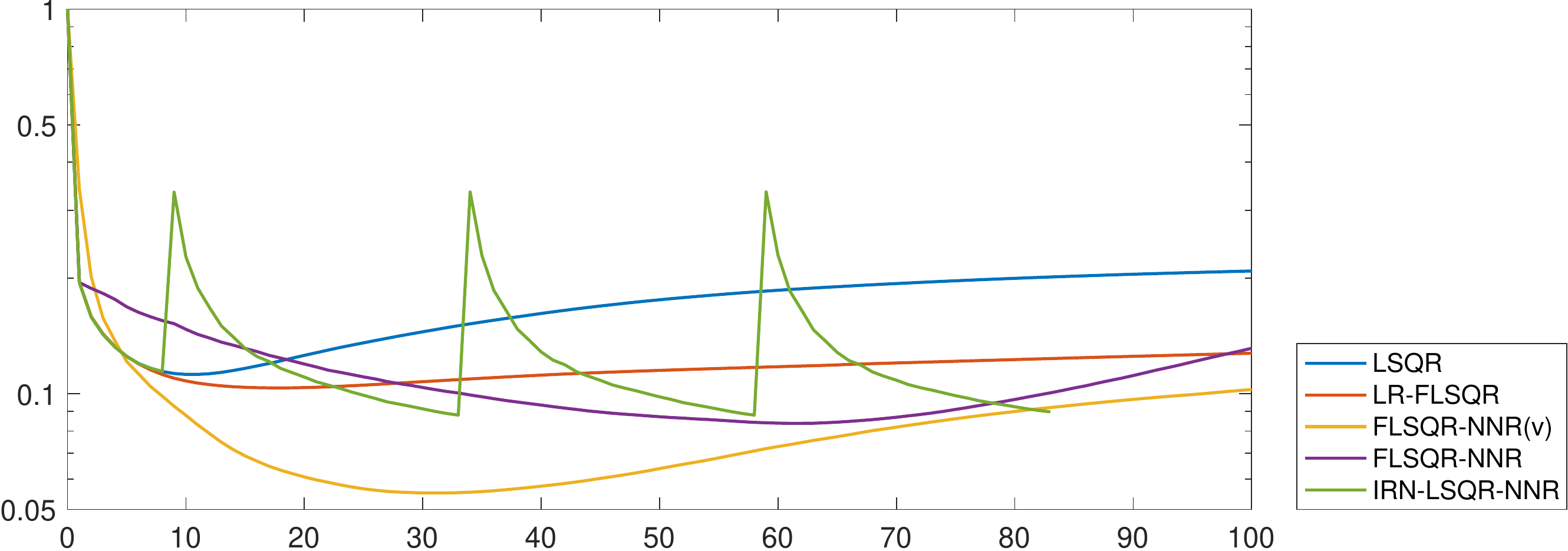}
\end{center}
\caption{{\em Example 3} (\texttt{peppers}). Relative errors vs. number of iterations for different solvers.}

\label{fig:relerrinp}
\end{figure}

\begin{figure}[H]
\begin{center}
\begin{tabular}{ccc}
\hspace{-1.5cm}\includegraphics[height=1.0in]{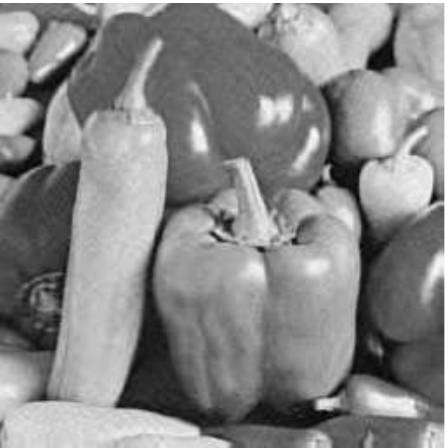}  &
\hspace{0.5cm}\includegraphics[height=1.0in]{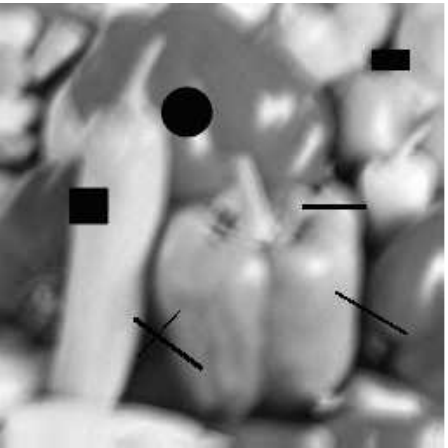} &
\hspace{0.5cm}\includegraphics[height=1.0in]{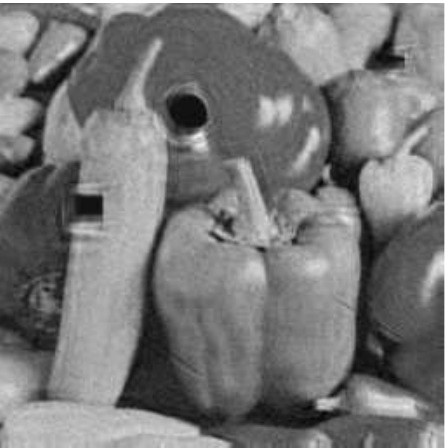} \\
\hspace{-1.5cm}{\small exact} & 
\hspace{0.5cm}{\small corrupted}&
\hspace{0.5cm}{\small LSQR}\vspace{0.3cm}
\end{tabular}
\begin{tabular}{cccc}
\hspace{-2.0cm}\includegraphics[height=1.0in]{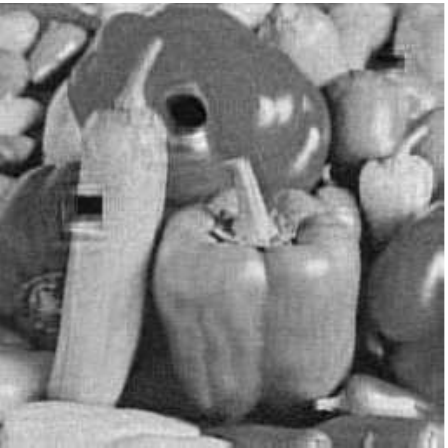} &
\hspace{0.8cm}\includegraphics[height=1.0in]{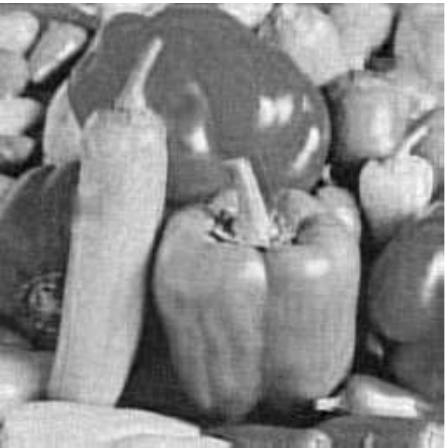}&
\hspace{0.8cm}\includegraphics[height=1.0in]{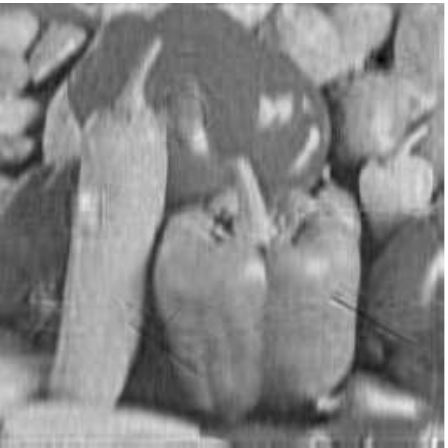}&
\hspace{0.8cm}\includegraphics[height=1.0in]{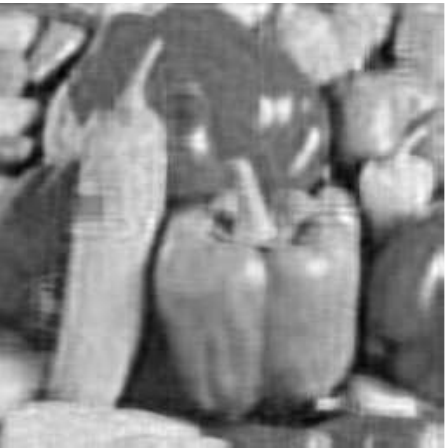}\\
\hspace{-2.0cm}{\small LR-FLSQR} & 
\hspace{0.8cm}{\small FLSQR-NNR(v)} &
\hspace{0.8cm}{\small FLSQR-NNR} & 
\hspace{0.8cm}{\small IRN-LSQR-NNR}
\end{tabular}
\end{center}
\caption{{\em Example 3} (\texttt{peppers}). Exact and corrupted images; best reconstructions obtained by standard and new solvers.}
\label{fig:recin}
\end{figure}

Secondly, we consider a test problem similar to the previous one, i.e., we take an exact image commonly known as \texttt{peppers}, we truncate it to rank 50, and we obtain the forward operator $\mxA$ by first applying a shaking blur, and by then undersampling the blurred image. Here the exact image has a total number of 65536 ($256\times 256$) pixels, and only around 1.3\% of pixels are missing and should be inpainted: differently from the previous problem, the missing pixels follow particular patterns (e.g., circles, squares, and rectangles), and this makes the inpainting task somewhat more challenging. 
Figure \ref{fig:relerrinp} displays the history of the relative errors for LSQR, LR-FLSQR (with truncation of the basis vectors for the solution, as well as the solution, to rank 50), FLSQR-NNR, FLSQR-NNR(v) and IRN-LSQR-NNR. 
Figure \ref{fig:recinp} displays the exact and corrupted images, together with the best reconstructions obtained by the methods listed above: these correspond to the 11th, 18th, 61st, 31st and 58th iterations of LSQR, LR-FLSQR, FLSQR-NNR, FLSQR-NNR(v) and IRN-LSQR-NNR, respectively (i.e., these are the iterations where the minimum relative error is attained over the total 100 iterations).

It is evident that FLSQR-NNR(v) achieves reconstructions of superior quality, including clarity, brightness, and smoothness. Its ability to fill-in missing spots with pixels that are of similar intensity to their surroundings is the best among all methods. The next best reconstructions are computed by IRN-LSQR-NNR for the \texttt{house} test image, and by FLSQR-NNR for the \texttt{pepper} test image: in both cases, these methods are also good at removing noise and restoring missing pixels. However, for both test images, the reconstructions obtained by IRN-LSQR-NNR lack clarity compared to ones obtained by both FLSQR-NNR and FLSQR-NNR(v) methods; compared to the reconstructions obtained by LSQR and LR-FLSQR, they are anyway more desirable in terms of recovered brightness and fill-in of the missing pixels.

%

\paragraph{A study of regularization parameters}

In the previous examples we have seen that the IRN-NNR methods and the flexible Krylov NNR methods perform exceptionally well on image deblurring, tomography, and inpainting problems, producing superior reconstructions compared to existing methods including SVT, RS-LR-GMRES and the low-rank flexible Krylov methods inspired by RS-LR-GMRES, even without the use of additional regularization. In this section, we explore the effect of additional regularization (i.e., we set $\widehat{\lambda} \neq 0$) on the reconstructed images and the corresponding relative errors. In particular, additional regularization allows the new methods to be used in a hybrid fashion. We are going to observe that the new IRN-NNR and flexible Krylov NNR methods are already very powerful, and give reconstructions of superb quality without the use of additional regularization: as a result, there is only little to negligible room for the methods to improve when they are used in a hybrid fashion.

We consider three different ways of choosing the regularization parameter $\hl$. (i) We take the ``secant method'' mentioned in Section \ref{sec:implementation}, which updates the regularization parameter at each iteration using the discrepancy. (ii) We select the optimal regularization parameter which minimizes the 2-norm of the difference between the exact solution and the regularized solution at each iteration. Namely, when using standard GMRES and LSQR, at the $m$th iteration we seek to minimize
\begin{equation}\nonumber
\|\bfx^{\mathrm{ex}} - \bfx_{m,\widehat{\lambda} }\| = \|\mxV_m^T\bfx^{\mathrm{ex}} - \mxV_m^T\bfx_{m,\widehat{\lambda} }\| = \|\mxV_m^T\bfx^{\mathrm{ex}} - \bfy_{m,\widehat{\lambda} }\| \,;
\end{equation}
when using the IRN methods we should incorporate the appropriate preconditioners $(\mxWpg)_k$ and $\mxS_k$ and, for all the iterations in the inner iteration cycle corresponding to the $k$th outer iteration, we seek to minimize
\begin{equation}\nonumber
\|\widehat{\bfx}^{\mathrm{ex}} - \mxV_m\bfy_{m,\widehat{\lambda} }\| = \|\mxV_m^T\widehat{\bfx}^{\mathrm{ex}} - \mxV_m^T\mxV_m\bfy_{m,\widehat{\lambda} }\| = \|\mxV_m^T\widehat{\bfx}^{\mathrm{ex}} - \bfy_{m,\widehat{\lambda} }\|,\ \mbox{where}\ \widehat{\bfx}^{\mathrm{ex}}=(\mxWpg)_k\mxS_k\bfx^{\mathrm{ex}}.
\end{equation}
It is intrinsically difficult to implement this strategy for flexible Krylov subspace methods, because of the complexity of changing preconditioners at each iteration. (iii) We perform a manual exhaustive search. Namely, we first run the solvers multiple times using various regularization parameters $\widehat{\lambda}$, starting with a larger range and narrowing down to a smaller range containing the best parameter; we then record the minimum relative errors among all iterations for all values of $\widehat{\lambda}$, and select the corresponding $\widehat{\lambda}$. This approach is the most expensive, and differs from the previous one in that the (optimal) regularization parameter $\widehat{\lambda}$ is fixed for all iterations. Of course, both the second and third approaches require the knowledge of the exact solution and we test them only to investigate the best possible performance of the hybrid approach.

Table \ref{tab:regchoices} compares the performances (in terms of minimum relative error achieved by each method) of standard Krylov methods (GMRES and LSQR) and their IRN-NNR and flexible NNR (F-NNR) counterparts, with and without using a hybrid approach. In this way we can understand  how the use of additional regularization affects each solver differently. The three parameter choice methods described above are dubbed ``Secant (i)", ``Optimal (ii)" and ``Fixed (iii)", respectively. All the previous examples are considered here. 
GMRES and its counterparts IRN-GMRES-NNR, FGMRES-NNR are used for Example 1, while LSQR and its counterparts IRN-LSQR-NNR, FLSQR-NNR(v) are used for Examples 2 and 3.

\begin{table}[htbp]
\footnotesize
\begin{center}
\begin{tabular}{|l|l|c|c|c|c|c|c|c|c|}
\cline{3-10}
\multicolumn{2}{c|}{} & $\hl=0$ & $\hl\neq0$ & $\hl=0$ & $\hl\neq0$ & $\hl=0$ & $\hl\neq0$ & $\hl=0$ & $\hl\neq0$\\
\cline{3-10}
\multicolumn{2}{c|}{} & \multicolumn{2}{c|}{\emph{Example 1}} & \multicolumn{2}{c|}{\emph{Example 2}} & \multicolumn{2}{c|}{\emph{Example 3} (\texttt{house})} & \multicolumn{2}{c|}{\emph{Example 3} (\texttt{peppers})}\\
    \hline
    \hline
    \multirow{3}{*}{Standard}& Secant (i) &0.2995 & 0.2528 &
    0.1201 & 0.1389 &
    0.2710 & 0.2713 &
    0.1123 & 0.1120\\ 
    \cline{2-10}
    & Optimal (ii) &0.2995 & 0.2268 &
    0.1201 & 0.1201&
    0.2710 & 0.2708&
    0.1123 & 0.1120\\
     \cline{2-10}
    & Fixed (iii) &0.2995 & 0.2268 &
    0.1201 & 0.1183 &
    0.2710 & 0.2708 &
    0.1123 & 0.1120\\
    \hline
    \hline
    \multirow{3}{*}{IRN-NNR}& Secant (i) &0.2081 & 0.2096 & 
    0.0685 & 0.0696 &
    0.1234 & 0.1234 &
    0.0880 & 0.0871\\ 
     \cline{2-10}
    & Optimal (ii) &0.2081 & 0.2292 &
    0.0685 & 0.0685 &
    0.1234 & 0.1234 &
    0.0880 & 0.0880\\
         \cline{2-10}
    & Fixed (iii) &0.2081 & X &
    0.0685 & 0.0660 &
    0.1234 & X &
    0.0880 & 0.0866\\
    \hline
    \hline
     \multirow{2}{*}{F-NNR}& Secant (i) &0.2829 & 0.2658 &
     0.0577 & 0.0684 &
     0.1021 & 0.1034 &
     0.0552 & 0.0549\\ 
     \cline{2-10}
    & Fixed (iii) &0.2829 & 0.2640 &
    0.0577 & 0.0568 &
    0.1021 & X &
    0.0552 & 0.0548\\
    \hline
    \hline
\end{tabular}
\end{center}
\caption{Minimum relative errors without ($\hl=0$) and with ($\hl\neq 0$) a hybrid approach. The mark ``X" means that the optimal regularization parameter found by the ``Fixed (iii)" method is $10^{-16}$, hence there is no need for regularization.}
\label{tab:regchoices}
\end{table}

\begin{figure}[htbp]
\begin{center}
\begin{tabular}{cccc}
\hspace{-2.0cm}\includegraphics[height=1in]{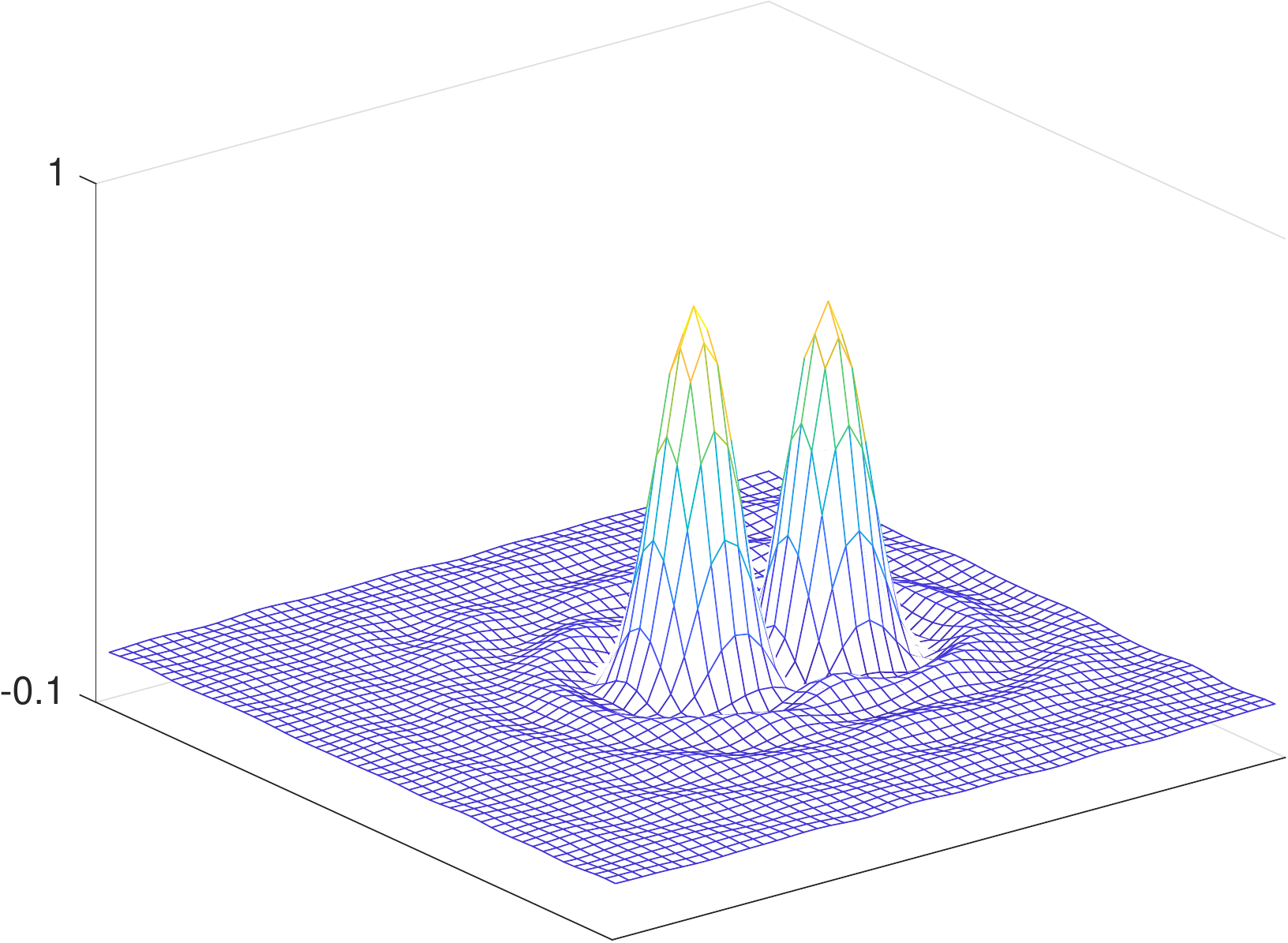}  &
\hspace{0.2cm}\includegraphics[height=1in]{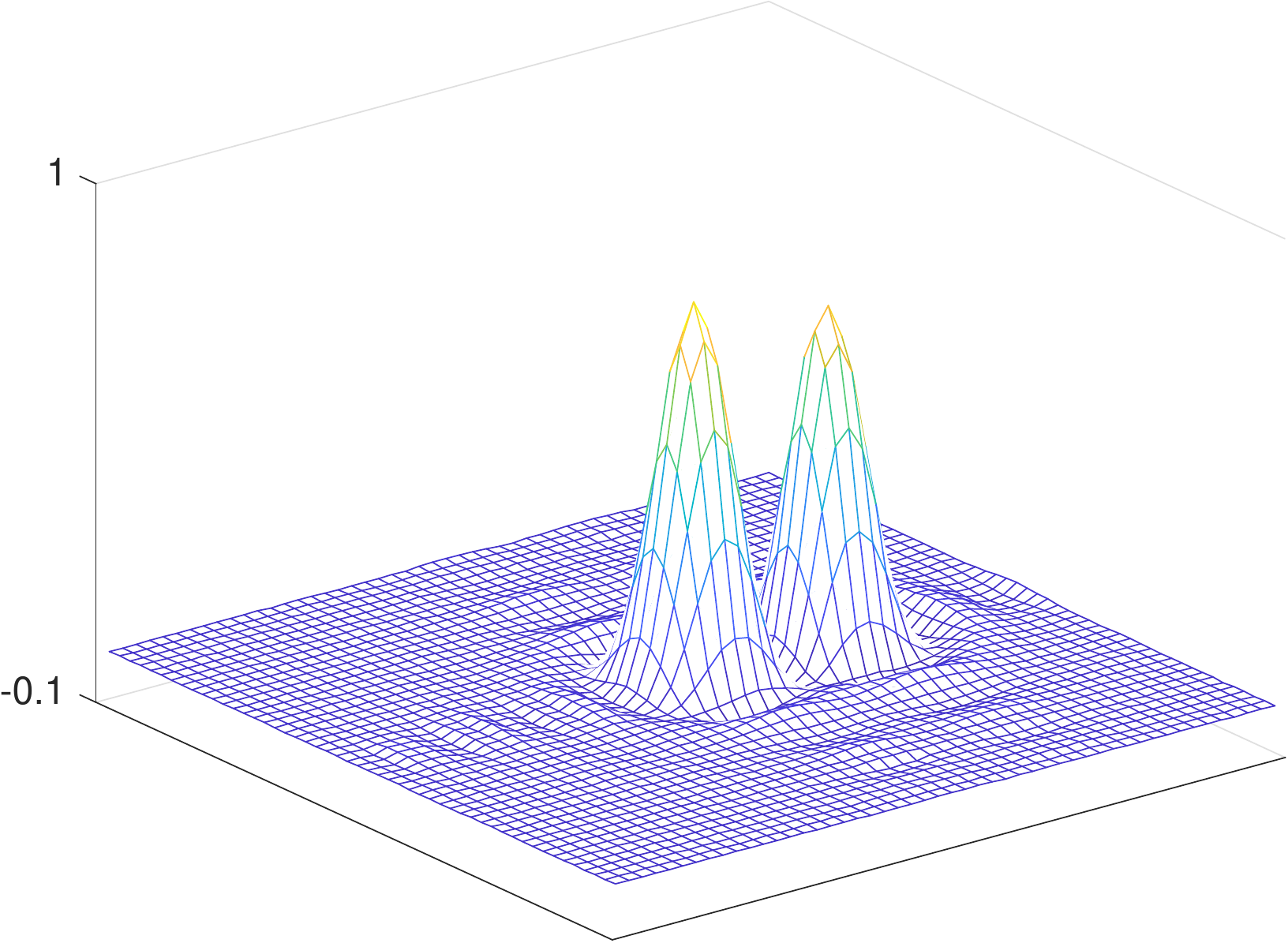}&
\hspace{0.8cm}\includegraphics[height=1in]{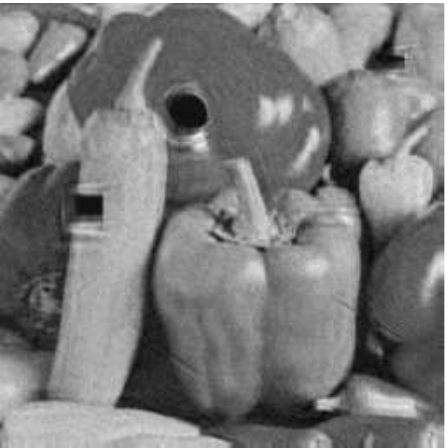}&
\hspace{0.8cm}\includegraphics[height=1in]{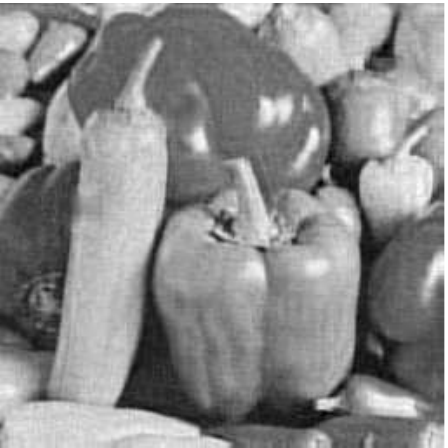} \\
\hspace{-2.0cm}{\small hybrid GMRES} & 
\hspace{0.2cm}{\small IRN-GMRES-NNR} & 
\hspace{0.8cm}{\small hybrid LSQR} & 
\hspace{0.8cm}{\small FLSQR-NNR(v)}\\
\end{tabular}
\end{center}
\caption{Reconstructions obtained by standard hybrid Krylov methods and by the new methods without using additional regularization. Left side: zoomed in surface plots of the reconstructions of \emph{Example 1}; right side: reconstructions of \emph{Example 3} (\texttt{peppers}).}
\label{fig:reg_compare}
\end{figure}

It is easy to observe that the use of additional regularization is most effective for the standard GMRES solver, where the minimum relative error is reduced significantly. 
However, for the other solvers, the hybrid approach does not have a notable advantage over not using regularization. At times the ``Fixed (iii)" parameter choice strategy delivers a regularization parameter of the order of $10^{-16}$, which is numerically equivalent to not having regularization. 
This indicates that our new IRN-NNR and F-NNR methods are successful in computing good reconstructions and, even without additional regularization, they perform much better than standard Krylov methods used in a hybrid fashion (comparing IRN-GMRES-NNR to GMRES in Example 1, and FLSQR-NNR(v) to LSQR in the other examples). Figure \ref{fig:reg_compare} shows a couple of such comparisons.

\section{Conclusions}\label{sec:conclusions}
This paper introduced new solvers, based on Krylov subspace methods, for the computation of approximate low-rank solutions to large-scale linear systems of equations. Our main goal was to apply the new methods to regularize inverse problems arising in imaging applications. The starting point of our derivations was an IRN approach to the NNR$p$ problem (\ref{eq:nnrp}). In this way, the original problem (\ref{eq:nnrp}) is reduced to the solution of a sequence of quadratic problems, where an appropriate regularized linear transformation is introduced to approximate the nondifferentiable nuclear norm regularization term. Our new methods make smart use of Kronecker product properties to reformulate each quadratic problem in the IRN sequence as a Tikhonov-regularized problem in standard form. 
We use both Krylov methods with fixed ``preconditioners'' within an inner-outer iteration scheme (namely, IRN-LSQR-NNR$p$ and IRN-GMRES-NNR$p$), and Krylov methods with flexible iteration-dependent ``preconditioners'' within a single iteration scheme (namely, FLSQR-NNR$p$, FGMRES-NNR$p$, LR-FGMRES, and LR-FLSQR). Some of these methods (namely, IRN-LSQR-NNR$p$, IRN-GMRES-NNR$p$, FLSQR-NNR$p$, and FGMRES-NNR$p$) can be used in a hybrid framework, so that the Tikhonov regularization parameter can be efficiently, effectively, and adaptively chosen. These new solvers are shown to perform exceptionally well on the test problems described in Section \ref{sec:experiments}, and they give reconstructions of significantly improved quality over existing methods.

Future work includes the extension of the present methods to handle cases where the solution of (\ref{eq:Axeqb}) is low-rank but rectangular, i.e., $\vect^{-1}(\bfx)=\mxX\in\RR^{m\times n}$ with $m\neq n$. Also, while a solid theoretical justification is provided for IRN-LSQR-NNR$p$ and IRN-GMRES-NNR$p$, the same is not true for FGMRES-NNR$p$ and FLSQR-NNR$p$: further analysis will be needed to deeply understand the regularization properties of these flexible solvers. Finally, the new IRN-LSQR-NNR$p$ and IRN-GMRES-NNR$p$ methods can be reformulated to work with well-posed problems and in the framework of matrix equations, possibly providing a valid and principled alternative to the current popular methods based on low-rank-projected and restarted Krylov solvers. 

\bibliographystyle{siamplain}
\bibliography{references}
\end{document}